\documentclass[10pt,a4paper]{article}

\usepackage[hidelinks]{hyperref}
\usepackage{lineno}

\usepackage{amssymb, bm}
\usepackage{amsmath}
\usepackage{amsthm}
\usepackage{mathrsfs}
\usepackage{siunitx}
\usepackage[margin=2cm]{geometry}
\usepackage{subfig}
\usepackage{float}
\usepackage{xcolor}
\usepackage{subfloat}

\usepackage{mathtools}
\usepackage{bbm}
\usepackage{multirow}
\usepackage{stmaryrd}

\usepackage{algorithm}     
\usepackage{algpseudocode} 

\usepackage[utf8]{inputenc} 
\usepackage[T1]{fontenc}

\usepackage{tikz}
\usetikzlibrary{decorations.markings}
\usetikzlibrary{shapes.geometric}
\usetikzlibrary{shapes.arrows}
\usetikzlibrary{fit,calc}
\usetikzlibrary{quotes,angles}
\usetikzlibrary{arrows.meta}
\usepackage{animate}  
\usepackage{graphicx}
\usepackage{pgf}

\usepackage{pgfplots}
\usepackage[normalem]{ulem}

\modulolinenumbers[5]

\newtheorem{remark}{Remark}
\newtheorem{lemma}{Lemma}

\newcommand{\vertiii}[1]{{\left\vert\kern-0.25ex\left\vert\kern-0.25ex\left\vert #1 
   \right\vert\kern-0.25ex\right\vert\kern-0.25ex\right\vert}}

\usepackage{geometry}
\makeatletter
\def\blfootnote{\xdef\@thefnmark{$\star$}\@footnotetext}
\makeatother
\newenvironment{Authors}%
  {\begin{center}\begin{bfseries}}%
  {\end{bfseries}\end{center}}
\newenvironment{Addresses}%
  {\begin{flushleft}\begin{itshape}}%
  {\end{itshape}\end{flushleft}}
  \newcommand{\email}[1]{\hspace*{\stretch{1}}\emph{\texttt{#1}}}

 \usepackage{fancyhdr}  
  \fancypagestyle{plain}{
\fancyhead{}
\fancyhead[C]{\hfill Submitted to Elsevier, December 2024}
 }

\begin{document}

\title{Optimization-based model order reduction of fluid-structure interaction problems} 
 \date{}
 \maketitle
\vspace{-50pt} 
 
\begin{Authors}
Tommaso Taddei$^{1}$,  Xuejun Xu$^{2,3}$,
Lei Zhang$^2$.
\end{Authors}

\begin{Addresses}
$^1$
Univ. Bordeaux, CNRS, Bordeaux INP, IMB, UMR 5251, F-33400 Talence, France\\ Inria Bordeaux Sud-Ouest, Team MEMPHIS, 33400 Talence, France, \email{tommaso.taddei@inria.fr} \\[3mm]
$^2$
School of Mathematical Sciences, Tongji University, Shanghai 200092 , China, \email{22210@tongji.edu.cn} \\[3mm]
$^3$
Institute of Computational Mathematics, AMSS, Chinese Academy of Sciences, Beijing 100190, China, \email{xxj@lsec.cc.ac.cn} \\
\end{Addresses}

\begin{abstract}
We introduce  optimization-based 
full-order and reduced-order formulations of  fluid structure interaction problems.
We study the flow of an incompressible Newtonian fluid which interacts with an elastic    body:
we consider an   arbitrary Lagrangian Eulerian formulation of the fluid problem and a fully Lagrangian formulation of the solid problem;
we rely on a finite element discretization of both fluid and solid equations.
The distinctive feature of our approach is an implicit coupling of fluid and structural problems that relies on the solution to a constrained optimization problem with equality constraints.
We discuss the application of projection-based model reduction to both fluid and solid subproblems: we rely on Galerkin projection for the solid equations and on least-square Petrov-Galerkin projection for the fluid equations.
Numerical results for three model problems illustrate the many features of the formulation.
\end{abstract}

\medskip

\section{Introduction}
\label{sec:introduction}

\subsection{Approximation of fluid structure interaction problems}

Accurate and reliable simulations of fluid structure interaction (FSI) problems is of utmost importance in science and engineering:
applications range from
aeroelasticity
\cite{cinquegrana2021,farhat1998load,Farhat_Geuzaine_2004},
bridge aeroelasticity
\cite{Sangalli_Braun_2020}, to the simulation of 
blood flow  in arteries  
\cite{Barker2010,formaggia2010cardiovascular,Wu_Cai_2014}, sailing boats  \cite{Lombardi2012}, and  wave energy conversion devices \cite{Agamloh2008}.
Model order reduction (MOR) methods aim to reduce the marginal cost (i.e., the cost in the limit of many queries) associated with the solution to parametric systems, for many-query and real-time applications: despite the many recent contributions to the field, significant advancements are still needed in the development of MOR methods for FSI problems. The aim of this work is to develop a new formulation for FSI problems, with emphasis on the coupling between 
full-order models (FOMs) and reduced-order models (ROMs).

The formulation of 
FSI problems comprises a system of balance laws for the fluid domain, an equation of motion for the solid, and a set of coupling (or interface) conditions that ensure the continuity of displacement, velocity, and normal stresses.
FSI methods can be classified as 
body-fitted (or interface tracking) methods, 
interface capturing methods, and fully-Eulerian methods.
Body-fitted methods are based either on 
a fully-Lagrangian formulation 
\cite{cremonesi2020state} or 
on 
an arbitrary Lagrangian Eulerian (ALE) formulation 
\cite{donea2004arbitrary}:
they  ensure a precise representation of the geometry and enable a straightforward imposition of the interface conditions.
Interface capturing techniques comprise the immersed boundary method \cite{kim2019immersed},
the embedded boundary method \cite{farhat2014ale},
the fictitious boundary method \cite{van2004combined}, and the 
cut-cell method \cite{schott2019monolithic}:
they are   based on a combination of Lagrangian and Eulerian formalism;
as opposed to the body-fitted methods, 
they avoid the need for frequent remeshing and can readily handle topological variations (splitting or fusion) and surface contacts (collisions) at the price of a more involved enforcement of transmission conditions.
Finally, 
fully-Eulerian methods
\cite{bergmann2022eulerian,cottet2008eulerian,richter2013fully,sugiyama2011full,valkov2015eulerian}
consider the fluid and the structure as a multiphase viscous material governed by the same equations with a space- and time-dependent constitutive law;
interface conditions  are  enforced weakly in the formulation.
Eulerian methods   are potentially less accurate than the previous two classes of methods, particularly for thin structures; however, their use  enables a more transparent handling of the mesh and thus a
drastic reduction of the computational setup overhead.

Following taxonomy introduced in \cite{fernandez2011coupling}, we distinguish between implicit, semi-implicit, and explicit coupling schemes for FSI problems.
Explicit coupling methods rely on the explicit treatment of one or two interface conditions:
they are simple to implement and relatively inexpensive to evaluate, but they might be unstable under certain choices of the physical parameters; in particular, they are ill-suited to simulate  FSI problems for which densities of the fluid and the solid are comparable
\cite{causin2005added}.
For this reason, 
they are broadly used in aeroelasticity (e.g., 
\cite{farhat2006provably}), but they are not appropriate, e.g., for  blood-flow simulations.
Conversely, implicit schemes
\cite{Deparis2006_steklov,fernandez2005newton} ensure energy balance and stability, at the price of a larger cost per iteration.
Finally, semi-implicit coupling methods
(see, e.g., \cite[section 4]{fernandez2011coupling}) seek a compromise between stability and effectiveness for select applications.
In this work, we propose an implicit body-fitted ALE formulation for FSI; our formulation enables
a seamless coupling of full-order and reduced-order models through the vehicle of an optimization formulation.

\subsection{Optimization-based model reduction of fluid structure interaction problems}

We consider the optimization-based method proposed by Gunzburger and coauthors
\cite{Gunzburger_Lee_2000,Gunzburger_Peterson_Kwon_1999}, which was later extended to FSI problems in 
\cite{Kuberry_Lee_2013,Kuberry_Lee_2015,Kuberry_Lee_2016};
this work extends our previous contribution \cite{Taddei2024optimization} on 
model reduction of parametric incompressible flows in parametric geometries.
The formulation reads as an optimal control problem where the control is given by the flux at the fluid-structure interface and the state variables are the velocity and the pressure in the fluid, and the displacement in the solid. The formulation automatically ensures continuity of the normal stresses and the displacement; on the other hand, the  continuity of the velocity is weakly ensured through the objective function.

Following \cite{Taddei2024optimization},
we solve the minimization problem using 
sequential quadratic programming (SQP, \cite{bonnans2006numerical}), which guarantees rapid convergence without resorting to Lagrange multipliers. 
We employ separate low-dimensional approximation spaces for the fluid state, the solid state, and the control using proper orthogonal decomposition
(POD,  \cite{volkwein2011model});
we consider a Galerkin ROM for the structure and a least-square Petrov-Galerkin 
(LSPG, \cite{carlberg2013gnat,taddei2021space}) ROM for the fluid subproblem. Furthermore, we exploit an enrichment strategy for the state spaces that ensure stability of the coupled problem.
Finally, we investigate the coupling of a Galerkin ROM for the structure with a full-order model for the fluid.

Our method is also related to recent works in component-based MOR that considered similar optimal control formulations. First, we recall the 
one-shot overlapping Schwartz method 
 \cite{iollo2023one}, which also relies on an optimal control formulation. Second, we recall the optimization-based method
 \cite{d2015optimization}
 for the coupling of non-local and local diffusion models and the domain decomposition LSPG ROM in
 \cite{hoang2021domain}, which considers a minimum residual formulation with constraints given by the interface conditions.
 Finally, we recall the work by Prusak
 \cite{prusak2023optimisation}, which is also based on the minimization formulation of 
 \cite{Gunzburger_Lee_2000} and was recently extended to unsteady problems \cite{prusak2024optimisation}.
The minimization framework considered in this work avoids the need for compatibility conditions among the approximation spaces  at the fluid-structure interface: in this respect, it is related to recent efforts to the development of overlapping Schwartz methods for coupled problems
\cite{wentland2024role}.

For completeness, we also list select works on model reduction of FSI problems.
Early works on POD-Galerkin ROMs for linearized FSI problems appeared in
\cite{lieu2006reduced,kalashnikova2013stable};
Xiao and collaborators \cite{xiao2016non} proposed an alternating  technique for FSI problems that exploit local ROMs based on radial basis function interpolation;
finally, the authors of 
\cite{Ballarin2017_fsi,nonino2023projection} 
exploit Galerkin ROMs to simulate the interaction of a parametric 
incompressible Newtonian flow and an elastic body through a body-fitted  semi-implicit  partitioned formulation.
 
The outline of the paper is as follows.
In section \ref{sec:governing_equations}, we review the ALE formulation of FSI problems and we introduce the constitutive laws considered in the numerical experiments.
In section \ref{sec:discrete_equations}, we discuss the finite element (FE) approximation of the fluid and the solid equations, and we introduce the optimal control formulation.
Then, in section \ref{sec:ROM}, we present the model reduction formulation: we discuss the construction of the reduced spaces, the definition of the local ROMs, and the solution to the global problem.
In section \ref{sec:numerical_results}, we present extensive numerical investigations for three two-dimensional model problems.
Section \ref{sec:conclusions} completes the paper.
 
\section{Governing equations}
\label{sec:governing_equations}

\subsection{Strong form of the equations}

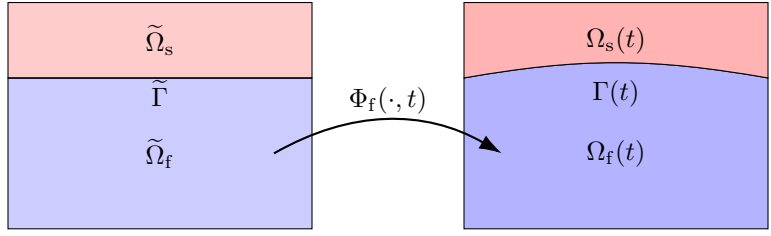
\begin{figure}[H]
    \centering
\begin{tikzpicture}
    \draw [fill=blue!20] (0,0) rectangle (4,2);
    \node at (2,1) {$\widetilde{\Omega}_{\rm f}$};

    \draw [fill=red!20] (0,2) rectangle (4,3);
    \node at (2,2.5) {$\widetilde{\Omega}_{\rm s}$};

    \draw (0,2) -- (4,2);
    \node at (2,1.8) {$\widetilde{\Gamma}$};

    \draw [fill=blue!30] (6,0) -- (10,0) -- (10,2) to [out=170,in=10] (6,2) -- cycle;
    \node at (8,1) {$\Omega_{\rm f}(t)$};

    \draw [fill=red!30] (6,2) to [out=10,in=170] (10,2) -- (10,3) -- (6,3) -- cycle;
    \node at (8,2.5) {$\Omega_{\rm s}(t)$};

    \node at (8,1.8) {$\Gamma(t)$};

    \draw[-{Latex[length=3mm, width=2mm]}, thick] (3.5,1) to [out=30,in=150] (6.5,1);
    	\node at (5.0,1.7) {${\Phi}_{\rm f}(\cdot, t)$};
\end{tikzpicture}

\caption{referential  and current   domains of FSI problems. Here, $\Phi_{\rm f}(t)$ is the ALE bijection that maps the reference fluid domain $\widetilde{\Omega}_{\rm f}$ into the current domain $\Omega_{\rm f}(t)$.}
\label{fig:FSI_illustration}
\end{figure}

We consider an ALE formulation of the FSI problem
\cite{Barker2010,donea2004arbitrary,formaggia2010cardiovascular,Wu_Cai_2014}. We denote by $\Omega_{\rm f}(t)$ and $\Omega_{\rm s}(t)$ the fluid and solid domains at time $t$ and by $\Gamma(t)$ the fluid/solid interface. We denote by $d_{\rm s}: \widetilde{\Omega}_{\rm s} \times \mathbb{R}_+ \to \mathbb{R}^d$ the solid displacement field such that
$\Omega_{\rm s}(t) = \{\widetilde{x}+ d_{\rm s}(\widetilde{x}, t)   : \widetilde{x} \in \widetilde{\Omega}_{\rm s} \}$, by $\Phi_{\rm f}: \widetilde{\Omega}_{\rm f} \times \mathbb{R}_+ \to \Omega_{\rm f}(t)$ the ALE deformation map such that 
$\Omega_{\rm f}(t) = \Phi_{\rm f}(\widetilde{\Omega}_{\rm f}, t )$ and 
$\Gamma(t) = \Phi_{\rm f}(\widetilde{\Gamma}, t )$ where $\widetilde{\Gamma}= \partial \widetilde{\Omega}_{\rm f} \cap \partial \widetilde{\Omega}_{\rm s}$. 
Figure \ref{fig:FSI_illustration} illustrates the fluid and solid, reference and Eulerian domains.
We further define the fluid velocity and pressure 
$u_{\rm f}: \Omega_{\rm f}(t)  \times \mathbb{R}_+ \to \mathbb{R}^{\texttt{d}}$ (with $\texttt{d}=2$ or $\texttt{d}=3$) and 
$p_{\rm f}: \Omega_{\rm f}(t)  \times \mathbb{R}_+ \to \mathbb{R}$, 
the ALE velocity $\widetilde{\omega}_{\rm f}:=\frac{\partial \Phi_{\rm f}}{\partial t}$, and the ALE time derivative for the Eulerian field $q: \Omega_{\rm f}(t) \times \mathbb{R}_+ \to \mathbb{R}$:
\begin{equation}
\label{eq:ALE_time_derivative}
\frac{\partial q}{\partial t}\big|_{\Phi_{\rm f}}(x,t):=
\frac{d}{dt} \left(
q\left(\Phi_{\rm f}(\widetilde{x},t), t \right)
\right),
\quad
{\rm with}
\;\;
x=\Phi_{\rm f}(\widetilde{x},t).
\end{equation}
Note that 
\eqref{eq:ALE_time_derivative} measures the rate of change of the field $q$ in a point that moves 
with the computational domain.

Exploiting the previous definitions, we introduce the fluid equations in the ALE frame:
\begin{subequations}
\label{eq:strong_equations}
\begin{equation}
\label{eq:fsi_ns}
\left\{
\begin{array}{ll}
\displaystyle{
\rho_{\rm f}\frac{\partial {u}_{\rm f}}{\partial t}  \Big|_{\Phi_{\rm f}} + \rho_{\rm f}(({u}_{\rm f}-\omega_{\rm f})\cdot \nabla){u}_{\rm f}  -\nabla\cdot {\sigma}_{\rm f} =0, }&  \mbox{in}\; \Omega_{\rm f}(t), \\[2mm]
\nabla \cdot {u}_{\rm f} = 0, & \mbox{in}\; \Omega_{\rm f}(t), \\[2mm]
{u}_{\rm f}|_{\Gamma^{\rm dir}_{\rm f}}={u}_{\rm f}^{\rm dir},
\;\;
{\sigma}_{\rm f} {n}_{\rm f}
|_{\Gamma^{\rm neu}_{\rm f}}={g}^{\rm neu}_{\rm f},
&
\end{array}
\right.
\end{equation} 
where
$\rho_{\rm f}$ is the fluid density, 
${n}_{\rm f}$ is the outward normal to $\Omega_{\rm f}(t)$,
${\omega}_{\rm f}=\widetilde{\omega}_{\rm f}\circ \Phi_{\rm f}^{-1}$ is the ALE velocity in the ALE frame,
$\sigma_{\rm f}$ is the stress tensor that is a function of 
fluid velocity and pressure, 
and 
${u}_{\rm f}^{\rm dir}$, ${g}^{\rm neu}_{\rm f}$ are suitable boundary data. Then, we introduce the solid equations in the Lagrangian frame:
\begin{equation}
\left\{
\begin{array}{ll}
\displaystyle{
\rho_{\rm s}\frac{\partial^2 {{d}}_{\rm s}}{\partial t^2}  -{\nabla} \cdot  {\sigma}_{\rm s} = f_{\rm s},} &\quad \mbox{in}\; \widetilde{\Omega}_{\rm s}, \\[3mm]
{d}_{\rm s}|_{\widetilde{\Gamma}^{\rm dir}_{\rm s}}={{d}}_{\rm s}^{\rm dir},
\;\;
\widetilde{{\sigma}}_{\rm s} \widetilde{{n}}_{\rm s}
|_{\widetilde{\Gamma}^{\rm neu}_{\rm s}}={{g}}_{\rm s}^{\rm neu},
&
\\
\end{array}
\right.
\label{eq:fsi_solid}
\end{equation}
where $\rho_{\rm s}$ is the solid density, 
$\sigma_{\rm s}$ is the stress tensor,
$\widetilde{{n}}_{\rm s}$ is the outward normal to $\widetilde{\Omega}_{\rm s}$,
$f_{\rm s}$ is the source term, and 
 ${{d}}_{\rm s}^{\rm dir}$,
 ${{g}}_{\rm s}^{\rm neu}$ are suitable boundary data.
Finally, we define the coupling conditions:
\begin{equation}
\label{eq:fsi_coupling}
\left\{
\begin{array}{l}
\displaystyle{
\Phi_{\rm f}  = \texttt{id} + {{d}}_{\rm s} }
  \\[2mm]
\displaystyle{
u_{\rm f} \circ \Phi_{\rm f} 
=
\frac{\partial {{d}}_{\rm s}}{\partial t}
}
\\[2mm]
\displaystyle{
{\sigma}_{\rm s} \widetilde{n}_{\rm s} +  J_{\Phi} \left( \sigma_{\rm f} \circ \Phi_{\rm f}  \right) \nabla \Phi_{\rm f}^{-\top} \widetilde{n}_{\rm f}
= 0
}
\\
\end{array}
\right.
\mbox{on}\; \widetilde{\Gamma},
\end{equation}
where
$\texttt{id}$ is the identity map and 
$J_{\Phi}  = {\rm det} (\nabla \Phi_{\rm f})$.
Our goal is to devise a numerical strategy for the solution 
$(u_{\rm f}, p_{\rm f}, d_{\rm s})$
to the coupled system \eqref{eq:fsi_ns}-\eqref{eq:fsi_solid}-\eqref{eq:fsi_coupling} for a proper choice of the constitutive models for $\sigma_{\rm s}$ and $\sigma_{\rm f}$.
\end{subequations}

In this work, we consider an incompressible Newtonian model for the fluid:
\begin{equation}
\label{eq:newtonian_model}
\sigma_{\rm f} := 2\mu_{\rm f} \varepsilon_{\rm f} 
- p_{\rm f} \mathbbm{1},
\end{equation}
where $\mu_{\rm f}$ is the dynamic viscosity and 
$\varepsilon_{\rm f}= \frac{1}{2} (\nabla u_{\rm f} + \nabla u_{\rm f}^\top)$ is the strain tensor in the ALE frame. On the other hand, we consider two distinct constitutive laws for the solid 
\begin{subequations} 
\label{eq:solid_model}
\begin{equation}
\sigma_{\rm s}
=
\left\{
\begin{array}{ll}
2\mu_{\rm s} \varepsilon_{\rm s} + \lambda_{\rm s} (\nabla \cdot d_{\rm s}) \mathbbm{1}   & {\rm ``linear \; \; model''} \\[2mm]
 F_{\rm s} \Sigma_{\rm s}    &   {\rm ``Saint \, Venant \;\; Kirchhoff \;\; model''} \\
\end{array}
\right.  
\end{equation}
where $\lambda_{\rm s}, \mu_{\rm s}$ are the first and  the second Lam{\'e} coefficients, $\varepsilon_{\rm s}= \frac{1}{2} (\nabla d_{\rm s} + \nabla d_{\rm s}^\top)$ is the strain tensor in the Lagrangian frame, 
$F_{\rm s} = \mathbbm{1} + \nabla d_{\rm s}$ is the deformation gradient, 
$\Sigma_{\rm s} =  2\mu_{\rm s} E_{\rm s} + \lambda_{\rm s} {\rm tr}(E_{\rm s}) \mathbbm{1}$ is the second Piola-Kirchhoff stress tensor and $E_{\rm s} = \frac{1}{2}(F_{\rm s}^\top  F_{\rm s}  - \mathbbm{1})$ is the Green-Lagrange tensor. 
In the numerical examples, we consider two-dimensional problems under the plane-strain assumption, that is we consider
\begin{equation}
\lambda_{\rm s}
= \frac{E \nu}{(1+\nu)(1-2\nu)},
\quad
\mu_{\rm s}
= \frac{E }{2(1+\nu)},    
\end{equation}
where $E$ is the Young's modulus and $\nu$ is the Poisson's ratio.
\end{subequations}

Following \cite{Shamanskiy2021_mesh,Wick2011}, we define the ALE map
$\Phi_{\rm f}$ using a pseudo-elasticity model.
In more detail, we define $\Phi_{\rm f} = \texttt{id}+d_{\rm f}$ where $d_{\rm f}:\widetilde{\Omega}_{\rm f} \times \mathbb{R}_+ \to \mathbb{R}^{\texttt{d}}$ satisfies:
\begin{equation}
\label{eq:fsi_mesh_pseudo_elasticity}
\left\{
\begin{array}{ll}  
\displaystyle{-{\nabla} \cdot \left(
2\mu_{\rm m} \varepsilon_{\rm m} + \lambda_{\rm m} (\nabla \cdot d_{\rm f}) \mathbbm{1} 
\right) = 0}, & \mbox{in}\; \widetilde{\Omega}_{\rm f}, \\
\widetilde{d}_{\rm f} =  {d}_{\rm s},  & \mbox{on } \widetilde{\Gamma}, \\
\widetilde{d}_{\rm f} = 0,  & \mbox{on } \partial \widetilde{\Omega}_{\rm f} \setminus \widetilde{\Gamma},
\end{array}
\right.
{\rm where} \;\;
\varepsilon_{\rm m}= \frac{1}{2} (\nabla d_{\rm f} + \nabla d_{\rm f}^\top).
\end{equation}
Here, the  pseudo-material parameters $\mu_{\rm m}$ and $\lambda_{\rm m}$ should be tuned to avoid mesh deterioration; we discuss their choice in  section \ref{sec:numerical_results}.
We remark that 
the model \eqref{eq:fsi_mesh_pseudo_elasticity} is appropriate for moderate deformations, whereas 
nonlinear elasticity models should be employed
for large deformations  \cite{froehle2015nonlinear}. 

\subsection{Weak formulation}
We derive a weak formulation that explicitly includes the control variable that will be used in the MOR formulation. We introduce the spaces:
\begin{equation}
\label{eq:fluid_spaces_continuous}
V_{\rm f}^{\rm dir}(t):=\{v\in [H^1(\Omega_{\rm f}(t)
)]^{\texttt{d}}  \, : \, v\big|_{\Gamma_{\rm f}^{\rm dir}(t)} = u_{\rm f}^{\rm dir}(t)  \},
\quad
V_{{\rm f},0}(t):= [H_{0,\Gamma_{\rm f}^{\rm dir}(t)}^1(\Omega_{\rm f}(t)
)]^{\texttt{d}} ,
\quad
Q_{\rm f}(t):= L^2  (\Omega_{\rm f}(t)
),
\end{equation}
and
\begin{equation}
\label{eq:solid_spaces_continuous}
V_{\rm s}^{\rm dir}(t):=\{v\in [H^1(\widetilde{\Omega}_{\rm s})]^{\texttt{d}}  \, : \, v\big|_{\widetilde{\Gamma}_{\rm s}^{\rm dir}} = d_{\rm s}^{\rm dir}(t) \},
\quad
V_{{\rm s},0} := [H_{0,\widetilde{\Gamma}_{\rm s}^{\rm dir}(t)}^1(\widetilde{\Omega}_{\rm s})]^{\texttt{d}}.
\end{equation}
Finally, we introduce the space $G:= 
[L^2  (\widetilde{\Gamma})]^{\texttt{d}}$ for the control.

We have now the elements to introduce the weak formulation --- below, we omit dependence on the time $t$ to shorten notation. For any $t>0$, we seek $u_{\rm f} \in V_{\rm f}^{\rm dir}$, $p_{\rm f} \in Q_{\rm f}$,  
$d_{\rm s} \in V_{\rm s}^{\rm dir}$ and
$g \in G$ such that
\begin{subequations}
\label{eq:weak_formulation}
\begin{equation}
\label{eq:weak_formulation_a}
\left\{
\begin{array}{ll}
\displaystyle{
\int_{ \Omega_{\rm f}   }
\rho_{\rm f} 
\frac{\partial u_{\rm f} }{\partial t}
\Big|_{\Phi_{\rm f}} \cdot v \, dx
+
R_{\rm f}(u_{\rm f}, p_{\rm f}, \omega_{\rm f}, v)
+
E_{\rm f}(g, v) = 0
}     &  \forall  v\in V_{{\rm f},0}  ;
\\[3mm]
\displaystyle{
b_{\rm f}(u_{\rm f}, q)
 = 0
}     &  \forall  q\in Q_{{\rm f}}  ;
\\[3mm]
\displaystyle{
\int_{ \widetilde{\Omega}_{\rm s}   }
\rho_{\rm s} 
\frac{\partial^2 d_{\rm s} }{\partial t^2}
\cdot w \, dx
+
R_{\rm s}(d_{\rm s},  w)
+
E_{\rm s}(g, w) = 0
}     &  \forall  w\in V_{{\rm s},0}  ;
\\[3mm]
\displaystyle{
\Phi_{\rm f}  \big|_{\widetilde{\Gamma}} = \texttt{id} + {{d}}_{\rm s}\big|_{\widetilde{\Gamma}},
\quad
u_{\rm f} \circ \Phi_{\rm f} \big|_{\widetilde{\Gamma}}
=
\frac{\partial d_{\rm s}}{\partial t}
\big|_{\widetilde{\Gamma}}
}.
&
  \\
\end{array}
\right.
\end{equation}
Here, $R_{\rm f},R_{\rm s},b_{\rm f}$ are given by
\begin{equation}
\label{eq:weak_formulation_b}
\left\{
\begin{array}{l}
\displaystyle{
R_{\rm f}(u, p, \omega, v)
=
\int_{\Omega_{\rm f}} 
\left(
\sigma_{\rm f}(u,p): \nabla v  \, + \,
\rho_{\rm f} (u - \omega) \cdot \nabla u \cdot v
+
\frac{\rho_{\rm f}}{2}
(\nabla \cdot u) u\cdot v 
\right) \, dx
-
\int_{\Gamma_{\rm f}^{\rm neu}} g_{\rm f}^{\rm neu} \cdot v \, dx;
} \\[3mm]
\displaystyle{
b_{\rm f}(u, q)
= - 
\int_{\Omega_{\rm f}}  
(\nabla \cdot u)  \, q  \, dx
;} \\[3mm]
\displaystyle{
R_{\rm s}(d, w)
=
\int_{\widetilde{\Omega}_{\rm s}} 
\left(
\sigma_{\rm s}(d): \nabla w  \, - \,
f_{\rm s}   \cdot  w
\right) \, dx
-
\int_{\widetilde{\Gamma}_{\rm s}^{\rm neu}} g_{\rm s}^{\rm neu} \cdot w \, dx;
} \\
\end{array}
\right.    
\end{equation}
while $E_{\rm f}$ and $E_{\rm s}$ are defined as follows:
\begin{equation}
\label{eq:weak_formulation_c}
E_{\rm f}(g,v) = - \int_{\widetilde{\Gamma}} g \cdot v\circ \Phi_{\rm f} \, dx,
\qquad
E_{\rm s}(g,w) = \int_{\widetilde{\Gamma}} g \cdot w\, dx.
\end{equation}
\end{subequations}

Some comments are in order. The strongly consistent term 
$\frac{\rho_{\rm f}}{2}
(\nabla \cdot u) u\cdot v $ in \eqref{eq:weak_formulation_b}$_1$ is added to the fluid residual to ensure the energy balance of Lemma \ref{th:stability_properties_continuous} at the discrete level  
\cite{nobile2001numerical}.
We notice that \eqref{eq:weak_formulation_c} implies that
$g=- {\sigma}_{\rm s} \widetilde{n}_{\rm s}$. Therefore, recalling the coupling condition \eqref{eq:fsi_coupling}$_3$ and Nanson's formula, we find
\begin{equation}
\label{eq:coupling_explained}
\int_{\widetilde{\Gamma}} g\cdot v\circ \Phi_{\rm f} \, dx
=
-
\int_{\widetilde{\Gamma}} ({\sigma}_{\rm s} \widetilde{n}_{\rm s}) \cdot v\circ \Phi_{\rm f} \, dx
=
\int_{\widetilde{\Gamma}} 
\left(
J_{\Phi} \left( \sigma_{\rm f} \circ \Phi_{\rm f}  \right) \nabla \Phi_{\rm f}^{-\top} \widetilde{n}_{\rm f}
 \right) \cdot v\circ \Phi_{\rm f} \, dx
=
\int_{{\Gamma}} 
\left(
 \sigma_{\rm f} {n}_{\rm f}
 \right) \cdot v  \, dx,
\end{equation}
which shows the consistency of the weak formulation
\eqref{eq:weak_formulation}
with the strong governing equations.

\section{Discrete formulation}
\label{sec:discrete_equations}

We resort to a nodal $\mathbb{P}_{\kappa}-\mathbb{P}_{\kappa-1}$ Taylor-Hood finite element (FE) discretization for the fluid and a $\mathbb{P}_{\kappa}$ FE discretization for the solid and the control; on the other hand, we consider a backward differentiation formula (BDF) for time integration. In section \ref{sec:FEM}, we introduce the FE discretization and we comment on the enforcement of the coupling conditions.
In section \ref{sec:semi-discrete}, we present the semi-discrete formulation; in section \ref{sec:fully-discrete_formulation} we introduce the fully-discrete formulation; in section \ref{sec:optimal_control_formulation}, we discuss the optimal control formulation that is the point of departure for the MOR formulation of section \ref{sec:ROM}. Finally, in section \ref{sec:discussionHF}, we discuss the connection of the optimization-based statement with related partitioned methods.

\subsection{Finite element spaces}
\label{sec:FEM}
We introduce nodal FE bases $\{   \varphi_j^{\rm f} \}_{j=1}^{N_{\rm u}}$ and 
$\{   \vartheta_j^{\rm f} \}_{j=1}^{N_{\rm p}}$  for fluid velocity and pressure associated to the reference domain $\widetilde{\Omega}_{\rm f}$, respectively.
We define  
$\widetilde{\mathcal{V}}_{\rm f} = {\rm span} \{   \varphi_j^{\rm f} \}_{j=1}^{N_{\rm u}}$ 
and
 $\widetilde{\mathcal{Q}}_{\rm f} = {\rm span} \{   \vartheta_j^{\rm f} \}_{j=1}^{N_{\rm p}}$; we use notation 
$\mathbf{u}_{\rm f}: (0,T] \to \mathbb{R}^{N_{\rm u}}$, 
$\mathbf{p}_{\rm f}: (0,T] \to \mathbb{R}^{N_{\rm p}}$,  to indicate the vectors of coefficients of 
fluid velocity and pressure
\begin{equation}
\label{eq:fluid_field2vec_ref}
\widetilde{u}_{\rm f}(\widetilde{x},t)
=\sum_{j=1}^{N_{\rm u}}
\left( \mathbf{u}_{\rm f}(t) \right)_j \phi_j^{\rm f} (\widetilde{x}),
\qquad
\widetilde{p}_{\rm f}(\widetilde{x},t)
=\sum_{j=1}^{N_{\rm p}}
\left( \mathbf{p}_{\rm f}(t) \right)_j \theta_j^{\rm f} (\widetilde{x}),
\quad
\widetilde{x}\in \widetilde{\Omega}_{\rm f};
\end{equation}
since we consider nodal bases
$\{   \phi_j^{\rm f} \}_{j=1}^{N_{\rm u}}$ and 
$\{   \theta_j^{\rm f} \}_{j=1}^{N_{\rm p}}$ of the FE spaces, the vector of coefficients  correspond to the evaluation of the corresponding field in the nodes of the mesh.
We further define 
the vector of coefficients 
$\boldsymbol{\Phi}_{\rm f}: (0,T] \to \mathbb{R}^{N_{\rm u}}$ associated with the geometric mapping
 such that
\begin{equation}
\label{eq:Phif_field2vec_ref}
\Phi_{\rm f}(\widetilde{x},t)
=\sum_{j=1}^{N_{\rm u}}
\left( \boldsymbol{\Phi}_{\rm f}(t) \right)_j \varphi_j^{\rm f} (\widetilde{x}),
\quad
\widetilde{x}\in \widetilde{\Omega}_{\rm f}.
\end{equation}

Exploiting the previous definitions, we introduce the approximations of the velocity and the pressure field at time $t$.
First, we introduce the time-dependent bases
$\{   \phi_j^{\rm f} (x,t) = \varphi_j^{\rm f}( \Phi_{\rm f}^{-1}(x,t) ) \}_{j=1}^{N_{\rm u}}$ and 
$\{   
\theta_j^{\rm f}(x,t) = 
\vartheta_j^{\rm f}( \Phi_{\rm f}^{-1}(x,t) )
 \}_{j=1}^{N_{\rm p}}$ and we define the corresponding  FE spaces
 ${\mathcal{V}}_{\rm f}(t) = {\rm span} \{   \phi_j^{\rm f}(\cdot, t) \}_{j=1}^{N_{\rm u}}$ and
 ${\mathcal{Q}}_{\rm f} = {\rm span} \{   \theta_j^{\rm f}(\cdot, t) \}_{j=1}^{N_{\rm p}}$; finally, we define the FE velocity and pressure fields
 \begin{equation}
\label{eq:fluid_field2vec_def}
{u}_{\rm f}(x,t)
=\sum_{j=1}^{N_{\rm u}}
\left( \mathbf{u}_{\rm f}(t) \right)_j \phi_j^{\rm f} (x,t),
\qquad
{p}_{\rm f}(x,t)
=\sum_{j=1}^{N_{\rm p}}
\left( \mathbf{p}_{\rm f}(t) \right)_j \theta_j^{\rm f} (x,t),
\quad
x\in  {\Omega}_{\rm f}(t).
\end{equation}
Note that the vectors of coefficients $ \mathbf{u}_{\rm f}$ and  $\mathbf{p}_{\rm f}$ are the same as in \eqref{eq:fluid_field2vec_ref}. Recalling the definition of ALE derivative and ALE velocity, we then find
\begin{equation}
\label{eq:discreteALE_derivative}
\frac{\partial {u}_{\rm f}}{\partial t}
 \Big|_{\Phi_{\rm f}} (x,t)
=\sum_{j=1}^{N_{\rm u}}
\left( \dot{\mathbf{u}}_{\rm f}(t) \right)_j \phi_j^{\rm f} (x,t),
\quad
\omega_{\rm f}(x,t)
=\sum_{j=1}^{N_{\rm u}}
\left( \dot{\boldsymbol{\Phi}}_{\rm f}(t) \right)_j \phi_j^{\rm f} (x,t),
\quad
x\in  {\Omega}_{\rm f}(t).
\end{equation}

\begin{remark}
\label{remark:shape_derivatives}
\textbf{Alternative expression of the ALE mapping.}
We  denote by $\{\widetilde{x}_j^{\rm f} \}_{j=1}^{N_{\rm f,v}}$
and  
$\{ \texttt{D}_k^{\rm f} \}_{k=1}^{N_{\rm f,e}}$
the nodes  and the elements
of the mesh for the velocity field in $\widetilde{\Omega}_{\rm f}$,
and by $\texttt{T}^{\rm f} \in \mathbb{N}^{N_{\rm e,f} \times  n_{\rm lp}}$ the connectivity matrix, where $n_{\rm lp}$ is the number of    nodes in each element.
We further denote by $\Psi_k^{\rm f}: \widehat{\texttt{D}} \to \texttt{D}_k^{\rm f} $ the FE map from the reference element
$\widehat{\texttt{D}}$ to the $k$-th element of the mesh, for $k=1,\ldots, N_{\rm f,e}$.
For any $\widetilde{x}\in \texttt{D}_k^{\rm f}$, we can express $\Phi_{\rm f}$ as follows:
\begin{equation}
\label{eq:Phif_field2vec_ref_ver2}
\Phi_{\rm f} (\widetilde{x},t) = 
\sum_{j \in \{ \texttt{T}_{k,i}^{\rm f} \}_{i=1}^{n_{\rm lp}}}
\left( \boldsymbol{\Phi}_{\rm f}(t) \right)_j \varphi_j^{\rm f} (\widetilde{x})
=
\sum_{i=1}^{n_{\rm lp}} 
\Phi_{\rm f}
\left( 
\widetilde{x}_{\texttt{T}_{k,i}^{\rm f}}^{\rm f} , 
t \right)
\ell_i \left(  
\left( \Psi_k^{\rm fe} \right)^{-1}
(\widetilde{x})  \right),
\end{equation}
where $\ell_1,\ldots,\ell_{n_{\rm lp}}$ is the Lagrangian basis of the polynomial space $\mathbb{P}_{\kappa}$
The expression \eqref{eq:Phif_field2vec_ref_ver2} enables a straightforward  derivation of the shape derivatives.
\end{remark}

We introduce the FE mesh
$\mathcal{T}_{\rm s}$ of the reference solid domain
$\widetilde{\Omega}_{\rm s}$;
we denote by $\{x_j^{\rm s} \}_{j=1}^{N_{\rm s,v}}$
and  
$\{   \texttt{D}_k^{\rm s} \}_{k=1}^{N_{\rm s,e}}$
the nodes  and the elements
of the mesh,  and by
$\texttt{T}^{\rm s}\in \mathbb{N}^{ N_{\rm s,e} \times n_{\rm lp} }$ the connectivity matrix. Given the 
 nodal FE bases $\{   \varphi_j^{\rm s} \}_{j=1}^{N_{\rm s}}$  for the solid displacement, we define  
$ {\mathcal{V}}_{\rm s} = {\rm span} \{   \varphi_j^{\rm s} \}_{j=1}^{N_{\rm s}}$ and we introduce the discrete representation of the displacement field
\begin{equation}
\label{eq:solid_field2vec_ref}
\widetilde{d}_{\rm s}(x,t)
=\sum_{j=1}^{N_{\rm s}}
\left( \mathbf{d}_{\rm s}(t) \right)_j \varphi_j^{\rm s} (x),
\quad
x\in \widetilde{\Omega}_{\rm s},
\;\;
t\in (0,T].
\end{equation}

\begin{remark}
\label{remark:coupling_discrete}
\textbf{Coupling conditions.}
We define the space for the control variable  
$\mathcal{G}$ as follows:
\begin{equation}
\label{eq:conforming_mesh_assumption}
{\mathcal{G}}  := 
\left\{
  \varphi ^{\rm f} \big|_{\widetilde{\Gamma}} \, : \, 
  \varphi ^{\rm f} \in \widetilde{\mathcal{V}}_{\rm f} 
\right\}.
\end{equation}
We assume that the   solid and the  fluid meshes are conforming at the fluid-solid interface: this implies that ${\mathcal{G}}  = 
\left\{
  \varphi ^{\rm s} \big|_{\widetilde{\Gamma}} \, : \, 
  \varphi ^{\rm s} \in \widetilde{\mathcal{V}}_{\rm s} 
\right\}$.
Furthermore, if we introduce the mask matrices $\mathbf{P}_{\rm f} \in \{0,1 \}^{N_{\Gamma}\times N_{\rm u}}$ and $\mathbf{P}_{\rm s} \in \{0,1 \}^{N_{\Gamma}\times N_{\rm s}}$ that extract the degrees of freedom  of the fluid and the solid FE spaces  on the interface $\Gamma$ --- provided that the numbering is consistent --- we find the
discrete counterparts of \eqref{eq:fsi_coupling}$_1$ and \eqref{eq:fsi_coupling}$_2$:
\begin{equation}
\label{eq:coupling_strong_discrete}
\mathbf{P}_{\rm f} \mathbf{d}_{\rm f} =     \mathbf{P}_{\rm s} \mathbf{d}_{\rm s},
\qquad
\mathbf{P}_{\rm f} \mathbf{u}_{\rm f} =     \mathbf{P}_{\rm s} \dot{\mathbf{d}}_{\rm s}.
\end{equation}
\end{remark}

\subsection{Semi-discrete formulation}
\label{sec:semi-discrete}
Exploiting the previous definitions, we can introduce the semi-discrete FE formulation.
First, we define the spaces
$\mathcal{V}_{\rm f}^{\rm dir} := \{v\in
\mathcal{V}_{\rm f} 
 \, : \, v\big|_{{\Gamma}_{\rm f}^{\rm dir}} = u_{\rm f}^{\rm dir}  \}$,
 $\mathcal{V}_{{\rm f}}^0 := \{v\in
\mathcal{V}_{\rm f} 
 \, : \, v\big|_{\widetilde{\Gamma}_{\rm f}^{\rm dir}} = 0 \}$ and
$\mathcal{V}_{\rm s}^{\rm dir} := \{v\in
\mathcal{V}_{\rm s} 
 \, : \, v\big|_{\widetilde{\Gamma}_{\rm s}^{\rm dir}} = d_{\rm s}^{\rm dir}(t) \}$,
 $\mathcal{V}_{{\rm s}}^0 := \{v\in
\mathcal{V}_{\rm s} 
 \, : \, v\big|_{\widetilde{\Gamma}_{\rm s}^{\rm dir}} = 0 \}$. Then,  
for any $t>0$, we seek $u_{\rm f}\in \mathcal{V}_{\rm f}^{\rm dir}$, $p_{\rm f}\in \mathcal{Q}_{\rm f}$,
$d_{\rm s}\in \mathcal{V}_{\rm s}^{\rm dir}$ and 
$g \in \mathcal{G}$ such that
\begin{equation}
\label{eq:weak_formulation_semidiscrete}
\left\{
\begin{array}{ll}
\displaystyle{
\int_{ \Omega_{\rm f}   }
\rho_{\rm f} 
\frac{\partial u_{\rm f} }{\partial t}
\Big|_{\Phi_{\rm f}} \cdot v \, dx
+
R_{\rm f}(u_{\rm f}, p_{\rm f}, \omega_{\rm f}, v)
+
E_{\rm f}(g, v) = 0
}     &  \forall  v\in  \mathcal{V}_{{\rm f}}^0  ;
\\[3mm]
\displaystyle{
b_{\rm f}(u_{\rm f}, q)
 = 0
}     &  \forall  q\in  \mathcal{Q}_{{\rm f}}  ;
\\[3mm]
\displaystyle{
\int_{ \widetilde{\Omega}_{\rm s}   }
\rho_{\rm s} 
\frac{\partial^2 d_{\rm s} }{\partial t^2}
\cdot w \, dx
+
R_{\rm s}(d_{\rm s},  w)
+
E_{\rm s}(g, w) = 0
}     &  \forall  w\in  \mathcal{V}_{\rm s}^0  ;
\\[3mm]
\displaystyle{
\Phi_{\rm f}  \big|_{\widetilde{\Gamma}} = \texttt{id} + {{d}}_{\rm s}\big|_{\widetilde{\Gamma}},
\quad
u_{\rm f} \circ \Phi_{\rm f} \big|_{\widetilde{\Gamma}}
=
\frac{\partial d_{\rm s}}{\partial t}
\big|_{\widetilde{\Gamma}}
};
&
  \\
\end{array}
\right.
\end{equation}
where the linear and nonlinear forms
$R_{\rm f}, E_{\rm f},b_{\rm f},
R_{\rm s}, E_{\rm s}$ are defined in 
\eqref{eq:weak_formulation_b} and 
\eqref{eq:weak_formulation_c}.

Exploiting the same argument as in \cite[Prop. 9.1]{formaggia2010cardiovascular}, we can prove the following stability result for the semi-discrete formulation; we postpone the proof to \ref{sec:proofs}. We remark that the stability result relies on the introduction of the   
density of elastic energy $W$ that depends on the prescribed constitutive law: we refer to 
\cite{formaggia2010cardiovascular} for the expression of $W$ for the two constitutive laws of  
\eqref{eq:solid_model}.
 
\begin{lemma}
\label{th:stability_properties_continuous}
Let the coupled fluid-structure system \eqref{eq:weak_formulation_semidiscrete} be isolated; i.e., ${u}_{\rm f}=0$ on $\partial \Omega_{\rm f}(t)\setminus\Gamma(t)$, $ \widetilde{\sigma}_{\rm s} \widetilde{n}_{\rm s}=0$ on $\partial \widetilde{\Omega}_{\rm s}\setminus\widetilde{\Gamma}$.
Assume that $ \widetilde{\sigma}_{\rm s} = \frac{\partial W}{\partial F_{\rm s}}(F_{\rm s})$ for a given density of elastic energy $W$.
Then the following energy balance holds:
\begin{equation}
\label{eq:stability_properties_continuous}
\frac{d}{dt}\left[\int_{\Omega_{\rm f}} \frac{\rho_{\rm f}}{2} \vert{u}_{\rm f} \vert^2\,dx  + \int_{\widetilde{\Omega}_{\rm s}} 
\frac{\rho_{\rm s} }{2}
 \vert  \frac{\partial d_{\rm s}}{\partial t}     \vert^2 \, dx +
\int_{\widetilde{\Omega}_{\rm s}} 
 W (F_{\rm s})\, dx\right] + \int_{\Omega_{\rm f}} 
2\mu_{\rm f}   \epsilon_{\rm f}({u}_{\rm f}) :  \epsilon_{\rm f}({u}_{\rm f})
=0.
\end{equation}
\end{lemma}

\subsection{Fully-discrete formulation}
\label{sec:fully-discrete_formulation}
We denote by  $\{ t^{(k)} = \Delta t (k-1)\}_{k=1}^K$ the time grid and by 
 ${d}_{\rm s}^{(k)}$,
  ${u}_{\rm f}^{(k)}$, $p_{\rm f}^{(k)}$  and
    ${g}^{(k)}$
 the state and the control estimates  at time $t^{(k)}$, respectively.
Following 
 \cite{chabannes2013high,deparis2016facsi},
 we consider the Newmark scheme for the solid system
 \begin{equation}
 \label{eq:newmark}
 \left\{
 \begin{array}{l}
 \displaystyle{
\texttt{D}_{{\rm s},\Delta t}^2 d_{\rm s}^{(k+1)}
:=
\frac{1}{\beta \Delta t^2}
\left(
d_{\rm s}^{(k+1)}
 -
d_{\rm s}^{(k)}
\right)
-
\frac{1}{\beta \Delta t}
\texttt{D}_{{\rm s},\Delta t} d_{\rm s}^{(k)}
-
\left(
\frac{1}{2\beta} - 1
\right)
\texttt{D}_{{\rm s},\Delta t}^2 d_{\rm s}^{(k)}
 }
 \\[3mm]
  \displaystyle{
\texttt{D}_{{\rm s},\Delta t} d_{\rm s}^{(k+1)}
:=
\frac{\gamma}{\beta \Delta t}
\left(
d_{\rm s}^{(k+1)}
 -
d_{\rm s}^{(k)}
\right)
-
\left(
\frac{\gamma}{\beta} - 1
\right)
\texttt{D}_{{\rm s},\Delta t} d_{\rm s}^{(k)}
-
\Delta t
\left(
\frac{\gamma}{2\beta} - 1
\right)
\texttt{D}_{{\rm s},\Delta t}^2 d_{\rm s}^{(k)}
},
 \\
 \end{array}
\right.
\end{equation}    
  where 
  $\gamma,\beta$ are positive constants. 
The Newmark method is 
unconditionally stable for   $2\beta \geq \gamma \geq \frac{1}{2}$ and 
second-order accurate   for   $\gamma=\frac{1}{2}$ and  $\beta=\frac{1}{4}$ (see, e.g., 
  \cite{rao2010finite}); increasing the value of $\gamma$ leads to an increase of numerical damping and might be beneficial for certain applications.
  For $k=0$, 
  since $D_{{\rm s},\Delta t}^2 d_{\rm s}^{(0)}$ is not available, we consider $\gamma=1$ and $\beta=\frac{1}{2}$.
This choice of the parameters leads to an  unconditionally stable and second-order accurate scheme. 
For the fluid subproblem, we consider the backward differentiation formulas of order one and two (BDF1, BDF2)
\begin{subequations}
\label{eq:BDF}
\begin{equation}
\texttt{D}_{{\rm f},\Delta t} u_{\rm f}^{(k+1)}(x)
=
\sum_{i=1}^{N_{\rm u}}
\left(
\texttt{D}_{{\rm f},\Delta t} \mathbf{u}_{\rm f}^{(k+1)}
\right)_j \varphi_{j,k+1}^{\rm f}(x),
\quad
x\in \Omega_{\rm f}^{(k+1)} = \Omega_{\rm f}(t^{(k+1)}),
\end{equation}
where $\varphi_{1,k+1}^{\rm f},\ldots, \varphi_{N_{\rm u},k+1}^{\rm f}$ are the Lagrangian bases in the domain $\Omega_{\rm f}^{(k+1)}$ and 
\begin{equation}
\texttt{D}_{{\rm f},\Delta t} \mathbf{u}_{\rm f}^{(k+1)} :=
 \left\{
 \begin{array}{ll}
 \displaystyle{
\frac{1}{\Delta t}
\left(
 {\mathbf{u}}_{\rm f}^{(k+1)}
 -
  {\mathbf{u}}_{\rm f}^{(k)}
\right)
 }
 &
 {\rm BDF1};
 \\[3mm]
  \displaystyle{
\frac{1}{2 \Delta t}
\left(
3  {\mathbf{u}}_{\rm f}^{(k+1)}
 -
 4 {\mathbf{u}}_{\rm f}^{(k)}
 +
 {\mathbf{u}}_{\rm f}^{(k-1)}
\right) 
}
&
 {\rm BDF2}.
 \\
 \end{array}
\right.
\end{equation}
\end{subequations}
Note that the time discretization formula is applied to the coefficients of \eqref{eq:fluid_field2vec_def} as prescribed by \eqref{eq:discreteALE_derivative}. Finally, we introduce the discretized ALE velocity as (see \eqref{eq:discreteALE_derivative}):
\begin{equation}
\label{eq:discreteALE_velocity}
\omega_{\rm f}^{(k+1)}(x)
=
\sum_{j=1}^{N_{\rm u}} \left( \boldsymbol{\omega}_{\rm f}^{(k+1)}  \right) \varphi_{j,k+1}^{\rm f}(x),
\quad
{\rm where} \;\;
\boldsymbol{\omega}_{\rm f}^{(k+1)}  = \frac{\partial \boldsymbol{\Phi}_{\rm f}^{(k+1)}}{\partial \mathbf{d}_{\rm s}} 
\left[\mathbf{d}_{\rm s}^{(k)} \right] \; 
 \texttt{D}_{{\rm s},\Delta t}  \mathbf{d}_{\rm s}^{(k+1)}.
\end{equation}
  Since we consider the pseudo-elastic model  \eqref{eq:fsi_mesh_pseudo_elasticity} for  the ALE morphing, we can express 
  $\boldsymbol{\Phi}_{\rm f}$ as
  $$
  \boldsymbol{\Phi}_{\rm f}^{(k)}  = \boldsymbol{\texttt{id}} + \mathbf{W} \mathbf{P}_{\rm s} \mathbf{d}_{\rm s}^{(k)}
  \quad
  {\rm and \;\; thus} \;\;
  \boldsymbol{\omega}_{\rm f}^{(k+1)}  = \mathbf{W} \mathbf{P}_{\rm s}  \texttt{D}_{{\rm f},\Delta t} \mathbf{d}_{\rm s}^{(k)},  
  $$
  for a properly-chosen matrix  $ \mathbf{W}  \in \mathbb{R}^{N_{\rm u} \times N_\Gamma}$ that is explicitly defined in Appendix \ref{sec:algebraic_formula}. 
We also tried to replace $\texttt{D}_{{\rm f},\Delta t} \mathbf{d}_{\rm s}^{(k)}$ (BDF approximation of the time derivative) with 
$\texttt{D}_{{\rm s},\Delta t} \mathbf{d}_{\rm s}^{(k)}$ (Newmark approximation of the time derivative) in \eqref{eq:discreteALE_velocity}: 
the two methods lead to similar performance for all the tests considered in section \ref{sec:numerical_results}.

We have now the elements to introduce the fully-discrete formulation: 
we introduce the affine  spaces  $\mathcal{V}_{{\rm f},k+1}^{\rm dir},
\mathcal{V}_{{\rm s},k+1}^{\rm dir}$, the  linear spaces 
$\mathcal{V}_{{\rm f},k+1}^0,\mathcal{Q}_{{\rm f},k+1}$ at time $t^{(k+1)}$, 
and the initial conditions 
$u_{\rm f}^{(0)}, d_{\rm s}^{(0)}, \frac{\partial d_{\rm s}}{\partial t}^{(0)}$; then,  
for $k=0,1,\ldots$, 
we seek $u_{\rm f}^{(k+1)}\in \mathcal{V}_{{\rm f},k+1}^{\rm dir}$, $p_{\rm f}\in \mathcal{Q}_{{\rm f},k+1}$,
$d_{\rm s}\in \mathcal{V}_{{\rm s},k+1}^{\rm dir}$ and
$g \in \mathcal{G}$ such that
\begin{equation}
\label{eq:weak_formulation_discrete}
\left\{
\begin{array}{ll}
\displaystyle{
\int_{ \Omega_{\rm f}^{(k+1)}   }
\rho_{\rm f} 
\texttt{D}_{{\rm f},\Delta t} u_{\rm f}^{(k+1)}  \cdot v \, dx
+
R_{\rm f}^{(k+1)}(u_{\rm f}^{(k+1)}, p_{\rm f}^{(k+1)}, \omega_{\rm f}^{(k+1)}, v)
+
E_{\rm f}(g^{(k+1)}, v) = 0
}     &  \forall  v\in  \mathcal{V}_{{\rm f},k+1}^0  ;
\\[4mm]
\displaystyle{
b_{\rm f}(u_{\rm f}^{(k+1)}, q)
 = 0
}     &  \forall  q\in  \mathcal{Q}_{{\rm f},k+1};
\\[4mm]
\displaystyle{
\int_{\widetilde{\Omega}_{\rm s}}
\rho_{\rm s} 
\texttt{D}_{{\rm s},\Delta t}^2  d_{\rm s}^{(k+1)}  
\cdot w \, dx
+
R_{\rm s}^{(k+1)}(d_{\rm s}^{(k+1)},  w)
+
E_{\rm s}(g, w) = 0
}     &  
\forall  w\in  \mathcal{V}_{\rm s}^0  ;
\\[4mm]
\displaystyle{
\Phi_{\rm f}^{(k+1)}  \big|_{\widetilde{\Gamma}} = \texttt{id} + {{d}}_{\rm s}^{(k+1)}\big|_{\widetilde{\Gamma}},
\quad
u_{\rm f}^{(k+1)} \circ \Phi_{\rm f}^{(k+1)} \big|_{\widetilde{\Gamma}}
=
\texttt{D}_{{\rm s},\Delta t}  d_{\rm s}^{(k+1)}  
\big|_{\widetilde{\Gamma}}
}.
&
  \\
\end{array}
\right.
\end{equation}

Some comments are in order.
The choice of the BDF scheme for the fluid subproblem and of the Newmark scheme for the solid subproblem are broadly used in the literature; however, we emphasize that our formulation can cope with virtually any time integrator scheme.
By tedious but trivial calculations, we can express \eqref{eq:weak_formulation_discrete} in the following algebraic form:
\begin{equation}
\label{eq:weak_formulation_discrete_algebraic}
\left\{
\begin{array}{ll}
\displaystyle{\mathbf{P}_{\rm f} \mathbf{u}_{\rm f}^{(k+1)} =  
\mathbf{P}_{\rm s}  \texttt{D}_{{\rm s},\Delta t} \mathbf{d}_{\rm s}^{(k+1)}
}
&
\\[3mm]
\displaystyle{
\mathbf{R}_{\rm f}^{(k+1)}
\left(
\mathbf{u}_{\rm f}^{(k+1)},
\mathbf{p}_{\rm f}^{(k+1)},
\mathbf{d}_{\rm s}^{(k+1)}
\right)
+\mathbf{E}_{\rm f} \mathbf{g}^{(k+1)}
 =  
0},
&
\displaystyle{\mathbf{P}_{\rm f,\rm dir} \mathbf{u}_{\rm f}^{(k+1)}  = \mathbf{u}_{\rm f}^{{\rm dir}, (k+1)}
};
\\[3mm]
\displaystyle{
\mathbf{R}_{\rm s}^{(k+1)}
\left(
\mathbf{d}_{\rm s}^{(k+1)}
\right)
+\mathbf{E}_{\rm s} \mathbf{g}^{(k+1)}
 =  
0
}
&
\displaystyle{\mathbf{P}_{\rm s,\rm dir} \mathbf{d}_{\rm s}^{(k+1)}  = \mathbf{d}_{\rm s}^{{\rm dir}, (k+1)}
}.
\\[3mm]
\end{array}
\right.
\end{equation}
To ease the presentation, we 
postpone the definition of the terms in \eqref{eq:weak_formulation_discrete_algebraic} to Appendix \ref{sec:algebraic_formula}.
A rigorous stability analysis for \eqref{eq:weak_formulation_discrete} is beyond the scope of the present work; we refer to \cite[Proposition 9.4]{formaggia2010cardiovascular} for a related result.

\subsection{Optimal control formulation}
\label{sec:optimal_control_formulation}

We introduce
notation $\texttt{D}_{{\rm s},\Delta t} \mathbf{d}_{\rm s}^{(k+1)} = \alpha_{\rm s,v}^{(k+1)} \mathbf{d}_{\rm s}^{(k+1)} + \mathbf{b}_{\rm s,v}^{(k+1)}$
for the finite difference approximation of the first-order time derivative,
for proper choices of $ \alpha_{\rm s,v}^{(k+1)}\in \mathbb{R}$ and $\mathbf{b}_{\rm s,v}^{(k+1)}\in \mathbb{R}^{N_{\rm s}}$.
Then, we consider the 
minimization statement 
\begin{equation}
\label{eq:optimal_control_algebraic}
\min_{\mathbf{u}_{\rm f}, \mathbf{p}_{\rm f}, \mathbf{d}_{\rm s},
 \mathbf{g}}
\frac{1}{2} \Big| 
 \mathbf{P}_{\rm f} \mathbf{u}_{\rm f}  - 
\mathbf{P}_{\rm s}  
\left(
\alpha_{\rm s,v}^{(k+1)} \mathbf{d}_{\rm s}^{(k+1)} + \mathbf{b}_{\rm s,v}^{(k+1)}
\right)
\Big|^2
\, + \,
\frac{\delta}{2} \vertiii{\mathbf{g}}^2,
\quad
{\rm s.t.} \;\;
\left\{
\begin{array}{l}
\displaystyle{
\mathbf{R}_{\rm f}^{(k+1)}
\left(
\mathbf{u}_{\rm f},
\mathbf{p}_{\rm f},
\mathbf{d}_{\rm s}
\right)
+\mathbf{E}_{\rm f} \mathbf{g}
 =  
0} ; \\[3mm]
\displaystyle{\mathbf{P}_{\rm f,\rm dir} \mathbf{u}_{\rm f}  = \mathbf{u}_{\rm f}^{{\rm dir}, (k+1)}
};
\\[3mm]
\displaystyle{
\mathbf{R}_{\rm s}^{(k+1)}
\left(
\mathbf{d}_{\rm s} 
\right)
+\mathbf{E}_{\rm s} \mathbf{g} 
 =  
0
};
\\[3mm]
\displaystyle{\mathbf{P}_{\rm s,\rm dir} \mathbf{d}_{\rm s} = \mathbf{d}_{\rm s}^{{\rm dir}, (k+1)}
}.
\\
\end{array}
\right.
\end{equation}
 The second term in the objective function of \eqref{eq:optimal_control_algebraic} is a regularizer that is designed to penalize controls of 
excessive size; the positive constant $\delta$ is chosen to control the relative importance of the two terms in the objective.
Following \cite{Taddei2024optimization}, we consider the  $H^1(\widetilde{\Gamma})$ seminorm for regularization, that is
$\vertiii{\mathbf{g}}^2 = \int_{\widetilde{\Gamma}} \Big|  \nabla_{\widetilde{\Gamma}} g \Big|^2 \, dx$:
computation of the $H^1(\widetilde{\Gamma})$ seminorm is performed by representing $g$ over a (\texttt{d}-1)-dimensional grid of 
$\widetilde{\Gamma}$. 

We observe that for $\delta=0$ any solution to 
\eqref{eq:weak_formulation_discrete_algebraic} is a global solution to \eqref{eq:optimal_control_algebraic}; furthermore, it is possible to show  that the solution to
\eqref{eq:optimal_control_algebraic} converges to 
\eqref{eq:weak_formulation_discrete_algebraic} 
 as the regularization parameter $\delta $  approaches $ 0$
\cite{Kuberry_Lee_2013,Kuberry_Lee_2015,Kuberry_Lee_2016}.
We further observe that
given $\mathbf{g}\in \mathbb{R}^{N_{\rm c}}$, we can determine the corresponding solution $\mathbf{d}_{\rm s}$ in the solid domain 
 by solving 
\eqref{eq:weak_formulation_discrete_algebraic}$_3$; then, we can compute the solution
$\mathbf{u}_{\rm f}, \mathbf{p}_{\rm f}$ in the fluid domain by solving
\eqref{eq:weak_formulation_discrete_algebraic}$_2$.
This observation motivates the static condensation approach discussed below.

We rely on sequential quadratic programming (SQP) to approximate the solution to \eqref{eq:optimal_control_algebraic}. At each sub-iteration $j$ of the algorithm, we solve the minimization problem
\begin{subequations}
\label{eq:SQP_HF}
\begin{equation}
    \min_{\mathbf{u}_{\rm f}, \mathbf{p}_{\rm f}, \mathbf{d}_{\rm f},
 \mathbf{g}}
\frac{1}{2} \Big| 
 \mathbf{P}_{\rm f} \mathbf{u}_{\rm f}  - 
\mathbf{P}_{\rm s}  
\left(
\alpha_{\rm s,v}^{(k+1)} \mathbf{d}_{\rm s}^{(k+1)} + \mathbf{b}_{\rm s,v}^{(k+1)}
\right)
\Big|^2
\, + \,
\frac{\delta}{2} \vertiii{\mathbf{g}}^2,
\end{equation}
subject to the linearized constraints
\begin{equation}
\label{eq:SQP_HF_b}
\left\{
\begin{array}{l}
\displaystyle{
\mathbf{R}_{\rm f}^{(k+1,j)}
+ 
\mathbf{J}_{\rm f,u}^{(k+1,j)}
\left( \mathbf{u}_{\rm f} - \mathbf{u}_{\rm f}^{(k+1,j)} \right)
+ 
\mathbf{J}_{\rm f,p}^{(k+1,j)}
\left( \mathbf{p}_{\rm f} - \mathbf{p}_{\rm f}^{(k+1,j)} \right)
+ 
\mathbf{J}_{\rm f,s}^{(k+1,j)}
\mathbf{P}_{\rm s}
\left( \mathbf{d}_{\rm s} - \mathbf{d}_{\rm s}^{(k+1,j)} \right)
+\mathbf{E}_{\rm f} \mathbf{g}
 =  
0} ,\\[3mm]
\displaystyle{\mathbf{P}_{\rm f,\rm dir} \mathbf{u}_{\rm f}  = \mathbf{u}_{\rm f}^{{\rm dir}, (k+1)}
},
\\[3mm]
\displaystyle{
\mathbf{R}_{\rm s}^{(k+1,j)}
+ 
\mathbf{J}_{\rm s,s}^{(k+1,j)}
\left(
\mathbf{d}_{\rm s} 
- 
\mathbf{d}_{\rm s}^{(k+1,j)} 
\right)
+\mathbf{E}_{\rm s} \mathbf{g} 
 =  
0
},
\\[3mm]
\displaystyle{\mathbf{P}_{\rm s,\rm dir} \mathbf{d}_{\rm s} = \mathbf{d}_{\rm s}^{{\rm dir}, (k+1)}
};
\\
\end{array}
\right.    
\end{equation}
\end{subequations}
where 
$\mathbf{J}_{\rm f,u}^{(k+1,j)}$
(resp., 
$\mathbf{J}_{\rm f,p}^{(k+1,j)}$)
denotes the Jacobian of the fluid residual with respect to
the velocity field (resp., the pressure field) evaluated at 
the current state;
similarly, $\mathbf{J}_{\rm f,s}^{(k+1,j)}$
(resp.,
$\mathbf{J}_{\rm s,s}^{(k+1,j)}$)
is the Jacobian of the fluid (resp., solid) residual with respect to the 
solid displacement.
First, we express\footnote{This operation is usually referred to as \emph{reduced system approach} in optimization or \emph{static condensation} in scientific computing.}  
$\mathbf{u}_{\rm f}, \mathbf{p}_{\rm f}, \mathbf{d}_{\rm f}$ as a function of $\mathbf{g}$; then, we solve a least-square problem to find $\mathbf{g}$.
We stop the iterations when the relative increment in the control is below a user-defined threshold
\begin{equation}
\label{eq:termination_condition}
\dfrac{\|  \mathbf{g}^{(j+1) } - \mathbf{g}^{(j) }   \|_{\rm c}}{\|    \mathbf{g}^{(j) }   \|_{\rm c}}
< {\rm tol}_{\rm sqp}.
\end{equation}

The   SQP procedure is elementary and might be significantly improved by resorting to line search and more sophisticated approximations of the Hessian of the Lagrangian in the objective function \cite{bonnans2006numerical}: the development of specialized SQP procedures for this class of problems is the subject of ongoing research.

We notice that if we neglect 
$\mathbf{J}_{\rm f,s}^{(k+1,j)}$, the two subproblems can be solved independently; we further observe that the matrix
 $\mathbf{J}_{\rm f,s}^{(k+1,j)}$ is dense and is hence expensive to evaluate and to store. For these reasons, in the numerical experiments, we investigate performance of SQP with both exact Jacobian and inexact Jacobian, where we neglect 
$\mathbf{J}_{\rm f,s}^{(k+1,j)}$ in \eqref{eq:SQP_HF_b}. The same study has been conducted for a related implicit method in 
\cite[section 3]{fernandez2011coupling}.

\subsection{Discussion}
\label{sec:discussionHF}
As extensively discussed in \cite{iollo2023one} and \cite{Taddei2024optimization}, we can interpret both the SQP algorithm \eqref{eq:SQP_HF} and the standard Dirichlet-to-Neumann (DtN) (e.g., \cite{formaggia2010cardiovascular}) algorithm as iterative solvers for the optimal control problem \eqref{eq:optimal_control_algebraic}.
\begin{itemize}
\item 
The DtN  algorithm requires Newton subiterations to solve the local problems, while the linearization of the constraints in SQP eliminates the need for Newton-type subiterations.
\item 
The DtN algorithm involves alternating solutions of the fluid and the structure problem and is hence inherently sequential; on the other hand, our approach with inexact Jacobian calculation enables the independent solution of the fluid and the structure subproblem. In this respect, our method shares features with implicit FSI methods based on the Steklov-Poincaré operator   \cite{Deparis2006_steklov}.
\item 
As opposed to DtN algorithms, our solution method relies on explicit computations of the sensitivity matrices
\begin{equation}
\label{eq:HF_sensitivity_matrix}
\left[
\begin{array}{cc}
\mathbf{J}_{\rm f,u}^{(k+1,j)}     &   \mathbf{J}_{\rm f,p}^{(k+1,j)}   
 \\
\mathbf{P}_{\rm f,dir}^{(k+1,j)}        &  0  \\
\end{array}
\right]^{-1}    
\left[
\begin{array}{c}
\mathbf{E}_{\rm f}    
 \\
0  \\
\end{array}
\right],
\quad
\left[
\begin{array}{c}
\mathbf{J}_{\rm s,s}^{(k+1,j)}        
 \\
\mathbf{P}_{\rm s,dir}^{(k+1,j)}   \\
\end{array}
\right]^{-1}    
\left[
\begin{array}{c}
\mathbf{E}_{\rm s}    
 \\
0  \\
\end{array}
\right],
\end{equation}
at each SQP iteration.
Clearly, computation of \eqref{eq:HF_sensitivity_matrix}
 is only feasible for moderate values of the number of control variables $N_{\rm c}$: this observation highlights the relevance of the present formulation for MOR applications or hybrid formulations involving dimensionality reduction of the control space.
 \end{itemize}

 In \cite{iollo2023one} and also \cite{Gunzburger_Lee_2000}, the authors rely on the Gauss-Newton method (GNM) to solve the optimization problem \eqref{eq:optimal_control_algebraic}. As discussed in \cite{Taddei2024optimization}, GNM requires Newton subiterations at each optimization step; it  might hence experience convergence issues due to the failure of Newton's method at early iterations of the optimization process. For this reason, we exclusively rely on the SQP method. 

As discussed in the introduction, our formulation does not explicitly require compatibility conditions between the fluid and the solid mesh: for non-compatible meshes --- or equivalently for non-compatible reduced spaces --- 
the optimization process  guarantees the approximate satisfaction of the velocity continuity constraint and the continuity of the normal stresses.
Recalling
\eqref{eq:fsi_mesh_pseudo_elasticity}, since  we consider the same mesh for the ALE displacement  as for the fluid states, the continuity of the displacement is also imposed strongly.
Thanks to these observations, we can   apply independent dimensionality reduction techniques for the two subproblems (cf. section \ref{sec:data_compression}).

\section{Reduced-order formulation}
\label{sec:ROM}

In this section, we define a partitioned ROM for the FSI problem \eqref{eq:strong_equations}
based on the optimal control formulation
\eqref{eq:optimal_control_algebraic}.
In this  work, we focus on the \emph{solution reproduction problem}: first, we simulate the FSI problem using the FOM; second, we build the reduced spaces for the fluid and the solid state, and for the control; third, we formulate the local ROMs; fourth,   we solve the ROM for the same parameters considered during the offline stage.
The problem is of little practical interest; however, it represents the first step towards the development of a predictive MOR framework for FSI problems. Algorithm \ref{alg:SRP} summarizes the procedure.

\begin{algorithm}[H]                      
\caption{Solution reproduction problem.}     
\label{alg:SRP}    
\begin{algorithmic}[1]
\State
Solve \eqref{eq:optimal_control_algebraic} and store $\{ ( u_{\rm f}^{(k)},p_{\rm f}^{(k)}, d_{\rm s}^{(k)}, g^{(k)}  )  \}_{k\in \texttt{I}_{\rm smp}}$ with 
$ \texttt{I}_{\rm smp} = \{ i_j  \}_{j=1}^{\rm train} \subset \{1,\ldots,K\}$.
\vspace{3pt}

\State
Generate the reduced spaces for states and control (cf. section \ref{sec:data_compression}).
\vspace{3pt}

\State
Formulate the local reduced-order models
(cf. section \ref{sec:local_roms}).
\vspace{3pt}

\State
Estimate the generalized coordinates by solving the global   reduced-order model
(cf. section \ref{sec:global_rom}).
\end{algorithmic}
\end{algorithm}

\subsection{Data compression}
\label{sec:data_compression}

We denote by ${w}_{\rm f}={\rm vec}(u_{\rm f}, p_{\rm f})$ the full state vector in the fluid domain;
we further define the tensor product space
$\widetilde{\mathcal{X}}_{\rm f}:=
\widetilde{\mathcal{V}}_{\rm f}^0 \times \widetilde{\mathcal{Q}}_{\rm f}$.
We denote by $\overline{w}_{\rm f}^{(k)}$ a lift function that satisfies the Dirichlet boundary conditions and 
$\overline{w}_{\rm f}^{(0)}=w_{\rm f}^{(0)}$. In the numerical experiments, we consider problems with separable Dirichlet conditions,   $u_{\rm f}^{\rm dir}(x,t)=c_{\rm f}^{\rm dir}(t) 
\bar{u}_{\rm f}^{\rm dir}(x)$: we might hence set 
$\widetilde{\overline{w}}_{\rm f}^{(k)} = \bar{c}_{\rm f}^{\rm dir}(t^{(k)}) 
\widetilde{\bar{w}}_{\rm f}$ with 
$\widetilde{\bar{w}}_{\rm f}:=\frac{1}{K} \sum_{k=1}^K \widetilde{w}_{\rm f}^{(k)}$
and
$\bar{c}_{\rm f}^{\rm dir}(t) = \frac{K}{\sum_{k=1}^K {c}_{\rm f}^{\rm dir}(t^{(k)})} {c}_{\rm f}^{\rm dir}(t).$ We apply the same argument to devise a lift function $\overline{d}_{\rm s}$ for the  solid state.
We refer to \cite{taddei2021discretize} for the corresponding ansatz for non-homogeneous Dirichlet conditions.

Exploiting the previous definitions, we consider the following ansatze for the states and 
the control variables:
\begin{equation}
\label{eq:ansatz}
\widetilde{w}_{\rm f}^{(k)}
\approx 
\widetilde{\overline{w}}_{\rm f}^{(k)} +  \widetilde{Z}_{\rm f} \widehat{\boldsymbol{\alpha}}_{\rm f}^{(k)},
\quad
d_{\rm s}^{(k)}
\approx 
d_{\rm s}^{(0)} + Z_{\rm s} \widehat{\boldsymbol{\alpha}}_{\rm s}^{(k)},
\quad
g^{(k)}
\approx 
 Z_{\rm c} \widehat{\boldsymbol{\beta}}^{(k)},
 \quad
 k=1,\ldots,K.
\end{equation}
Note that the  affine ansate 
\eqref{eq:ansatz}
for the states ensure the satisfaction of the initial condition and of the strong boundary conditions.
Here, $\widetilde{Z}_{\rm f} = [\widetilde{\zeta}_1^{\rm f}, \ldots, \widetilde{\zeta}_{n_{\rm f}}^{\rm f}]:\mathbb{R}^{n_{\rm f}} \to \widetilde{\mathcal{X}}_{\rm f}$,
 ${Z}_{\rm s} = [{\zeta}_1^{\rm s}, \ldots, 
 {\zeta}_{n_{\rm s}}^{\rm s}]:\mathbb{R}^{n_{\rm s}} \to  {\mathcal{V}}_{\rm s}^0$, 
 ${Z}_{\rm c} = [{\zeta}_1^{\rm c}, \ldots, 
 {\zeta}_{n_{\rm c}}^{\rm c}]:\mathbb{R}^{n_{\rm c}} \to  \mathcal{G}$ are suitable linear operators.
 We recall that for the fluid domain we use the symbol 
 $\widetilde{\bullet}$ to refer to quantities (fields, approximation spaces) defined in the reference domain;
 given the state  $\widetilde{u}$ defined over $\widetilde{\Omega}_{\rm f}$, we use notation $u$ to refer to the corresponding field in the physical domain  (cf. section \ref{sec:FEM}).

We rely on  POD   based on the method of snapshots 
 \cite{sirovich1987turbulence}
 to determine the reduced spaces; we determine  the size of the spaces  according to an energy criterion. 
We denote by $(\cdot, \cdot)_{\rm f},(\cdot, \cdot)_{\rm s}, (\cdot, \cdot)_{\rm c}$ the inner products for the fluid state, the solid state and the control,  we denote by $\mathring{{\widetilde{w}}}_{\rm f}^{(k)} =
{{\widetilde{w}}}_{\rm f}^{(k)} - {{\widetilde{w}}}_{\rm f}^{(0)}$  and
by $\mathring{d}_{\rm s}^{(k)} =
d_{\rm s}^{(k)} - d_{\rm s}^{(0)}$
the ``lifted'' snapshots for the fluid  and the solid domain.
If we  denote by  $\{\lambda_j^{\rm f} \}_{j=1}^{n_{\rm train}}$ the eigenvalues  of the Gramian matrix associated with the fluid lifted states, we can  express the energy criterion as follows:
\begin{equation}
\label{eq:energy_criterion}
n_{\rm f} = \left\{
n \, : \,
\sum_{i=1}^n \lambda_i^{\rm f} \geq \left(1 - {\rm tol}_{\rm pod} \right)
\sum_{i=1}^{{n}_{\rm train}}\lambda_i^{\rm f}
\right\}.
\end{equation}
The same criterion also applies to the solid state and to the control. We rely on 
a weighted $H^1-L^2$ inner product for the fluid state, and on the $H^1$ inner product for the solid state and the control,
\begin{equation}
\label{eq:inner_products}
(\widetilde{w},\widetilde{w}')_{\rm f} := 
\frac{1}{u_{\infty}^2} (\widetilde{u},\widetilde{u}')_{H^1(\widetilde{\Omega}_{\rm f})}
+
\frac{1}{(\rho_{\rm f}u_{\infty}^2)^2} (\widetilde{p},\widetilde{p}')_{L^2(\widetilde{\Omega}_{\rm f})},
\quad
(d,d')_{\rm s} := (d,d')_{H^1(\widetilde{\Omega}_{\rm s})}, 
\quad
(g,g')_{\rm c} := (g,g')_{H^1(\widetilde{\Gamma})}, 
\end{equation}
where $u_{\infty}$ is a characteristic velocity of the flow. The choice of the inner product for the fluid is motivated by the fact that the flow velocity and the flow pressure might have different magnitude and also different units.

The application of 
POD leads to  three separate spaces that are built independently; 
as discussed in 
\cite[section 4.4]{Taddei2024optimization}, it is important to enrich the state spaces to ensure that the sensitivity matrix is full-rank. Towards this end, we define the perturbed snapshots $\{ {w}_{k,j}^{\rm f,en} \}_{k,j}$ such that
\begin{equation}
\label{eq:enrichment_step}
J_{\rm f,w}^{(k)} \left[ w_{\rm f}^{(k)},
d_{\rm s}^{(k)}
 \right] ( {w}_{k,j}^{\rm f,en} , w) +
 E_{\rm f}(\zeta_j^{\rm c}, v) = 0 , \quad \forall \, 
 w={\rm vec}(v,q) \in \mathcal{V}_{{\rm f},k}^0 \times  \mathcal{Q}_{{\rm f},k},
 \quad
k\in \texttt{I}_{\rm smp}, 
j=1,\ldots,n_{\rm c};
\end{equation}
 where $J_{\rm f,w}^{(k)} $ denotes the combined velocity and pressure Jacobian at time $t^{(k)}$.
The evaluation of $J_{\rm f,w}^{(k)} \left[ w_{\rm f}^{(k)}  \right]$ requires to estimate the time derivative of the solid velocity field (which enters in the calculation of the ALE velocity): if the sampling time step is significantly larger than the characteristic time scales of the fluid velocity, estimation of the flow acceleration   might be challenging; here, we do not address this issue. Then, we define the enriching modes for the fluid state using POD,
 \begin{equation}
 \label{eq:enriching_modes}
 \left\{ \widetilde{\zeta}_{i}^{\rm f,en}  \right\}_{j=1}^{n_{\rm f,en}} = {\rm POD} \left(
\{  
\widetilde{w}_{k,j}^{\rm f,en}  - \Pi_{\widetilde{\mathcal{Z}}_{\rm f}^{\rm pod}}  \widetilde{w}_{k,j}^{\rm f,en}
\}_{k,j}, \| \cdot \|_{\rm f}, n_{\rm f}^{\rm en}
 \right),
 \end{equation}
where
$\widetilde{\mathcal{Z}}_{\rm f}^{\rm pod} = {\rm span}
\{ \widetilde{\zeta}_{i}^{\rm f,pod} \}_{i=1}^{n_{\rm f,pod}}$ is the space generated using POD, and 
$\Pi_{\widetilde{\mathcal{Z}}_{\rm f}^{\rm pod}}  : 
\widetilde{\mathcal{X}}_{\rm f}^0  
\to 
\widetilde{\mathcal{Z}}_{\rm f}^{\rm pod}
$ is the orthogonal projection operator. 
The number of   modes $n_{\rm f}^{\rm en}$ that is added to the POD basis is chosen using the criterion:
\begin{equation}
\label{eq:projection_criterion}
n_{\rm f}^{\rm en} = {\rm min} \left\{
n \, : \, \max_{k,j} \frac{\| \widetilde{w}_{k,j}^{\rm f,en}  
- 
\Pi_{\widetilde{\mathcal{Z}}_{\rm f}^{\rm pod}}  \widetilde{w}_{k,j}^{\rm f,en}
-
\Pi_{\widetilde{\mathcal{Z}}_{\rm f,en}^{\rm pod}(n)}  \widetilde{w}_{k,j}^{\rm f,en}
\|_{\rm f}}{\| \widetilde{w}_{k,j}^{\rm f,en}   \|_{\rm f}}
< {\rm tol}_{\rm en}
\right\}.    
\end{equation}
Criterion \eqref{eq:projection_criterion} enables a direct control of the in-sample projection 
error and is justified by the fact that the norm of the snapshots 
$\{  
\widetilde{w}_{k,j}  - \Pi_{\widetilde{\mathcal{Z}}_{\rm f}^{\rm pod}}  \widetilde{w}_{k,j}
\}_{k,j}$ might be small, which potentially makes the energy criterion  \eqref{eq:energy_criterion} ineffective.
In conclusion, the reduced space for the fluid is given by $\widetilde{\mathcal{Z}}_{\rm f} = \widetilde{\mathcal{Z}}_{\rm f}^{\rm pod} \bigoplus \widetilde{\mathcal{Z}}_{\rm f}^{\rm en}$.
The approach
outlined above
is applied 
\emph{as is}
to the solid state; we omit the details.

\subsection{Definition of the local reduced-order models}
\label{sec:local_roms}

We resort to Galerkin projection to determine a local ROM for the solid problem:
\begin{equation}
\label{eq:solidROM}
\widehat{\mathbf{R}}_{\rm s}^{(k+1)}
\left(
\widehat{\boldsymbol{\alpha}}_{\rm s}^{(k+1)}
\right)
+\widehat{\mathbf{E}}_{\rm s} 
\widehat{\boldsymbol{\beta}}^{(k+1)}
 = 0,
 \quad 
 {\rm where} \;\; 
 \widehat{\mathbf{R}}_{\rm s}^{(k+1)}
\left(
 \boldsymbol{\alpha} 
\right)
:=
\mathbf{Z}_{\rm s}^\top 
{\mathbf{R}}_{\rm s}^{(k+1)}
\left(
\mathbf{d}_{\rm s}^{(0)} +
\mathbf{Z}_{\rm s} \boldsymbol{\alpha} 
\right),
\;\;
\widehat{\mathbf{E}}_{\rm s}  = 
\mathbf{Z}_{\rm s}^\top 
\mathbf{E}_{\rm s}
\mathbf{Z}_{\rm c}.
\end{equation}
Note that \eqref{eq:solidROM} is a (possibly nonlinear) system with $n_{\rm s}+n_{\rm c}$ unknowns and $n_{\rm s}$ equations that implicitly defines 
$\widehat{\boldsymbol{\alpha}}_{\rm s}^{(k+1)}
$ as a function of 
$\widehat{\boldsymbol{\beta}}^{(k+1)}$.

For the fluid subdomain, we first introduce the
$j_{\rm es}$-dimensional
empirical test space $\widetilde{\mathcal{Y}}_{\rm f} = {\rm span} \{  \widetilde{\psi}_{j}^{\rm f} \}_{j=1}^{j_{\rm es}} 
\subset 
\widetilde{\mathcal{V}}_{\rm f}^0 \times
\widetilde{\mathcal{Q}}_{\rm f}
$ such that $j_{\rm es}\geq n_{\rm f}$ and 
$( \widetilde{\psi}_{j}^{\rm f} , \widetilde{\psi}_{j'}^{\rm f}    )_{\rm f} = \delta_{j,j'}$; then, we define the residual:
\begin{equation}
\label{eq:fluidROM}
\widehat{\mathbf{R}}_{\rm f}^{(k+1)}
\left(
 \boldsymbol{\alpha}_{\rm f},  \boldsymbol{\alpha}_{\rm s}
\right)
:=
\mathbf{Y}_{\rm f}^\top 
 {\mathbf{R}}_{\rm f}^{(k+1)}
\left(
\overline{\mathbf{w}}_{\rm f}^{(k)} +
\mathbf{Z}_{\rm f} \boldsymbol{\alpha}_{\rm f}  \; , \;
\overline{\mathbf{d}}_{\rm s}^{(k)} +
\mathbf{Z}_{\rm s} \boldsymbol{\alpha}_{\rm s} 
\right),
\quad
\widehat{\mathbf{E}}_{\rm f} 
:=
\mathbf{Y}_{\rm f}^\top 
\mathbf{E}_{\rm f}
\mathbf{Z}_{\rm c};
\end{equation}
with $\mathbf{Y}_{\rm f} = \left[\boldsymbol{\psi}_1^{\rm f}, \ldots, 
\boldsymbol{\psi}_{j_{\rm es}}^{\rm f}\right]$.
Following \cite{taddei2021discretize}, we exploit POD to generate the space $\widetilde{\mathcal{Y}}_{j_{\rm es}}^{\rm f}$. First, we compute the Riesz elements:
$$
\left(
\widetilde{\psi}_{k,i}, \widetilde{z}
\right)_{\rm f}
= J_{\rm f,w}^{(k)}\left[
w_{\rm f}^{(k)},
d_{\rm s}^{(k)}
\right]
\left(
\zeta_i^{\rm f}, z
\right),
\quad
\forall \, \widetilde{z}\in \widetilde{\mathcal{X}}_{\rm f}^0,
\;\;
i=1,\ldots,n_{\rm f}, \;\;
k\in \texttt{I}_{\rm smp}.
$$
Then, we apply POD to determine the reduced space
$$
\left\{ \widetilde{\psi}_{j}^{\rm f}  \right\}_{j=1}^{j_{\rm es}} = {\rm POD} \left(
\{  
\widetilde{\psi}_{k,i}
\}_{k,i}, \| \cdot \|_{\rm f}, j_{\rm es}
 \right),
$$
where the integer $j_{\rm es}$ is set equal to $2n_{\rm f}$ in the numerical experiments.

Note that the 
residual in \eqref{eq:solidROM} is a reduced-order counterpart of the full-order residual in \eqref{eq:weak_formulation_discrete}$_3$: it can hence be included in  \eqref{eq:optimal_control_algebraic} without changing the structure of the problem.
On the other hand, the 
 residual in \eqref{eq:fluidROM} leads to an underdetermined system with $j_{\rm es}$ equations and $n_{\rm f}+n_{\rm c}$ unknowns.  In the next section, we discuss how to determine an actionable ROM based on \eqref{eq:fluidROM}. 
We further observe that the local ROMs associated with  \eqref{eq:solidROM} and \eqref{eq:fluidROM} require hyper-reduction to achieve significant speedups. As discussed in \cite{iollo2023one}, standard hyper-reduction techniques for monolithic ROMs can be readily applied to \eqref{eq:solidROM} and \eqref{eq:fluidROM} ; in this work, we do not address this issue.
 
\subsection{Solution to the global reduced-order model}
\label{sec:global_rom}

Exploiting the expressions in 
\eqref{eq:optimal_control_algebraic} and 
\eqref{eq:ansatz} and the approximation
$\texttt{D}_{{\rm s}, \Delta t} \widehat{\boldsymbol{\alpha}}_{\rm s}^{(k+1)} = \alpha_{\rm s,v}^{(k+1)} \widehat{\boldsymbol{\alpha}}_{\rm s}^{(k+1)} + \widehat{\mathbf{b}}_{\rm s,v}$, we define the objective 
\begin{equation}
\label{eq:reduced_objective}
\mathfrak{f}^{(k+1)}( \boldsymbol{\alpha}_{\rm f},\boldsymbol{\alpha}_{\rm s},\boldsymbol{\beta})
=
\frac{1}{2} \Big|
\widehat{\mathbf{P}}_{\rm f} \boldsymbol{\alpha}_{\rm f}
-
\alpha_{\rm s,v}^{(k+1)} 
\widehat{\mathbf{P}}_{\rm s}  
\widehat{\boldsymbol{\alpha}}_{\rm s}^{(k+1)}  
- 
\widehat{\mathbf{b}}
\Big|^2
\, + \,
\frac{\delta}{2}
\vertiii{\mathbf{Z}_{\rm c} \boldsymbol{\beta}}^2,
\end{equation}
with $\widehat{\mathbf{P}}_{\rm f}  = 
{\mathbf{P}}_{\rm f} \mathbf{Z}_{\rm f}(1:
N_{\rm u}, :)$, $\widehat{\mathbf{P}}_{\rm s}  = 
{\mathbf{P}}_{\rm s} \mathbf{Z}_{\rm s}$,
and
$\widehat{\mathbf{b}} = 
-\widehat{\mathbf{P}}_{\rm f} \mathbf{w}_{\rm f}^{(0)}(1:N_{\rm u})
+ {\mathbf{P}}_{\rm s} 
\left(
\mathbf{d}_{\rm s}^{(0)} \, + \, 
\mathbf{Z}_{\rm s} 
\widehat{\mathbf{b}}_{\rm s,v} \right)$.
In view of the definition of the partitioned ROM, we introduce the LSPG ROM for $\widehat{\boldsymbol{\alpha}}_{\rm f}^{(k+1)}$, given $\widehat{\boldsymbol{\alpha}}_{\rm s}^{(k+1)}$ and $\widehat{\boldsymbol{\beta}}^{(k+1)}$:
$$
\widehat{\boldsymbol{\alpha}}_{\rm f}^{(k+1)} = 
{\rm arg} \min_{\boldsymbol{\alpha}\in \mathbb{R}^{n_{\rm f}}} \Big|
\widehat{\mathbf{R}}_{\rm f}^{(k+1)}
\left(
 \boldsymbol{\alpha}_{\rm f},  \widehat{\boldsymbol{\alpha}}_{\rm s}^{(k+1)}
\right)
+ \widehat{\mathbf{E}}_{\rm f} 
\widehat{\boldsymbol{\beta}}^{(k+1)}
\Big|.
$$
The latter  implies that $\widehat{\boldsymbol{\alpha}}_{\rm f}^{(k+1)}$ solves
\begin{equation}
\label{eq:LSPG_exact}
    \left(
\widehat{\mathbf{J}}_{\rm f,f}^{(k+1)}
    \right)
^\top
    \left(
\widehat{\mathbf{R}}_{\rm f}^{(k+1)}
\left(
\widehat{\boldsymbol{\alpha}}_{\rm f}^{(k+1)},  \widehat{\boldsymbol{\alpha}}_{\rm s}^{(k+1)}
\right)
+ \widehat{\mathbf{E}}_{\rm f} 
\widehat{\boldsymbol{\beta}}^{(k+1)}
    \right) = 0,
    \quad
    {\rm where} \;\;
 \widehat{\mathbf{J}}_{\rm f,f}^{(k+1)}
:=
\frac{\partial \widehat{\mathbf{R}}_{\rm f}^{(k+1)}   }{\partial \boldsymbol{\alpha}_{\rm f}}
\left(
\widehat{\boldsymbol{\alpha}}_{\rm f}^{(k+1)},  \widehat{\boldsymbol{\alpha}}_{\rm s}^{(k+1)}
\right).
\end{equation}
Linearization of \eqref{eq:LSPG_exact}
requires the estimation of the Hessian of the reduced residual of the fluid problem, which is expensive to compute; for this reason, we discuss two approximate strategies that 
circumvent the need for Hessian estimation.
The two methods show similar performance for the numerical experiments of the present work.

First, we consider the Petrov-Galerkin formulation:
\begin{equation}
\label{eq:global_rom_strategy1}
\min_{\boldsymbol{\alpha}_{\rm f}\in \mathbb{R}^{n_{\rm f}}, \boldsymbol{\alpha}_{\rm s}
\in \mathbb{R}^{n_{\rm s}},
 \boldsymbol{\beta} \in \mathbb{R}^{n_{\rm c}}} \;\; 
\mathfrak{f}^{(k+1)}( \boldsymbol{\alpha}_{\rm f},\boldsymbol{\alpha}_{\rm s},\boldsymbol{\beta}),
\qquad
{\rm s.t.} \;\;
\left\{
\begin{array}{l}
\displaystyle{
   \left(
 \widehat{\mathbf{J}}_{\rm f,f}^{(k)}
    \right)
^\top
    \left(
\widehat{\mathbf{R}}_{\rm f}^{(k+1)}
\left(
 \boldsymbol{\alpha}_{\rm f},   \boldsymbol{\alpha}_{\rm s}
\right)
+ \widehat{\mathbf{E}}_{\rm f} 
\boldsymbol{\beta}
    \right) = 0
} ; \\[3mm]
\displaystyle{
 \widehat{\mathbf{R}}_{\rm s}^{(k+1)}
\left(
 \boldsymbol{\alpha}_{\rm s} 
\right)
+\widehat{\mathbf{E}}_{\rm s} 
\boldsymbol{\beta}
 = 0
};
\\
\end{array}
\right.
\end{equation}
where $ 
 \widehat{\mathbf{J}}_{\rm f,f}^{(k)}
:=
\frac{\partial \widehat{\mathbf{R}}_{\rm f}^{(k)}   }{\partial \boldsymbol{\alpha}_{\rm f}}
\left(
\widehat{\boldsymbol{\alpha}}_{\rm f}^{(k)},  \widehat{\boldsymbol{\alpha}}_{\rm s}^{(k)}
\right) \in \mathbb{R}^{j_{\rm es} \times n_{\rm f}}$ depends on the generalized coordinates from the previous time step.
Note that \eqref{eq:global_rom_strategy1} can be solved using the same SQP strategy outlined in section \ref{sec:optimal_control_formulation} for the full-order formulation 
\eqref{eq:optimal_control_algebraic}.

The second approach, which was considered in
\cite{Taddei2024optimization},
consists in finding an approximate solution to 
\begin{equation}
\label{eq:global_rom_strategy2}
\min_{\boldsymbol{\alpha}_{\rm f}\in \mathbb{R}^{n_{\rm f}}, \boldsymbol{\alpha}_{\rm s}
in \mathbb{R}^{n_{\rm s}},
 \boldsymbol{\beta} \in \mathbb{R}^{n_{\rm c}}} \;\; 
\mathfrak{f}^{(k+1)}( \boldsymbol{\alpha}_{\rm f},\boldsymbol{\alpha}_{\rm s},\boldsymbol{\beta}),
\qquad
{\rm s.t.} \;\;
\left\{
\begin{array}{l}
\displaystyle{
   \left(
 \widehat{\mathbf{J}}_{\rm f,f}^{(k+1)}
    \right)^\top
    \left(
\widehat{\mathbf{R}}_{\rm f}^{(k+1)}
\left(
 \boldsymbol{\alpha}_{\rm f},   \boldsymbol{\alpha}_{\rm s}
\right)
+ \widehat{\mathbf{E}}_{\rm f} 
\boldsymbol{\beta}
    \right) = 0
} ; \\[3mm]
\displaystyle{
 \widehat{\mathbf{R}}_{\rm s}^{(k+1)}
\left(
 \boldsymbol{\alpha}_{\rm s} 
\right)
+\widehat{\mathbf{E}}_{\rm s} 
\boldsymbol{\beta}
 = 0
};
\\
\end{array}
\right.
\end{equation}
using the SQP method with an inexact Jacobian that neglects the 
gradient of $\widehat{\mathbf{J}}_{\rm f,f}^{(k+1)}$, which depends on the Hessian of the residual. In more detail,
at each subiteration $j$ of the method,
if we denote by
$\left(
\widehat{\boldsymbol{\alpha}}_{\rm f}^{(k+1,j)},
\widehat{\boldsymbol{\alpha}}_{\rm s}^{(k+1,j)},
\widehat{\boldsymbol{\beta}}^{(k+1,j)} \right)$ the estimate of the solution at time $t^{(k+1)}$,
 we solve \eqref{eq:global_rom_strategy2} with the linearized constraints:
\begin{equation}
\label{eq:SQP_reduced_strategy2}
\left\{
\begin{array}{l}
\displaystyle{
\left(
\widehat{\mathbf{J}}_{\rm f,f}^{(k+1,j)}
\right)^\top
\left(
\widehat{\mathbf{R}}_{\rm f}^{(k+1,j)}
+
\widehat{\mathbf{J}}_{\rm f,f}^{(k+1,j)}
\left( \boldsymbol{\alpha}_{\rm f} - 
\widehat{\boldsymbol{\alpha}}_{\rm f}^{(k+1,j)} \right)
+
\widehat{\mathbf{J}}_{\rm f,s}^{(k+1,j)}
\left( \boldsymbol{\alpha}_{\rm s} - 
\widehat{\boldsymbol{\alpha}}_{\rm s}^{(k+1,j)} \right)
+ 
\widehat{\mathbf{E}}_{\rm f} \boldsymbol{\beta}
\right)
 =  
0} ,\\[3mm]
\displaystyle{
\widehat{\mathbf{R}}_{\rm s}^{(k+1,j)}
+ 
\widehat{\mathbf{J}}_{\rm s,s}^{(k+1,j)}
\left( \boldsymbol{\alpha}_{\rm s} - 
\widehat{\boldsymbol{\alpha}}_{\rm s}^{(k+1,j)} \right)
+\widehat{\mathbf{E}}_{\rm s} 
\boldsymbol{\beta}
 =   0
},
\\
\end{array}
\right.    
\end{equation}
 where
 $ \widehat{\mathbf{J}}_{\rm f,f}^{(k+1,j)}
:=
\frac{\partial \widehat{\mathbf{R}}_{\rm f}^{(k+1)}   }{\partial \boldsymbol{\alpha}_{\rm f}}
\left(
\widehat{\boldsymbol{\alpha}}_{\rm f}^{(k+1,j)},  \widehat{\boldsymbol{\alpha}}_{\rm s}^{(k+1,j)}
\right) \in \mathbb{R}^{j_{\rm es} \times n_{\rm f}}$.

\section{Numerical results}
\label{sec:numerical_results}

We apply our method to three model problems: 
the simplified model problem of \cite{Astorino_Grandmont_2010} with known exact solution;
the interaction of a Newtonian  fluid with an elastic beam; and 
the benchmark  presented in  \cite{TurekHron2006} that features the  interaction of a Newtonian fluid  with an hyper-elastic body modeled by the Saint-Venant Kirchhoff law.
We assume  a $\mathbb{P}_2-\mathbb{P}_1$ Taylor-Hood discretization for the fluid and a $\mathbb{P}_2$ discretization for the solid,
 the Newmark time integration scheme with 
 $\gamma=1/2$ and $\beta=1/4$ in \eqref{eq:newmark} and the BDF2 method in  \eqref{eq:BDF}.
 Furthermore, we consider ${\rm tol}_{\rm en}=0.1$ in \eqref{eq:projection_criterion}.
 We denote by $\Delta t$ the time step size.
 
\subsection{Model problem with exact solution}
\label{sec:astorino_grandmont}

 We consider the configuration in Figure \ref{fig:FSI_illustration} with 
$\widetilde{\Omega}_{\rm f}=(0,\,1)\times (0,1)$ and $\widetilde{\Omega}_{\rm s}=(0,\,1)\times (1,\,1.25)$;
we neglect the deformation of the structure --- that is, we consider ${\Omega}_{\rm f}(t)\equiv \widetilde{\Omega}_{\rm f}$ for all times --- but we impose the transmission conditions. Following \cite{Astorino_Grandmont_2010},
we  consider the Newtonian model  \eqref{eq:newtonian_model}
and the linear model \eqref{eq:solid_model} for the structure with the parameters
$$
\rho_{\rm f} = 1, \quad
\mu_{\rm f} = 0.013, \quad 
\rho_{\rm s} = 1.9,  \quad
\mu_{\rm s} = 3, \quad \lambda_{\rm s} = 4.5.
$$
We prescribe Neumann boundary conditions on the left and the right boundary and Dirichlet conditions on the upper boundary of the solid domain;
we prescribe Neumann boundary conditions on the left and  right boundaries and Dirichlet conditions on the lower boundary of the fluid domain. The boundary and the initial conditions and the source terms are chosen to ensure that the exact solution is
$$
\begin{array}{l}
u_{\rm f}(x,y,t)
=\left[
\begin{array}{l}
\cos(x+t)\sin(y+t) + \sin(x+t)\cos(y+t)\\
-\sin(x+t)\cos(y+t) - \cos(x+t)\sin(y+t) \\
\end{array}
\right],\quad
d_{\rm s}(x,y,t)
=\left[
\begin{array}{l}
\sin(x+t)\sin(y+t)\\
\cos(x+t)\cos(y+t)\\
\end{array}
\right],
\\[5mm]
p_{\rm f}(x,y,t)
=
2\mu_{\rm f}(\sin(x+t)\sin(y+t) - \cos(x+t)\cos(y+t))+2\mu_{\rm s}\cos(x+t)\sin(y+t).
\\
\end{array}
$$
We integrate the coupled system in the time interval $(0,1)$, and we consider structured equispaced grids with characteristic spatial size  $h$ for both fluid and solid domains.

Table \ref{tab:astorino_table} illustrates the performance of the optimization-based approach for several spatial and time discretizations. In more details, we report the integral errors
\begin{equation}
\label{eq:error_metric_model1}
E_{\rm f}(h,\Delta t) :=
\sqrt{ 
\frac{1}{K} \sum_{k=1}^K  
\big\|
u_{\rm f}^{{\rm hf}, (k)}
-
u_{\rm f}^{{\rm true}, (k)}
\big\|_{H^1(\widetilde{\Omega}_{\rm f})}^2
},
\quad
E_{\rm s}(h,\Delta t) :=
\sqrt{ 
\frac{1}{K} \sum_{k=1}^K  
\big\|
d_{\rm s}^{{\rm hf}, (k)}
-
d_{\rm s}^{{\rm true}, (k)}
\big\|_{\rm s}^2
};    
\end{equation}
and we estimate the convergence rate as
$\eta_{\rm f}: =  -  \log_2 \left(  \frac{E_{\rm f}(h/2,\Delta t/2) }{E_{\rm f}(h,\Delta t)}    \right)  $ and
$\eta_{\rm s}: =  -  \log_2 \left(  \frac{E_{\rm s}(h/2,\Delta t/2) }{E_{\rm s}(h,\Delta t)}    \right)  $.
We empirically find that the global FOM is second-order accurate in space and time for BDF2 time integration, and first-order accurate for BDF1 time integration.

\begin{table}[H]
\begin{tabular}{|c|c|c|c|c|c|c|c|c|}
\hline
$(h,\Delta t)$ & 
\multicolumn{4}{|c|}{BDF1}
 &  
\multicolumn{4}{|c|}{BDF2} \\ \hline
 & error (solid) &    rate & error (fluid) &    rate & 
   error (solid) &    rate & error (fluid) &   rate  
   \\[3mm]
\hline
 $\big( \frac{1}{15},  0.01 \big)$  & $3.98\text{e}{-}4$ &  & $7.05\text{e}{-}2$ & & $8.28\text{e}{-}5$ &  & $5.54\text{e}{-}3$ &  \\[3mm]
  $\big( \frac{1}{30},  0.005 \big)$ & $1.98\text{e}{-}4$ & 1.01 & $3.51\text{e}{-}2$& 1.01 & $1.70\text{e}{-}5$ & 2.28 & $1.00\text{e}{-}3$ & 2.47 \\[3mm]
   $\big( \frac{1}{60},  0.0025 \big)$ &  $9.97\text{e}{-}5$ & 0.99 & $1.76\text{e}{-}2$ & 1.00 & $3.78\text{e}{-}6$ & 2.17 & $1.96\text{e}{-}4$ & 2.35 \\[3mm]
    $\big( \frac{1}{120},  0.00125 \big)$  & $5.01\text{e}{-}5$ & 0.99 & $8.78\text{e}{-}3$ & 1.00 & $8.85\text{e}{-}7$ & 2.09 & $4.13\text{e}{-}5$ & 2.25 \\[3mm]
 \hline
\end{tabular}

\caption{model problem with exact solution. Accuracy of the optimization-based partitioned method for multiple discretizations of varying sizes. The error metrics are introduced in \eqref{eq:error_metric_model1}.}
\label{tab:astorino_table}
\end{table}

\subsection{Elastic beam}
\label{sec:elastic_beam}
\begin{figure}[H]
\centering

\subfloat[]{\includegraphics[width=0.45\textwidth]{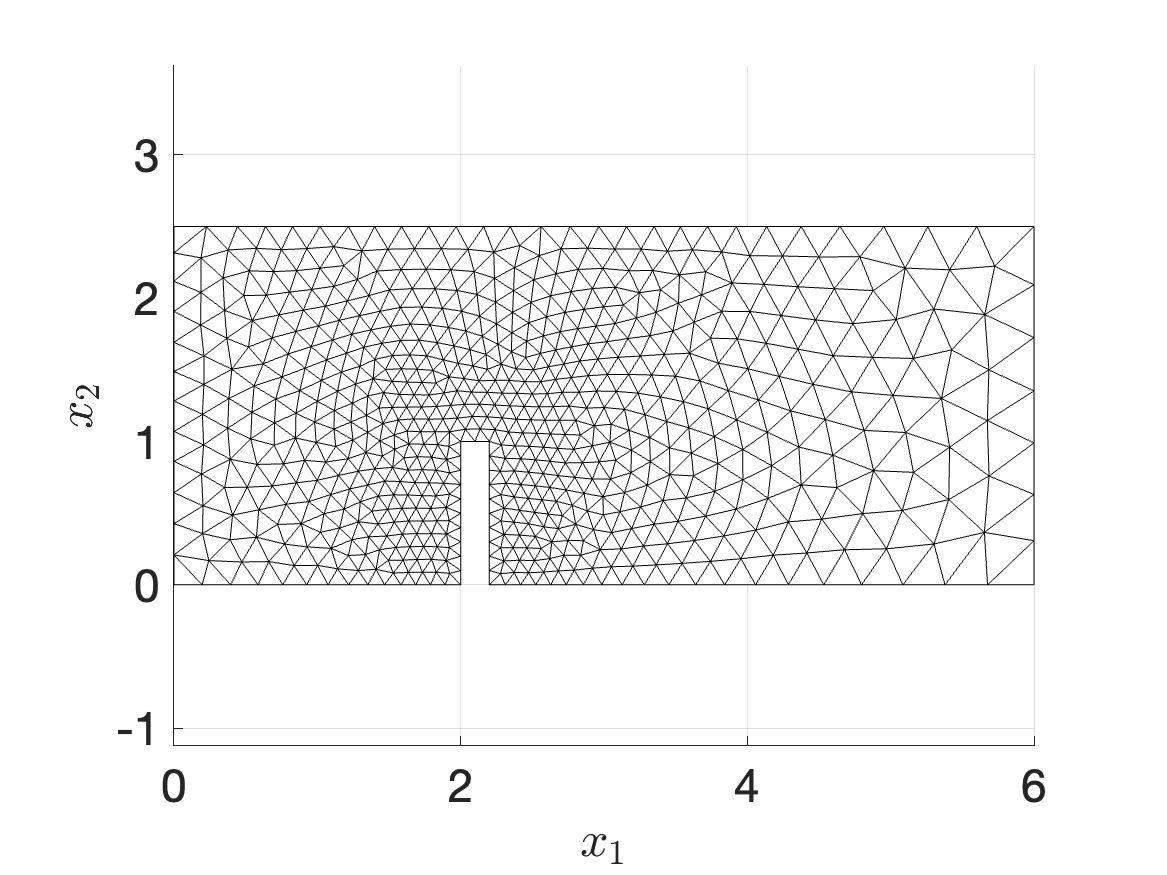}}
~~
\subfloat[]{\includegraphics[width=0.45\textwidth]{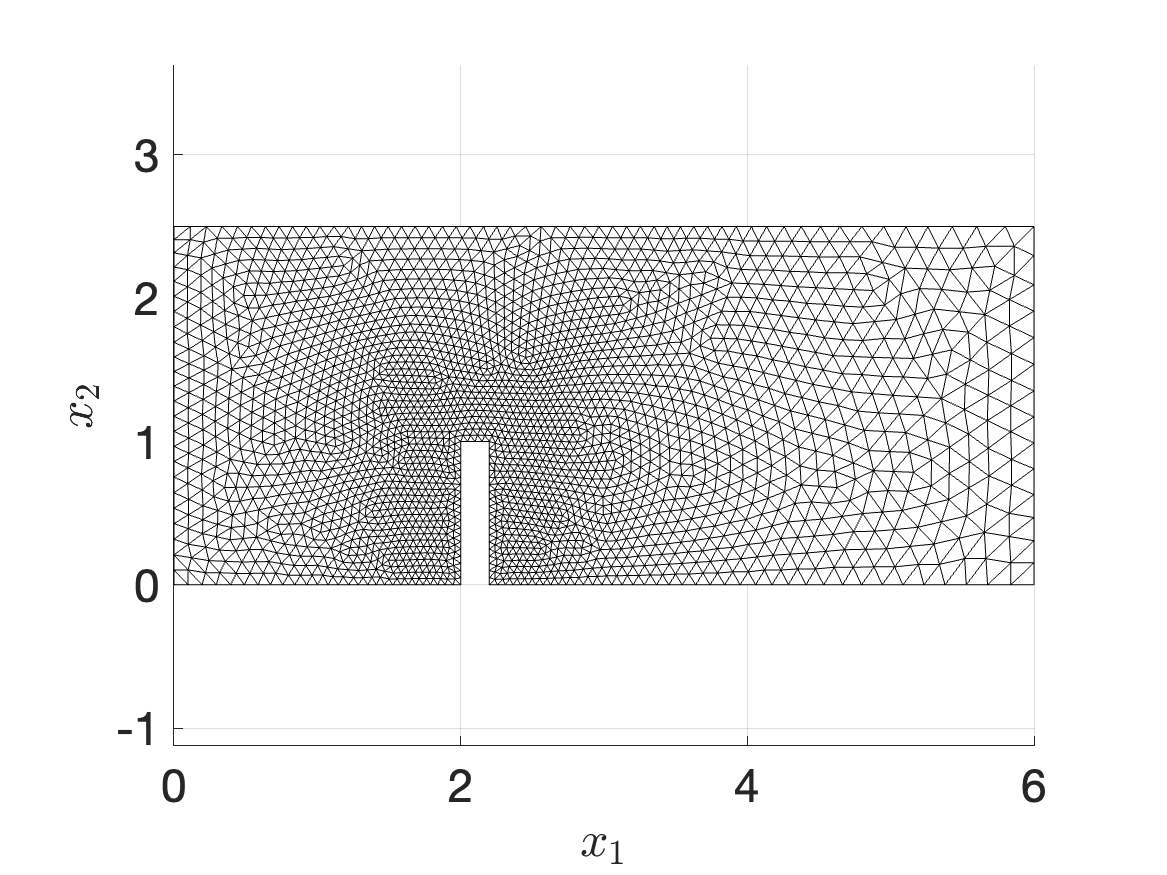}}
\caption{elastic beam; computational meshes.}
\label{fig:vbeam_meshes}
\end{figure}

\begin{figure}[H]
\centering

\subfloat[$t=2$]{\includegraphics[width=0.33\textwidth]{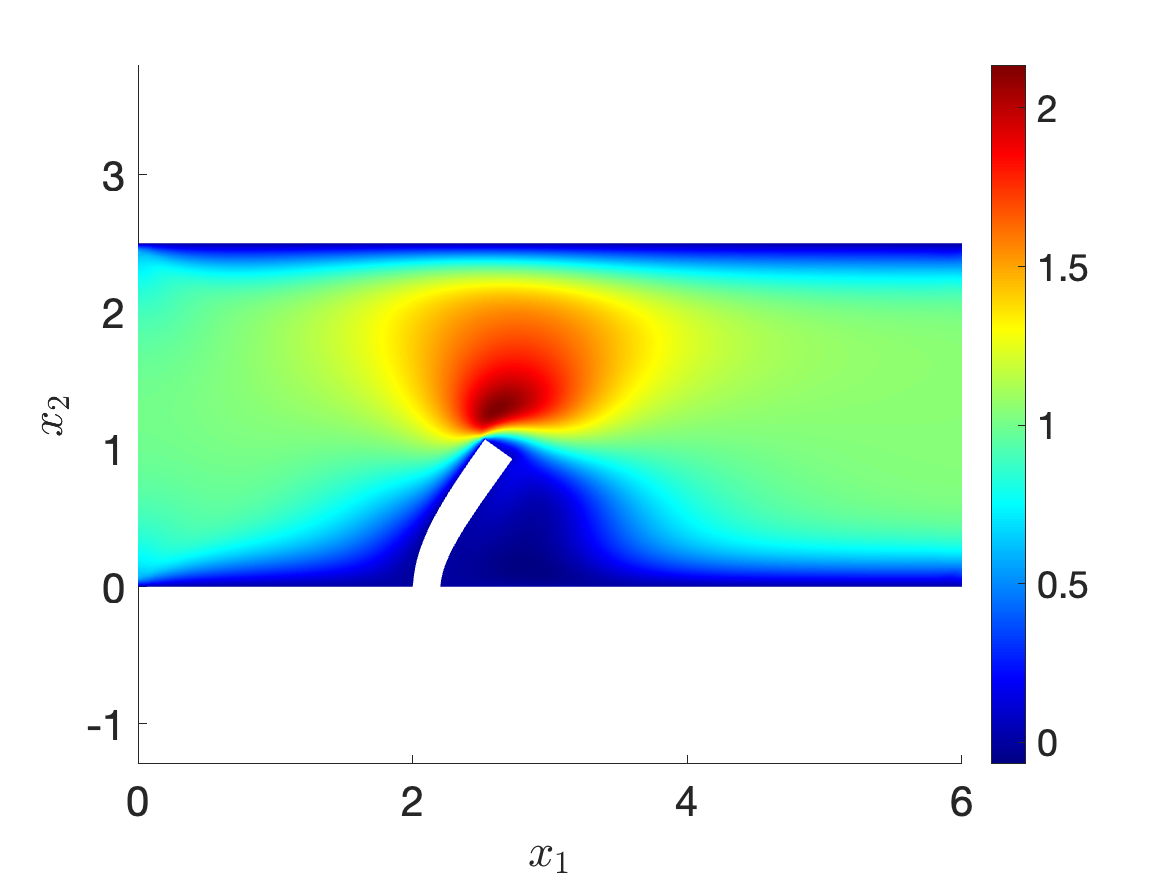}}
~~
\subfloat[$t=4$]{\includegraphics[width=0.33\textwidth]{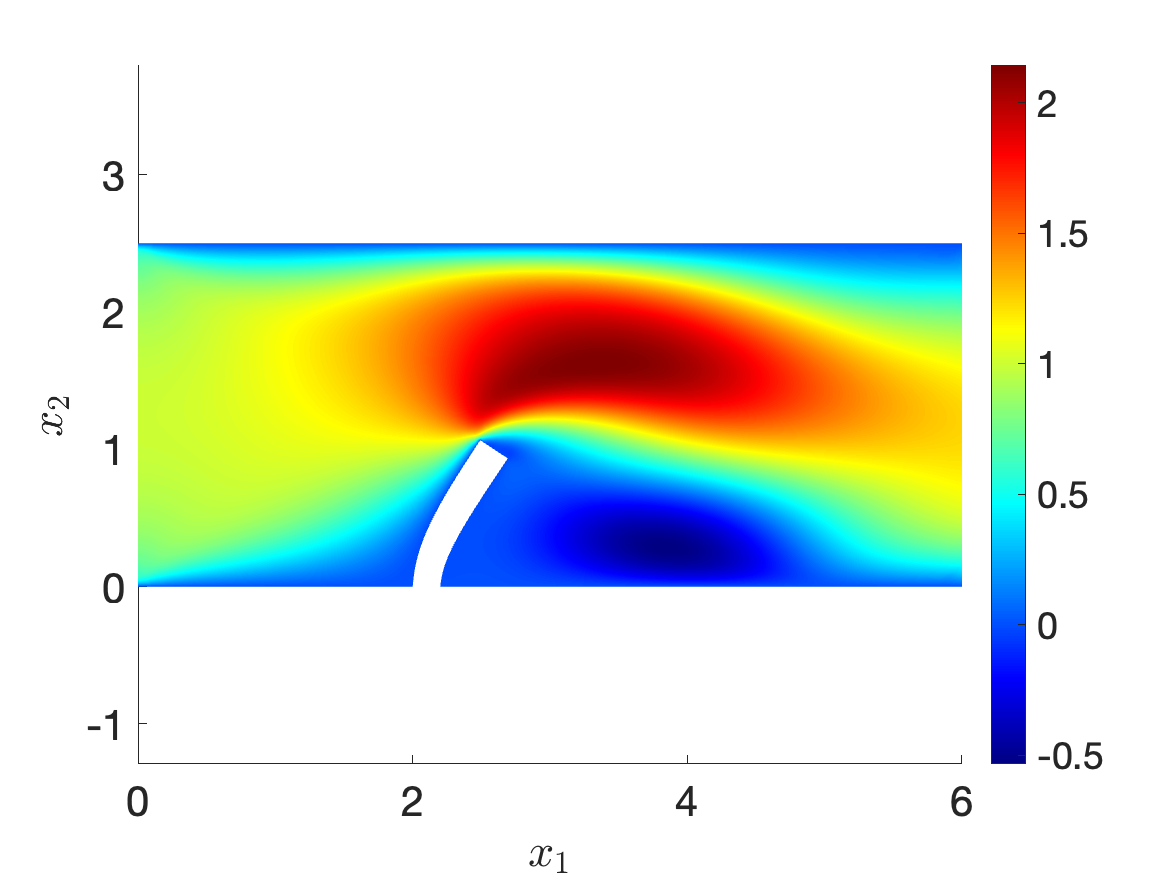}}
~~
\subfloat[$t=6$]{\includegraphics[width=0.33\textwidth]{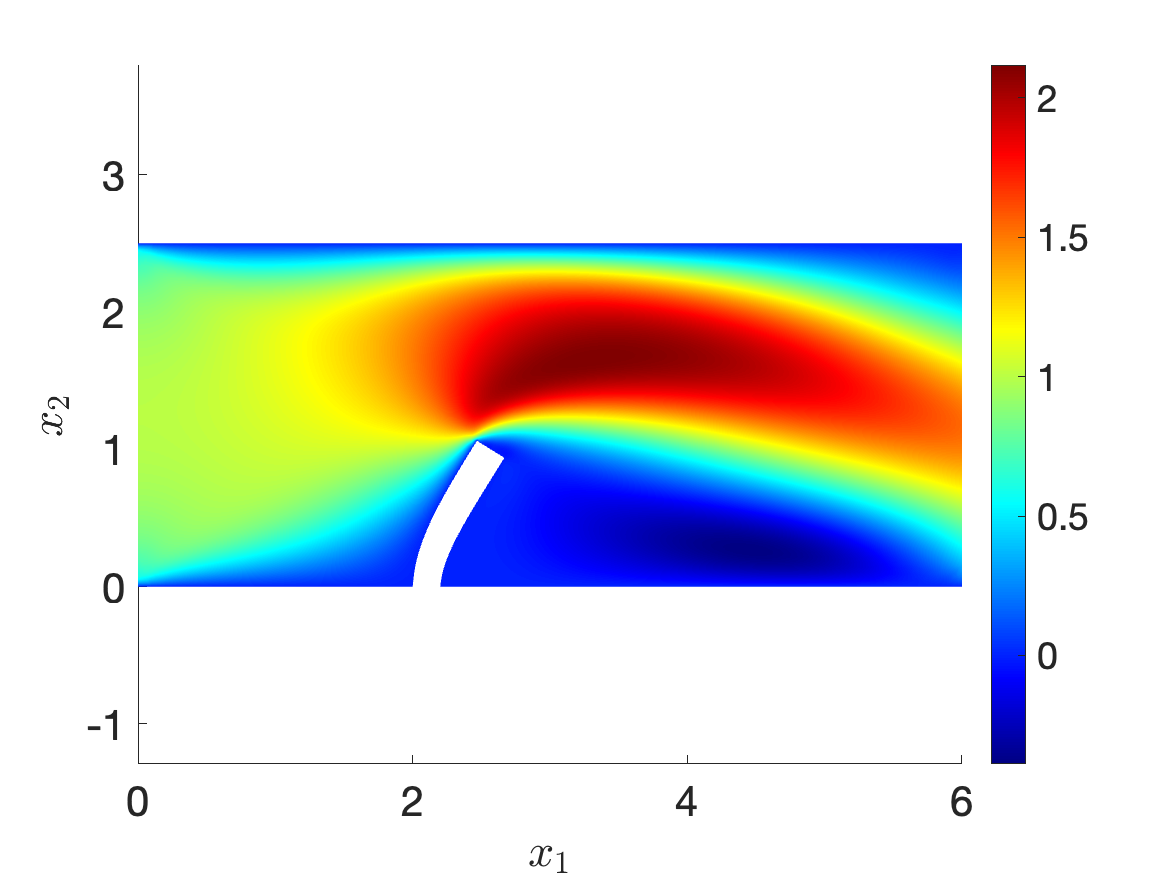}}
\caption{elastic beam; behavior of the streamwise velocity for three time instants.}
\label{fig:vbeam_velocity}
\end{figure}

We consider the interaction of a Newtonian fluid in the laminar regime with an elastic beam; similar test cases were considered in \cite{xiao2016non,nonino2023projection}. 
The computational domain is depicted in Figure \ref{fig:vbeam_meshes}.
We prescribe wall boundary conditions at the top and the bottom boundaries; we impose a parabolic velocity profile at the inlet with peak equal to $u_{\infty}=1$,
$$
u_{\rm f}^{\rm dir}(x,t) = \left\{
\begin{array}{ll}
\displaystyle{c_{\rm f}^{\rm dir}(t)
\frac{4u_{\infty}}{H^2}
{\rm vec} \left( x_2(H-x_2), 0  \right)   }  &  {\rm if} \, x_1=0\\
 0    &  {\rm otherwise};\\
\end{array}
\right.
\quad
{\rm where} \;\;
c_{\rm f}^{\rm dir}(t) = 
\left\{
\begin{array}{ll}
\frac{1 - \cos(\pi t / 2)}{2}   &  {\rm if} \,t \in [0,2],\\
1   &   {\rm if} \, t >2;\\
\end{array}
\right.
$$
and $H=2.5$;
and homogeneous Neumann conditions at the outflow. We consider the material parameters
$$
\rho_{\rm f} = 1, \quad \rho_{\rm s}=1.1, \quad E=10^3, \quad
\nu=0.3, \quad \mu_{\rm f}= 0.035;
$$
and we integrate the system for $T=6$ non-dimensional units. We 
consider two discretizations:
the coarse discretization with $N_{\rm u} = 6550 $ , $N_{\rm s} = 210$  and
$N_{\rm c} = 90$ spatial degrees of freedom, and $\Delta t= 5\cdot 10^{-2}$; and a finer discretization with $N_{\rm u} = 25576$ , $N_{\rm s} = 738$  and
$N_{\rm c} = 178$ spatial degrees of freedom, and $\Delta t = 2 \cdot 10^{-2}$. 
The fluid meshes are reported in Figures \ref{fig:vbeam_meshes}(a) and (b). 
Figure \ref{fig:vbeam_velocity} shows the behavior of the streamwise (horizontal) velocity for three different time instants: the elastic beam gradually bends to reach its maximal deflection and then slightly recovers to reach a steady state.

Table \ref{tab:vbeam_table} assesses performance of the FOM solver based on SQP with exact and inexact Jacobian, and on Dirichlet-to-Neumann (DtN) iterations with Aitken's relaxation \cite{Kuttler2008,Habchi2013} --- for DtN, we set the initial relaxation parameter equal to $\theta_0=0.05$. 
For both methods, we consider the termination condition \eqref{eq:termination_condition} with ${\rm tol}_{\rm sqp}=10^{-6}$.
We notice that SQP with inexact Jacobian is superior to DtN for both discretizations, due to the massive reduction in the number of iterations.
Interestingly, the benefit of SQP is less significant as we increase the size of the mesh due to the need to compute the sensitivity matrices \eqref{eq:HF_sensitivity_matrix} (10x speedup for the coarse mesh vs 6.7x speedup for the fine mesh). We further notice that the computation of the exact Jacobian of the constraints, which includes the shape derivatives, is not beneficial in terms of computational cost, even if it significantly reduces the average number of iterations.
We conjecture that more sophisticated implementation  of the shape derivatives might reduce the gap between the two methods.

\begin{table}[h!]
\centering
\begin{tabular}{|l|c|c|c|c|}
\hline
  &  
\multicolumn{2}{|c|}{avg nbr its}
 &  
\multicolumn{2}{|c|}{comp cost [s]} \\
 &  mesh 1 & mesh 2 &     
   mesh 1 & mesh 2  \\[3mm]
\hline
SQP (exact Jacobian) & 3.23 & 3.00&  627 & 12550\\[3mm]
 \hline
 SQP (inexact Jacobian) & 5.16& 4.68 &  116 & 2150 \\[3mm]
 \hline
DtN & 22.24 & 27.64 & 1158 & 14351  \\[3mm]
 \hline
\end{tabular}

\caption{elastic beam. Comparison of SQP and DtN  solvers for two FE discretizations of the local subproblems.}
\label{tab:vbeam_table}
\end{table}

\begin{figure}[H]
\centering

\subfloat[]{\includegraphics[width=0.33\textwidth]{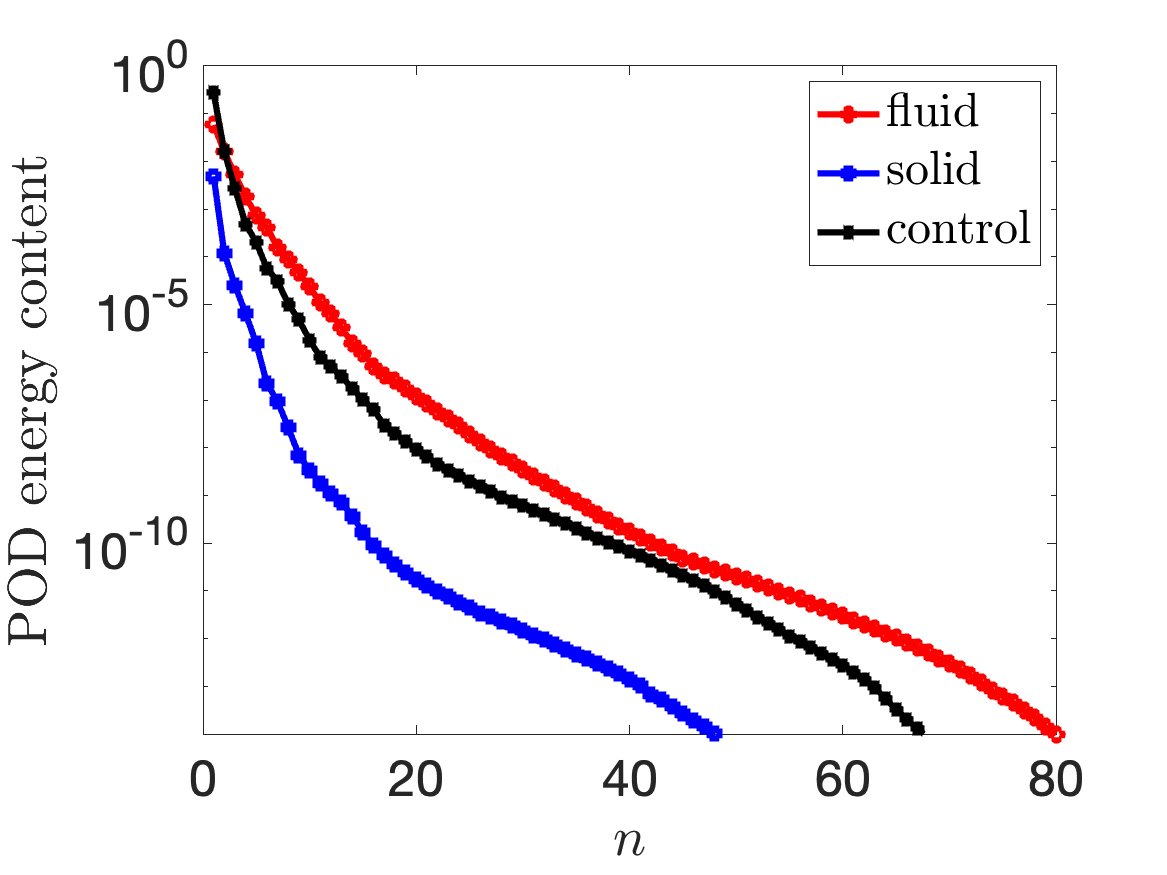}}
~~
\subfloat[]{\includegraphics[width=0.33\textwidth]{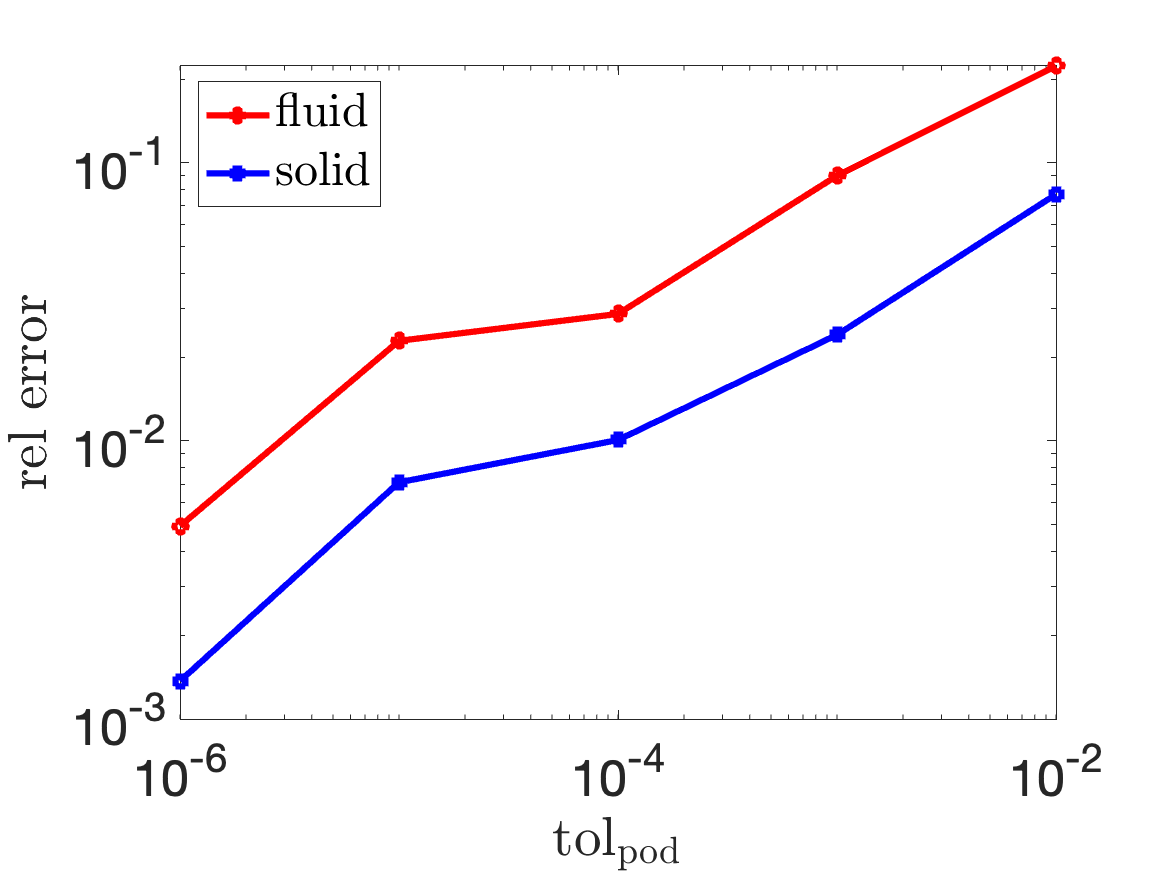}}
~~
\subfloat[]{\includegraphics[width=0.33\textwidth]{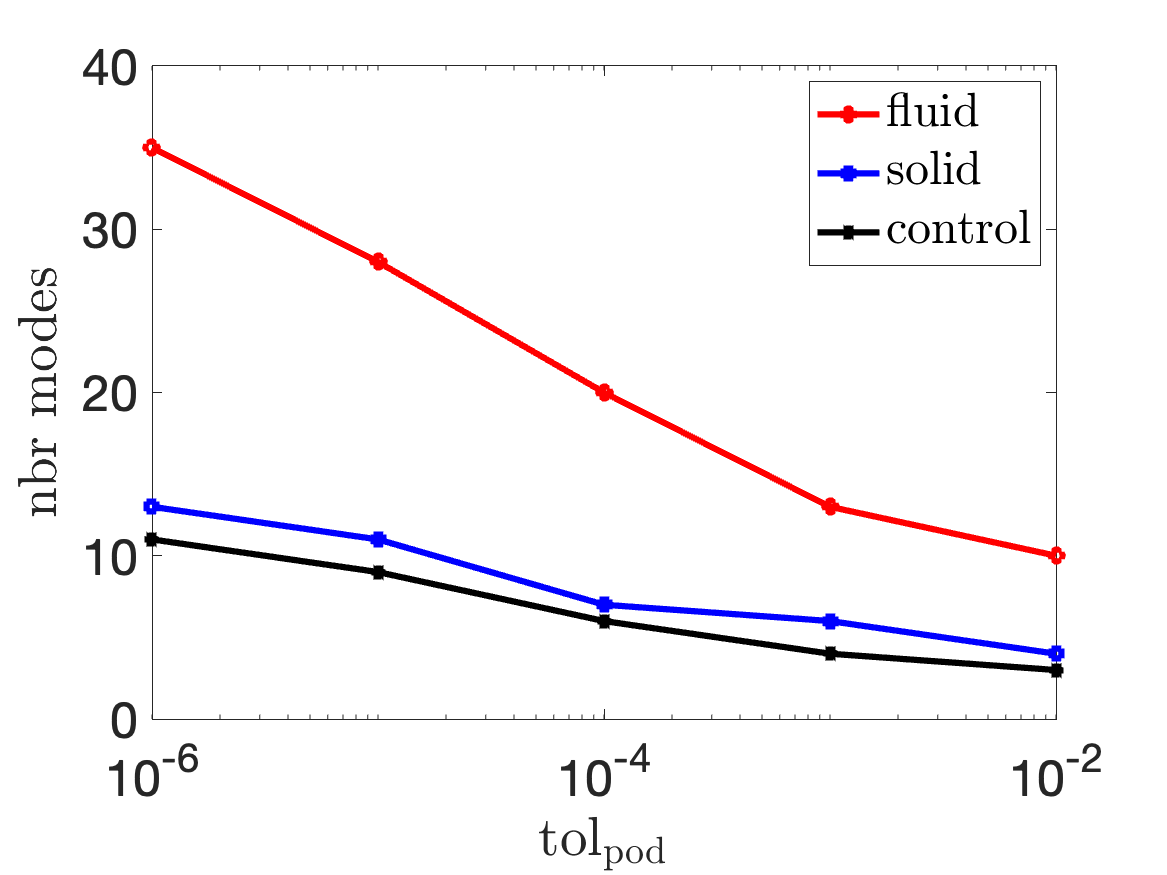}}
\caption{elastic beam; performance of the global ROM. (a) energy content of the discarded POD modes.
(b)-(c) relative error and number of modes for different choices of the tolerance ${\rm tol}_{\rm pod}$ with ${\rm tol}_{\rm en}=0.1$.}
\label{fig:vbeam_MOR1}
\end{figure}

\begin{figure}[H]
\centering
\subfloat[]{\includegraphics[width=0.33\textwidth]{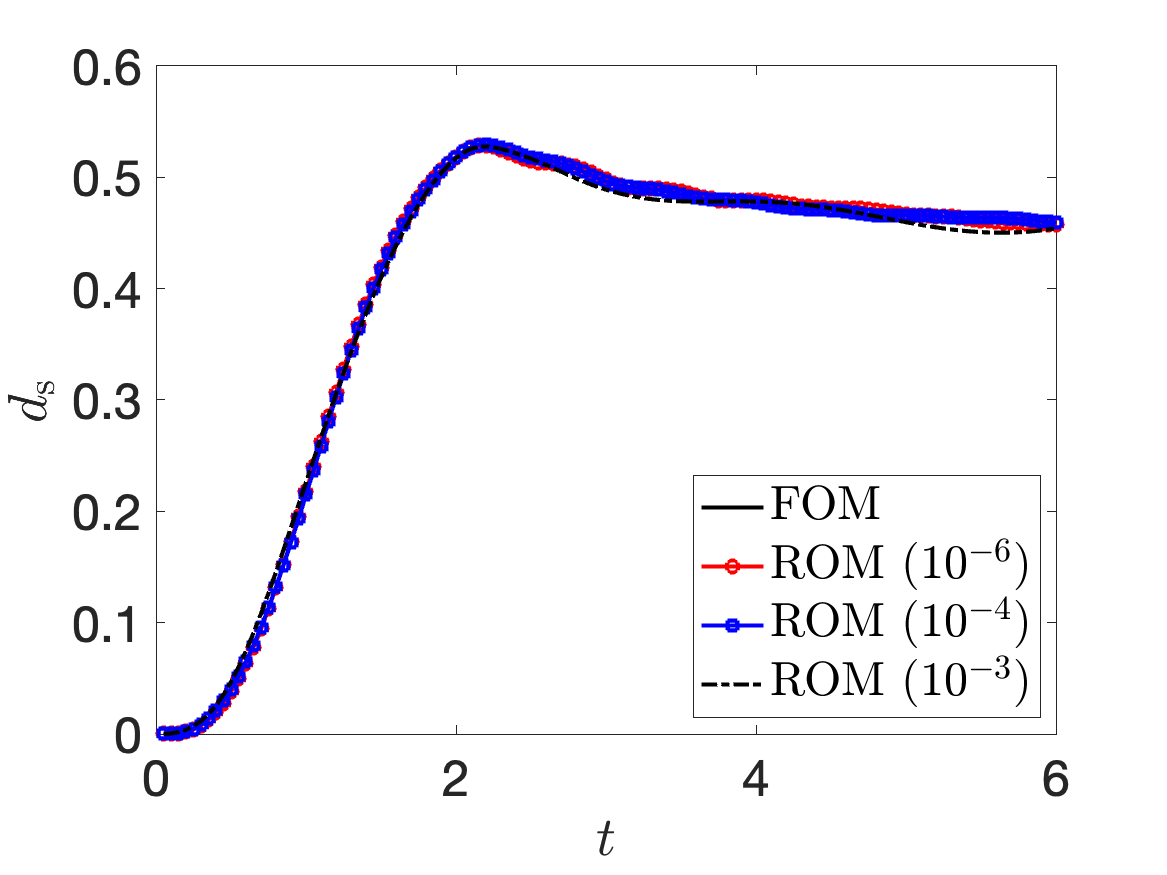}}
~~
\subfloat[]{\includegraphics[width=0.33\textwidth]{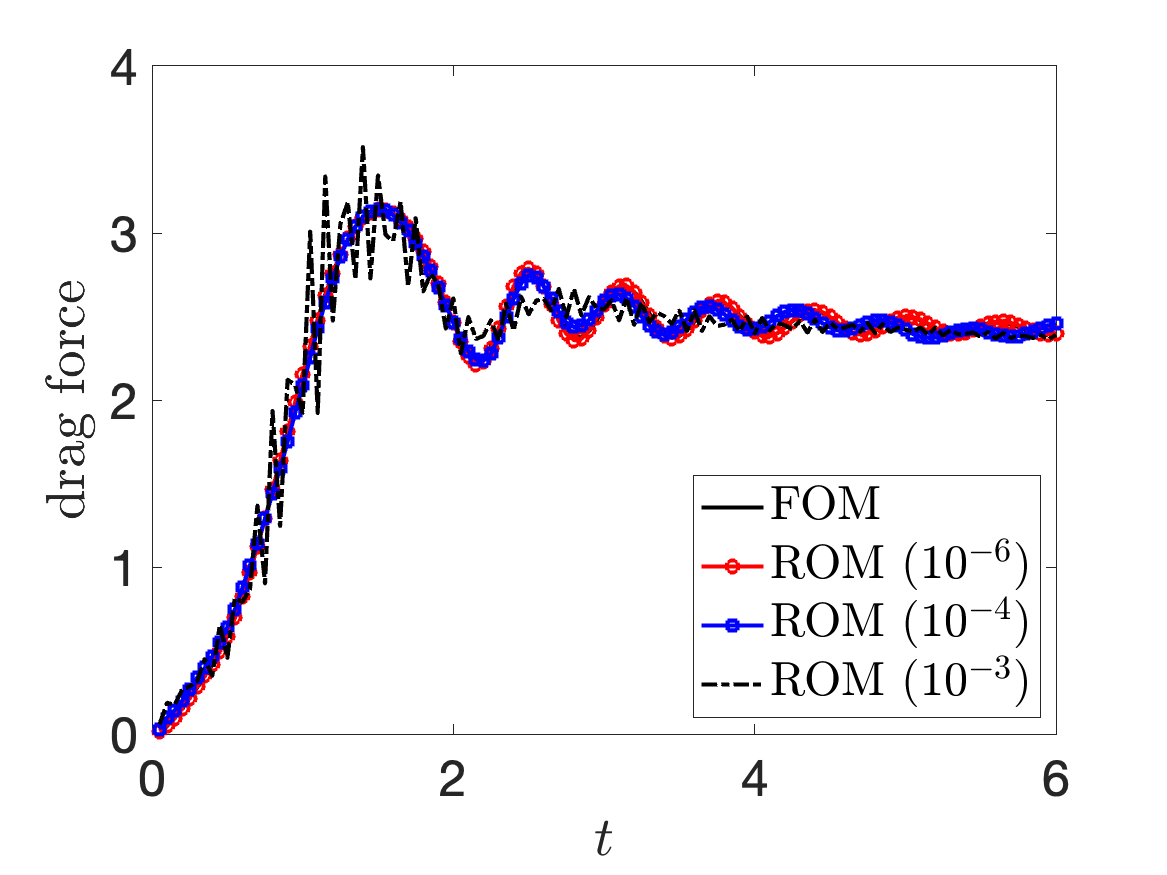}}
~~
\subfloat[]{\includegraphics[width=0.33\textwidth]{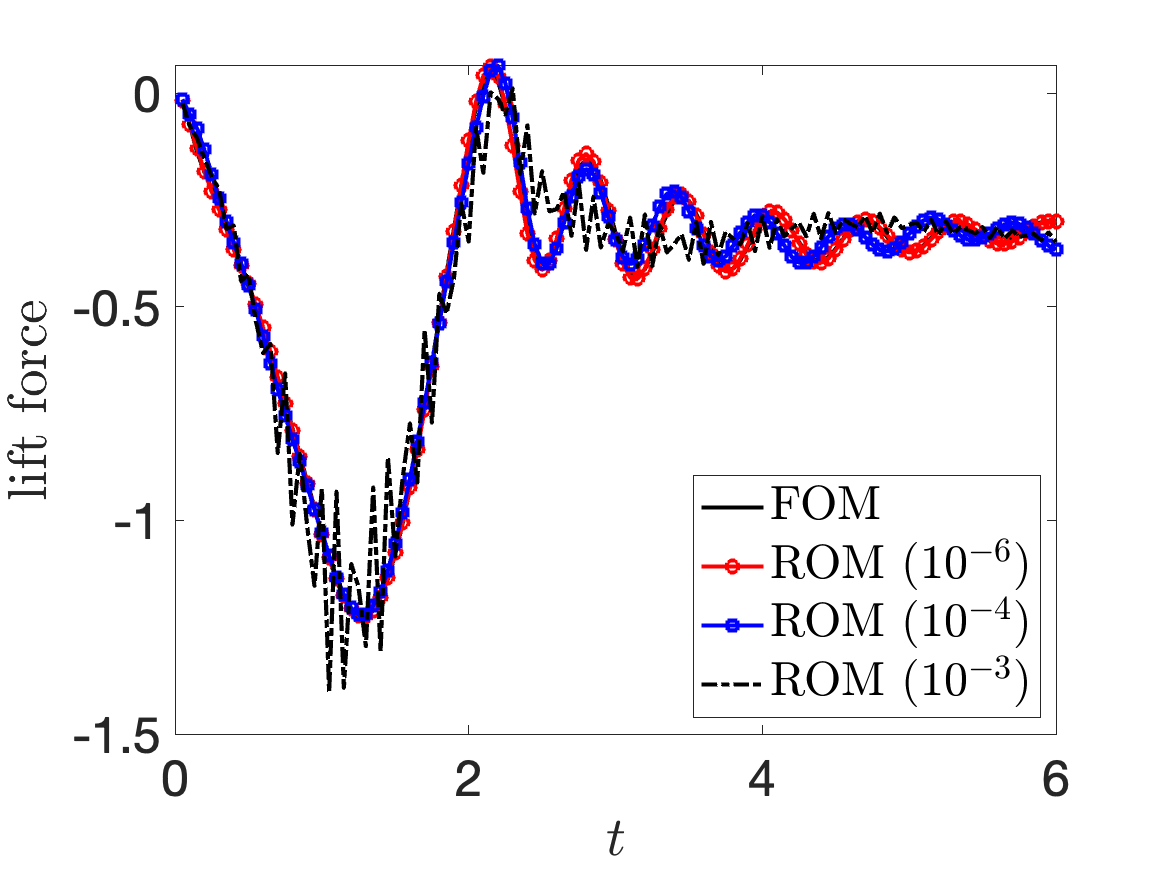}}
\caption{elastic beam; performance of the global ROM. 
(a) horizontal displacement of the top-right corner of the structure.
(b)-(c) drag and lift force.}
\label{fig:vbeam_MOR2}
\end{figure}

Figures \ref{fig:vbeam_MOR1} and \ref{fig:vbeam_MOR2} show the performance of the global ROM; here, we consider the HF data associated with the coarser discretization and we store the snapshots at each time step.
Figure \ref{fig:vbeam_MOR1}(a) shows the 
energy content of the discarded POD modes
$1 - \frac{\sum_{k=1}^{n} \lambda_k}{\sum_{j=1}^{n_{\rm train}} \lambda_j}$. We notice that the decay of the eigenvalues associated with the fluid state is considerably slower.
Figures \ref{fig:vbeam_MOR1}(b) and (c) show the behavior of the average relative error for the fluid and the solid states:
\begin{align}
E_{\rm f}^{\rm rel} :=
\dfrac{
\sqrt{ 
  \sum_{k=1}^K  
\big\|
\widehat{\widetilde{w}}_{\rm f}^{ (k)}
-
\widetilde{w}_{\rm f}^{{\rm hf}, (k)}
\big\|_{\rm f}^2
}
}{
\sqrt{ 
  \sum_{k=1}^K  
\big\|
\widetilde{w}_{\rm f}^{{\rm hf}, (k)}
\big\|_{\rm f}^2
}
},
\quad
E_{\rm s}^{\rm rel} :=
\dfrac{
\sqrt{ 
  \sum_{k=1}^K  
\big\|
\widehat{d}_{\rm s}^{ (k)}
-
d_{\rm s}^{{\rm hf}, (k)}
\big\|_{\rm s}^2
}
}{
\sqrt{ 
  \sum_{k=1}^K  
\big\|
d_{\rm s}^{{\rm hf}, (k)}
\big\|_{\rm s}^2
}
}, 
\label{eq:rel_error}
\end{align}
and the corresponding number of modes, for four different choices of ${\rm tol}_{\rm pod}$ in \eqref{eq:energy_criterion} and for 
${\rm tol}_{\rm en}=0.1$ in \eqref{eq:projection_criterion}.
The choice ${\rm tol}_{\rm en}=0.1$ leads to a number of enriching modes of the same order as $n_{\rm c}$ for all cases considered.

Figure \ref{fig:vbeam_MOR2} (a) shows the 
horizontal displacement of the top-right corner of the structure, while
Figures \ref{fig:vbeam_MOR2} (b) and (c) show the temporal behavior of the   drag force and the lift force exerted by the fluid on the structure,
\begin{equation}
\label{eq:drag_lift}
F_{\rm drag}(t) = \int_{\Gamma(t)} \sigma_{\rm f}(t) n_{\rm f}(t) \cdot e_1 \, dx,
\quad
F_{\rm lift}(t) = \int_{\Gamma(t)} \sigma_{\rm f}(t) n_{\rm f}(t) \cdot e_2 \, dx,
\quad
{\rm with} \;\;
e_1 = \left[\begin{array}{l}
1\\
 0 \\ 
\end{array} \right],
\;\;
e_2 = \left[\begin{array}{l}
0\\
1 \\ 
\end{array} \right].
\end{equation}
We notice that all the three ROMs reproduce  the displacement of the structure and the drag force, while 
only the ROM  associated with ${\rm tol}_{\rm pod}=10^{-6}$ accurately reproduces the lift force.

\begin{figure}[H]
\centering
\subfloat[]{\includegraphics[width=0.33\textwidth]{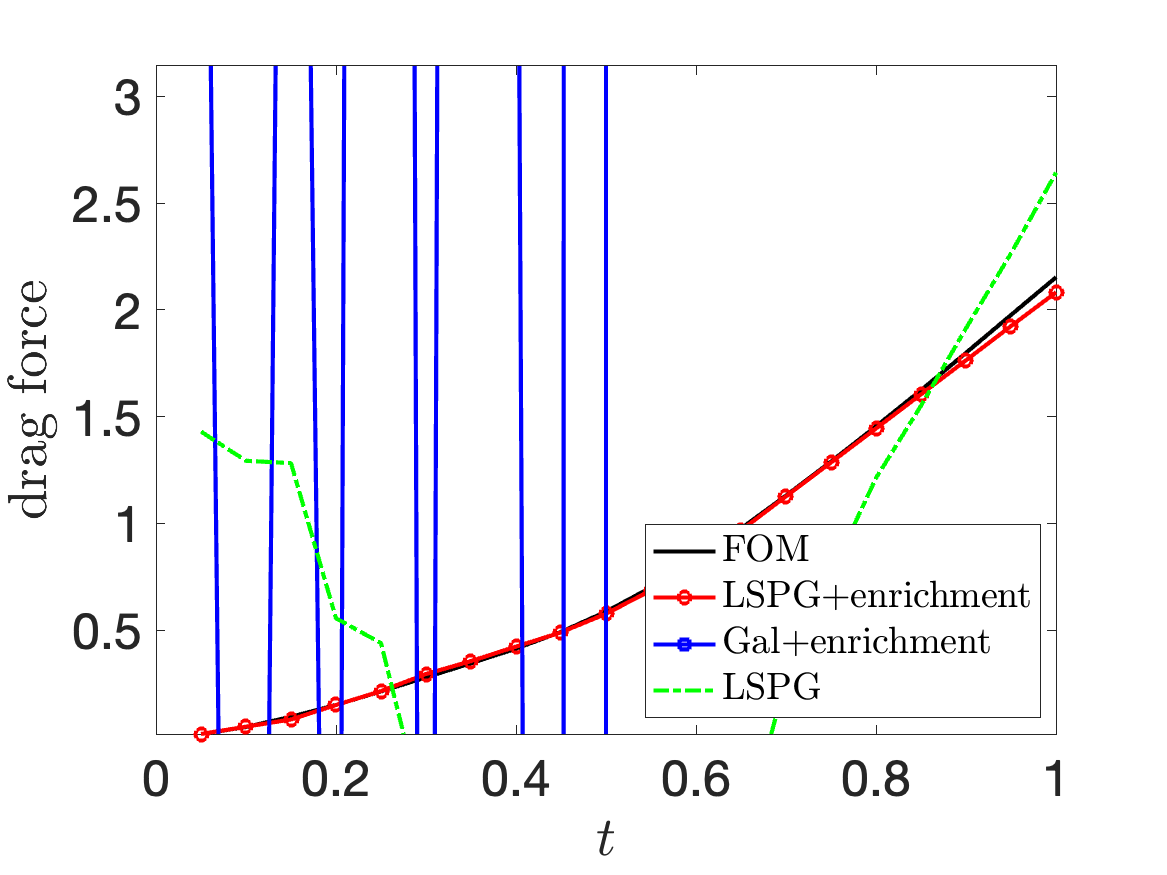}}
~~
\subfloat[]{\includegraphics[width=0.33\textwidth]{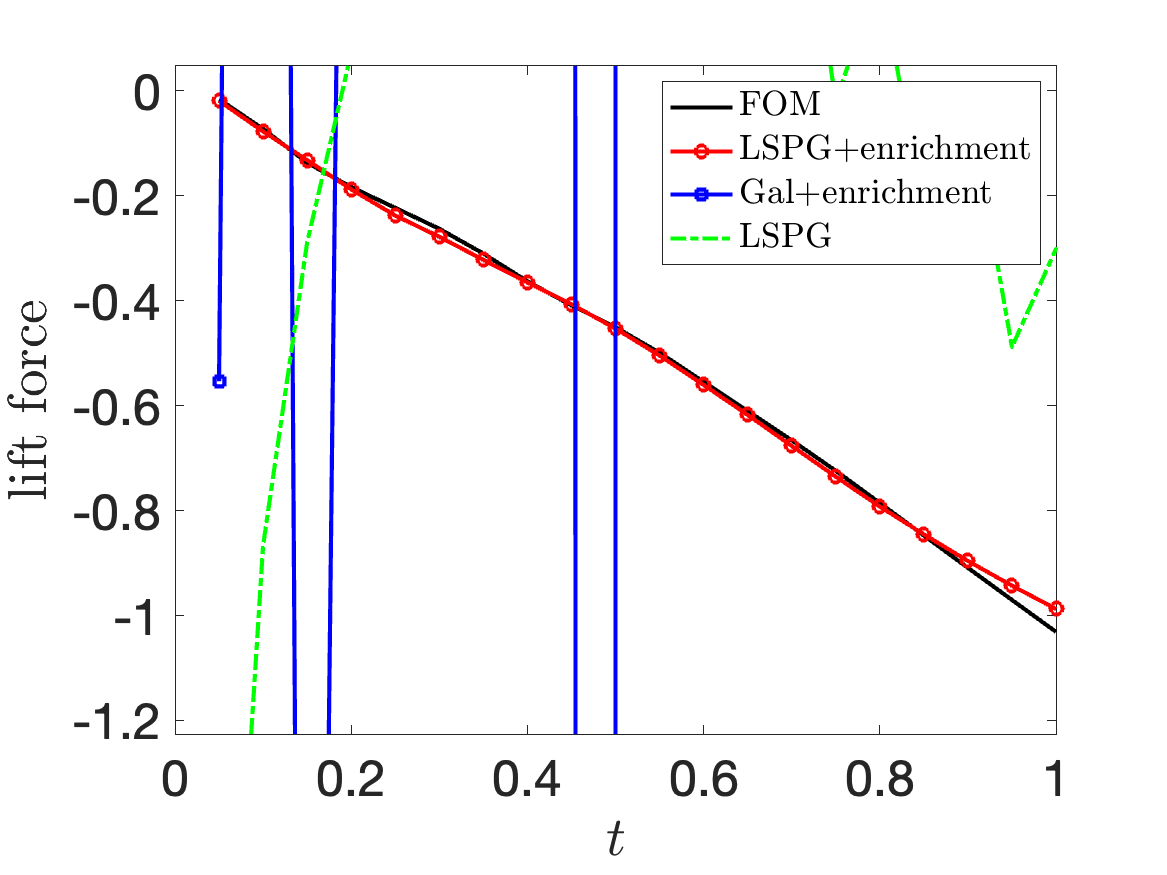}}
\caption{elastic beam; importance of enrichment and Petrov-Galerkin projection for the fluid subproblem. 
(a)-(b) drag and lift force (${\rm tol}_{\rm pod}=10^{-5}$).}
\label{fig:vbeam_MOR3}
\end{figure}

We conclude this section by investigating the impact of the enrichment 
strategy discussed in section \ref{sec:data_compression} and of least-square Petrov-Galerkin projection for the fluid subproblem (cf. section \ref{sec:local_roms}).
Figure \ref{fig:vbeam_MOR3} shows the behavior of the lift and the drag force for the global ROM with enrichment and LSPG projection of the fluid subproblem,
the global ROM based on Galerkin projection of the fluid subproblem, and the global ROM without enrichment; here,
 we consider the POD tolerance ${\rm tol}_{\rm pod}=10^{-5}$.
 We notice that both enrichment and Petrov-Galerkin projection are key to ensure convergence and also stability of the ROM.

\subsection{Turek problem}
\label{sec:turek}

\subsubsection{High-fidelity discretization}
\begin{figure}[H]
\centering

\begin{tikzpicture}[>={Stealth[length=2mm, width=1mm]}]  
    \draw[thick] (0,0) rectangle (10,2);  

    \draw[<->] (-0.3,0) -- (-0.3,2) node[midway, fill=white] {$H$};
    \draw[<->] (0,-0.3) -- (10,-0.3) node[midway, fill=white] {$L$};

    \draw[thick] (2,1) circle (0.5);  
    
    \filldraw [black] (2,1) circle (1pt) node[below right] {$C$};

    \pgfmathsetmacro{\rectHeight}{0.2} 
    \pgfmathsetmacro{\rectWidth}{3} 
    \pgfmathsetmacro{\rectStartY}{1 - \rectHeight/2} 

    \fill[gray] (2.5, \rectStartY) rectangle (2.5 + \rectWidth, \rectStartY + \rectHeight);
    \draw[thick] (2.5, \rectStartY) rectangle (2.5 + \rectWidth, \rectStartY + \rectHeight);

    \draw[<->]  
        (2.5 + \rectWidth + 0.3, \rectStartY) -- 
        (2.5 + \rectWidth + 0.3, \rectStartY + \rectHeight) node[midway,fill=white, rotate=90] {$h$};

\filldraw [black] (2.5 + \rectWidth, \rectStartY+\rectHeight/2) circle (1pt) node[above right] {$A$};

    \draw[<->]  
        (2.5, \rectStartY - 0.3) -- 
        (2.5 + \rectWidth, \rectStartY - 0.3) node[midway, fill=white] {$l$};

    \node at (0,0) [below left] {(0, 0)};
\end{tikzpicture}
\caption{Geometric setup (not scaled) of the Turek test case.
Here, we set 
$H=0.41$, $L=2.5$, $r=0.05$, $l=0.354$, $h=0.02$, $C=(0.2,0.2)$;
$A=(0.6,0.2)$ is the control point.}
\label{fig:turek_geo}
\end{figure}
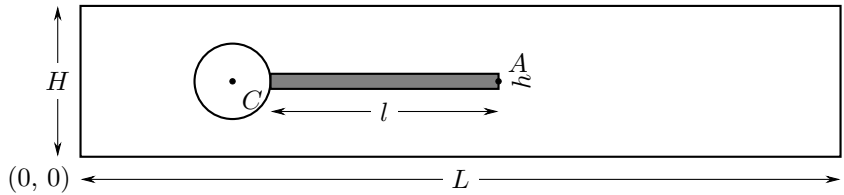

\begin{figure}[H]
  \centering
  
\subfloat[$t=30.2$]{\includegraphics[width=0.33\textwidth]{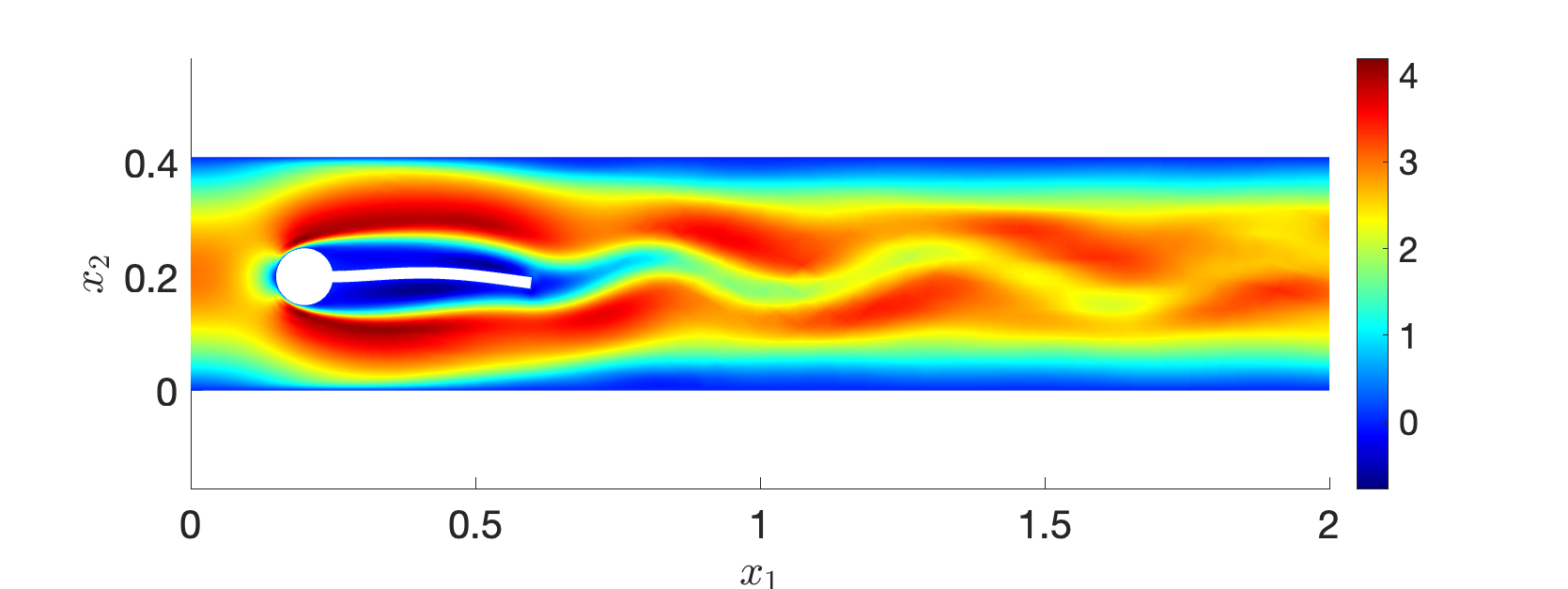}}
~~
\subfloat[$t=30.8$]{\includegraphics[width=0.33\textwidth]{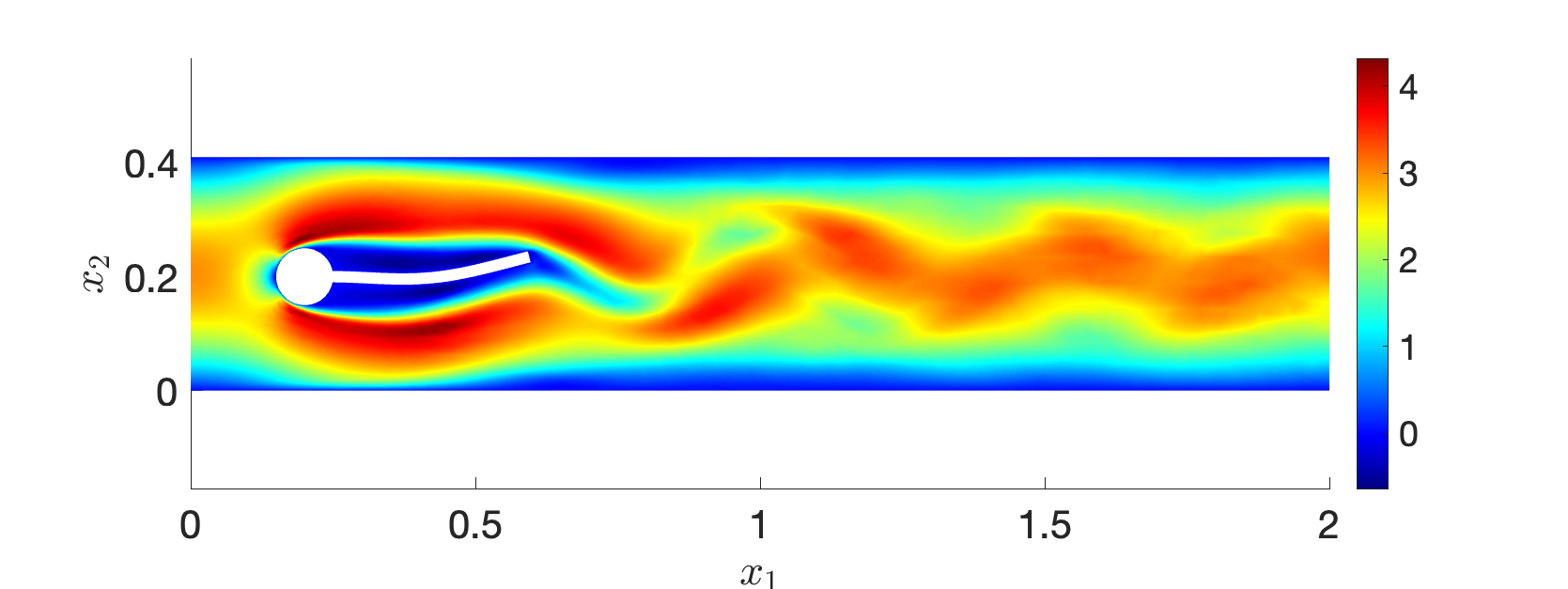}}
~~
\subfloat[$t=33$]{\includegraphics[width=0.33\textwidth]{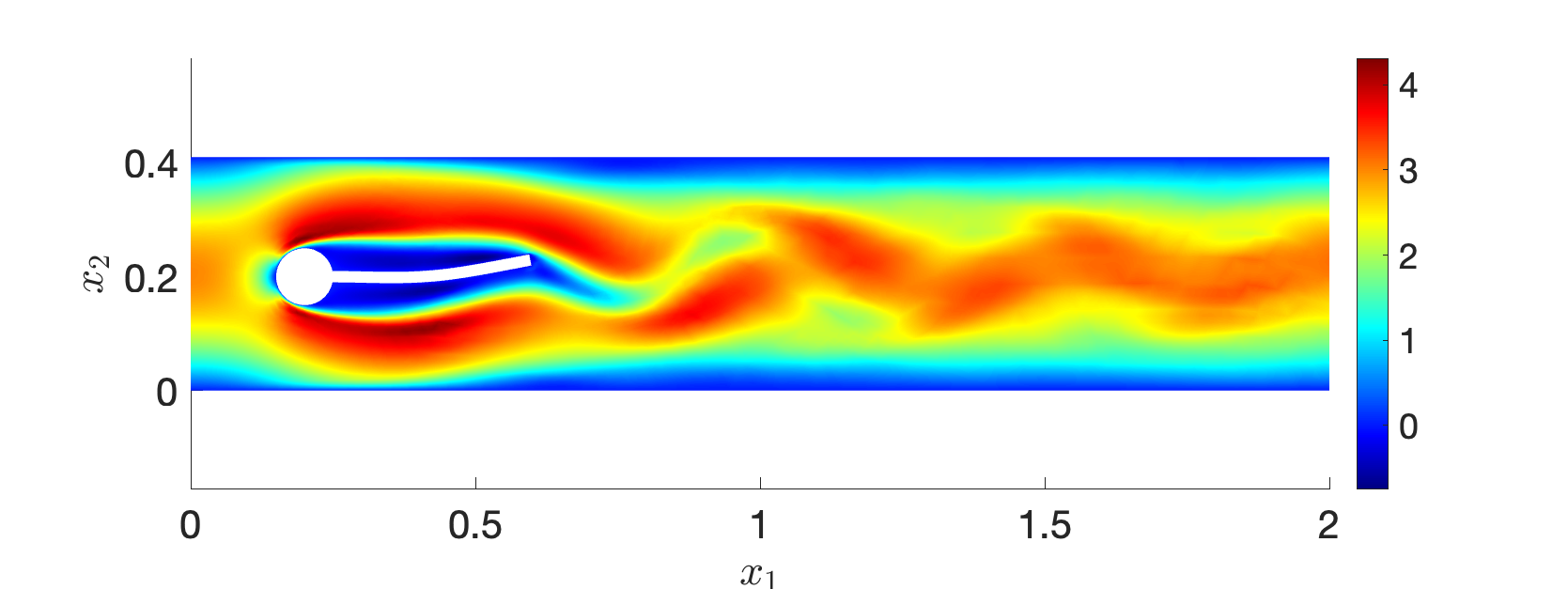}}
  \caption{Turek problem.
  Behavior of the streamwise (horizontal) velocity at  $3$ different time instants.}
  \label{fig:turek_hf_u}
\end{figure}

We study the viscous  flow past an elastic object in the laminar regime \cite{TurekHron2006}.  
Figure \ref{fig:turek_geo} illustrates the geometric configuration:  the height of the domain is
 $H =0.41$, and its length is  $L =2.5$; 
 the cylinder is centered at $C =(0.2, 0.2) $ and has a radius of $r =0.05$;
 the elastic beam is $l = 0.354$ in length  and  $h = 0.02$ in height, and is placed 
such that its right bottom corner is at $(0.6,\, 0.19)$. 
The top and bottom boundaries  of the fluid domain are treated as walls, the left boundary as an inlet with a specified flow velocity, and the right boundary as an outlet with a homogeneous Neumann condition. The rigid cylinder also acts as a fixed wall for the fluid subproblem.
The left  boundary of the beam is fully attached to the fixed cylinder, while the other three sides of the beam constitute the fluid-structure interface.
We consider the same physical parameters  as in 
the FSI3 test case of \cite{TurekHron2006}:
$$
\rho_{\rm f} = 10^3,
\quad
\rho_{\rm s} = 10^3,
\quad
E= 5.6 \cdot 10^6, 
\quad
\nu = 0.4,
\quad
\mu_{\rm f} = 1;
$$
 we impose a parabolic profile  at the inlet with a  maximum inflow velocity of $u_{\infty}=3$.
We initialize the FSI simulation by  holding the beam fixed for $T_0=30$ seconds to allow the vortex shedding to develop; then, we run the FSI problem for $T=10$ seconds. 
In our implementation, we divide the equations by $\rho_{\rm f}$ to improve the conditioning of the discrete problem.
Despite its simplicity, this setup leads to a complex fluid-structure interaction behavior that features significant deformations of the beam
(cf. Figure \ref{fig:turek_hf_u}).

\begin{figure}[H]
\centering

\subfloat[\texttt{mesh0}]{\includegraphics[width=0.45\textwidth]{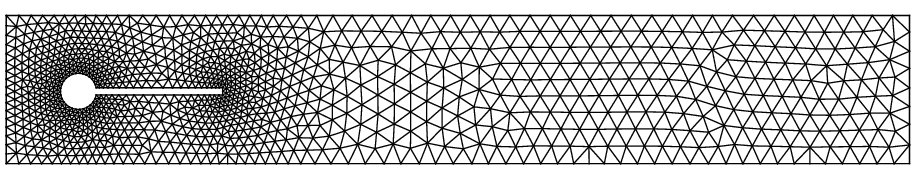}}
~~
\subfloat[\texttt{mesh1}]{\includegraphics[width=0.45\textwidth]{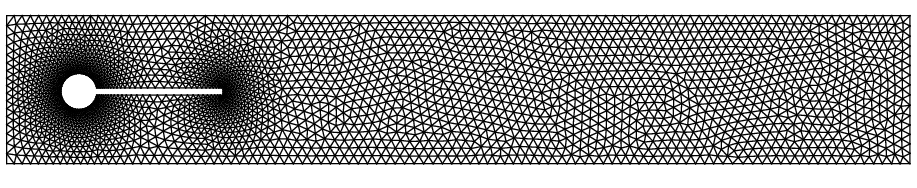}}
\\
\subfloat[\texttt{mesh2}]{\includegraphics[width=0.45\textwidth]{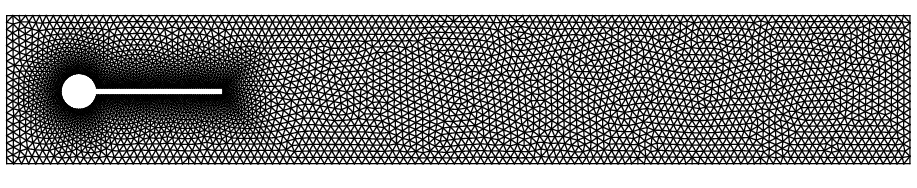}}
\caption{Turek; computational meshes.}
\label{fig:turek_meshes}
\end{figure}

We consider three different  discretizations with varying spatial degrees of freedom and time step sizes;
the  three spatial meshes  for the fluid subproblem are depicted in Figure \ref{fig:turek_meshes}. Table \ref{tab:turek_discretization} reports the size of the spatial meshes and of the temporal discretization; we also report the sampling time  associated with the collected snapshots.

\begin{table}[h!]
\centering

\begin{tabular}{lccccc}
\hline
\textbf{Discretization}        & $N_{\rm u}$ & $N_{\rm s}$ &$N_{\rm c}$ & $\Delta t$ & $\Delta t_s$ \\ \hline
\texttt{mesh0}  & $10379$ & $294$ & $162$ & $0.01$ & $0.01$ \\ \hline
\texttt{mesh1}  & $26891$ & $746$ & $258$ & $0.001$ & $0.005$\\ \hline
\texttt{mesh2}  & $39670$ & $2252$ & $578$ & $0.001$ & $0.005$ \\ \hline
\end{tabular}
\caption{Turek problem; HF discretization.}
\label{tab:turek_discretization}
\end{table}

\begin{figure}[H]
\centering

\subfloat[]{\includegraphics[width=0.45\textwidth]{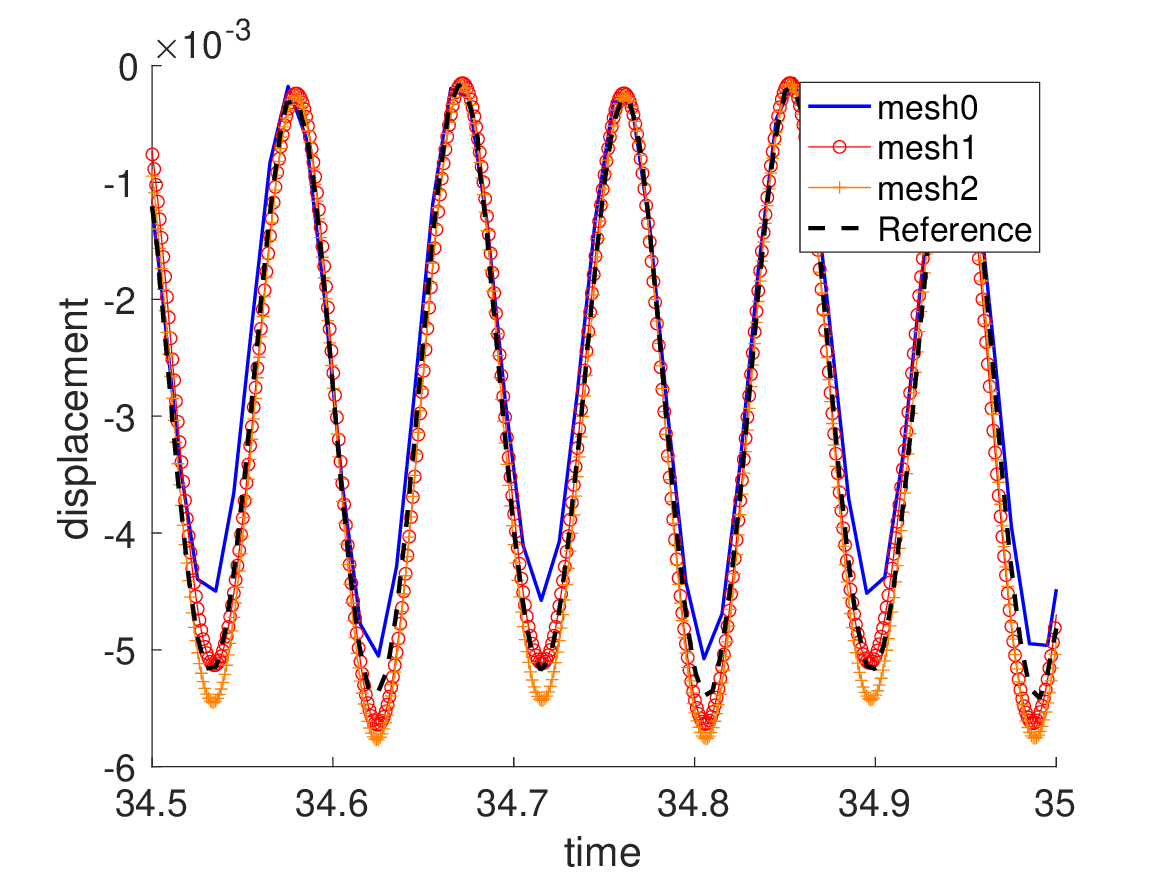}}
~~
\subfloat[]{\includegraphics[width=0.45\textwidth]{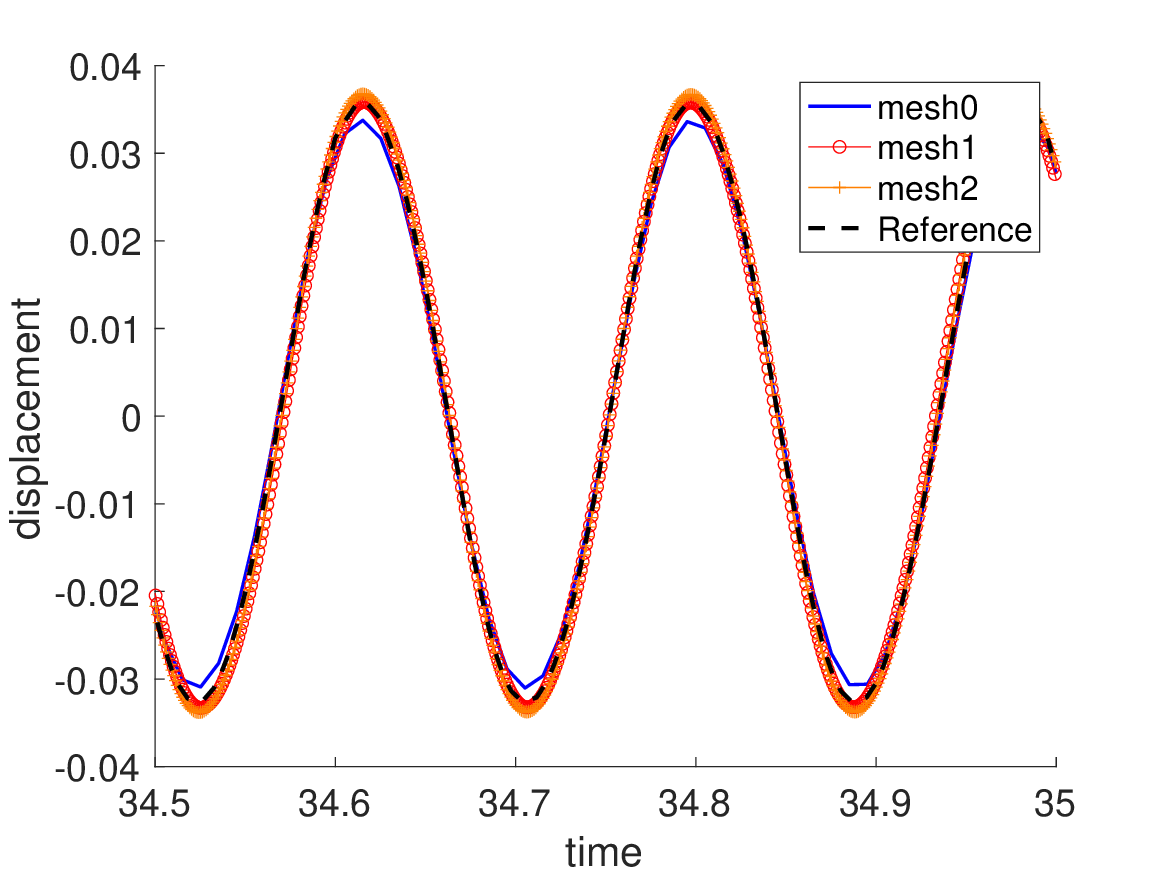}}
\\
\subfloat[]{\includegraphics[width=0.45\textwidth]{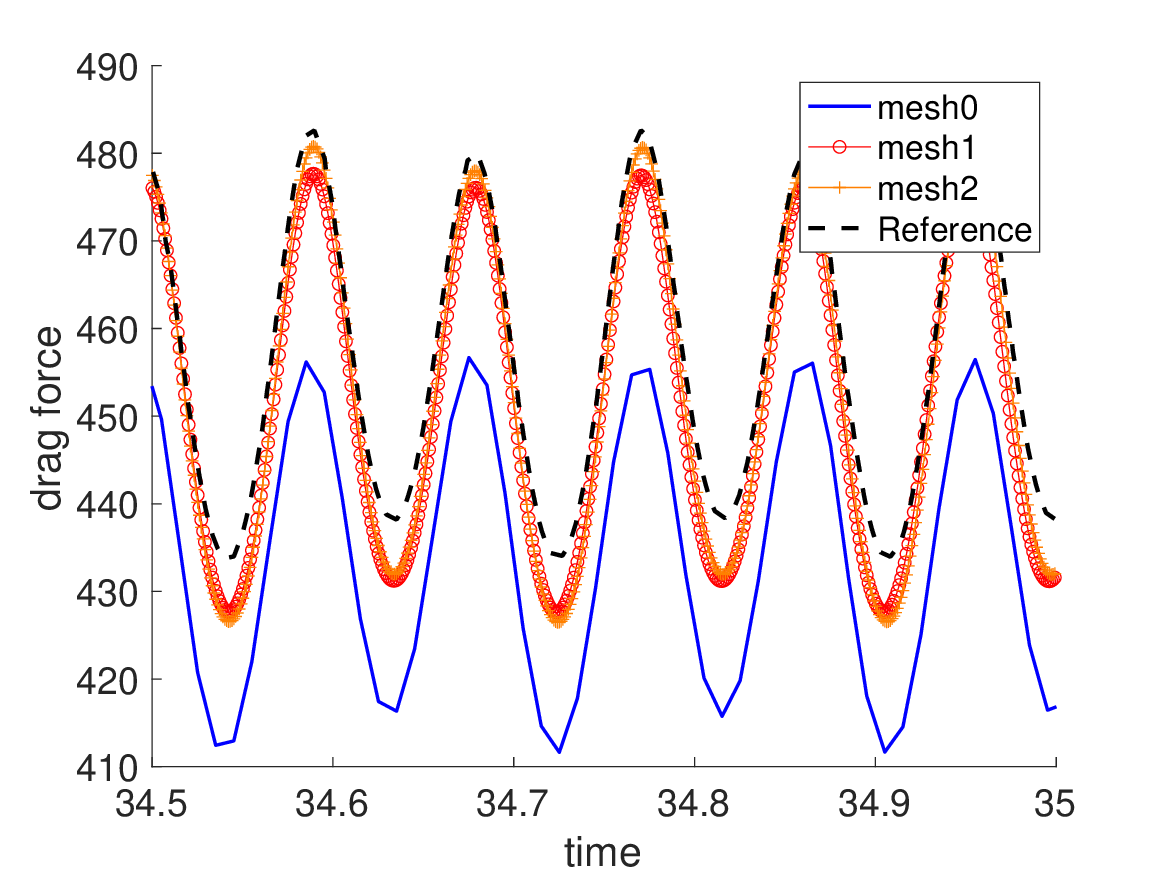}}
~~
\subfloat[]{\includegraphics[width=0.45\textwidth]{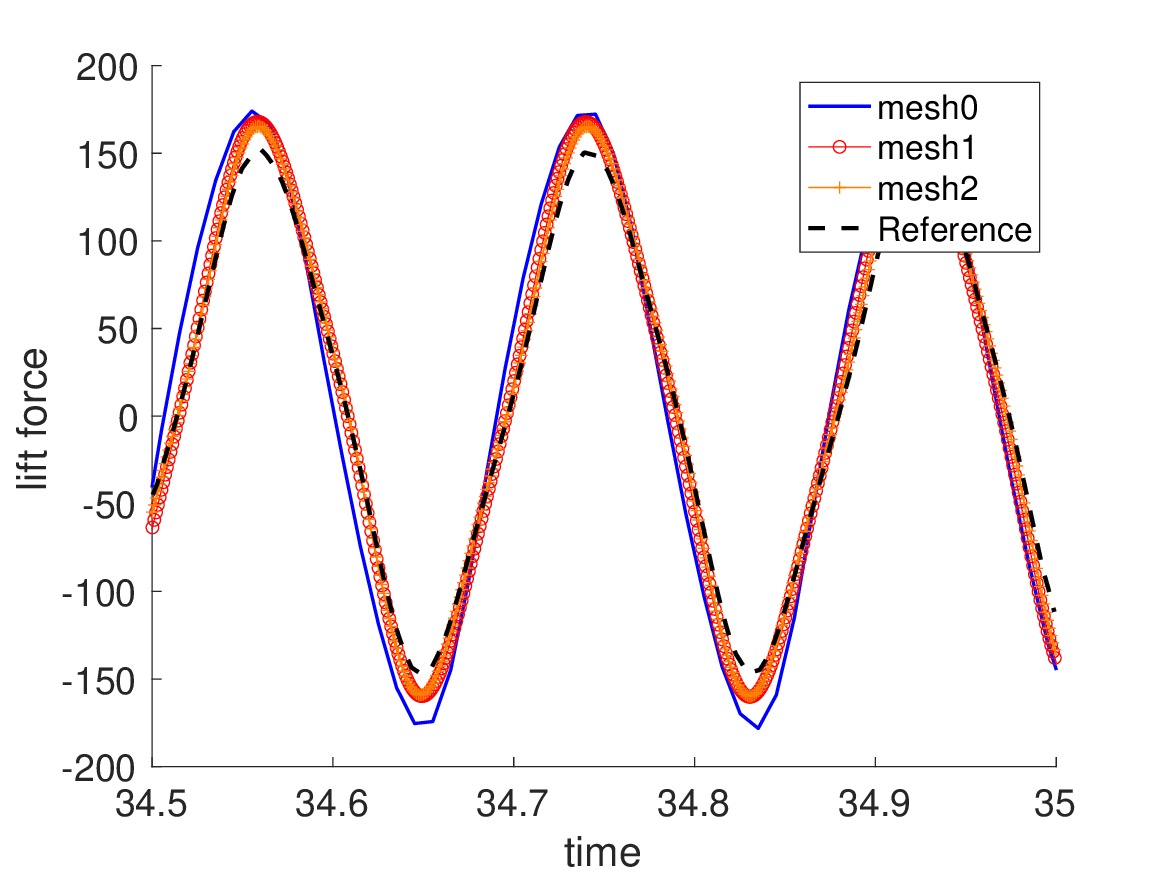}}
\caption{Turek problem; HF results.
(a)-(b) horizontal and vertical displacements of the control point $A$. 
(c)-(d) drag and lift forces  \eqref{eq:drag_lift}.}
\label{fig:turek_HF_results}
\end{figure}

Figure \ref{fig:turek_HF_results}
shows the numerical results of the HF partitioned model and compares them with the reference solution provided in \cite{TurekHron2006}.
Figures (a) and (b) show the horizontal and the vertical displacements of the control point $A$ depicted in Figure \ref{fig:turek_geo}.
Figures (c) and (d) show the drag and lift forces  as defined in  \eqref{eq:drag_lift}.
Additionally,
Table \ref{tab:turek_HF_results} reports 
the mean value, the  amplitude, and the frequency of  each measured quantity.

\begin{table}[h!]
\centering
\begin{tabular}{lcccc}
\hline
\textbf{Quantity}        & \textbf{\texttt{mesh0}} & \textbf{ \texttt{mesh1} } & \textbf{ \texttt{mesh2} } & \textbf{Reference} \\ \hline
$x$-displacement $[10^{-3}]$  & $-2.52 \pm 2.31$ & $-2.83 \pm 2.59$ & $-2.95 \pm 2.70$ & $-2.69 \pm 2.53$ \\
frequency (Hz)           & $11.00$            & $11.03$            & $11.03$            & $10.90$          \\ \hline
$y$-displacement $[10^{-3}]$  & $1.17 \pm 32.29$  & $1.25 \pm 34.52$  & $1.37 \pm 35.26$  & $1.48 \pm 34.38$ \\
frequency (Hz)           & $5.50$             & $5.52$             & $5.52$             & $5.30$           \\ \hline
drag force               & $434.6 \pm 21.13$  & $451.6 \pm 23.72$  & $452.7 \pm 25.08$  & $457.3 \pm 22.66$ \\
frequency (Hz)           & $11.00$            & $11.03$            & $11.03$            & $10.90$          \\ \hline
lift force               & $2.88 \pm 175.41$  & $5.39 \pm 163.80$  & $5.84 \pm 162.41$  & $2.22 \pm 149.78$ \\
frequency (Hz)           & $5.50$             & $5.52$             & $5.52$             & $5.30$           \\ \hline
\end{tabular}
\caption{Turek; mean, amplitude, and frequency of the measured quantities for \texttt{mesh0}, \texttt{mesh1}, and \texttt{mesh2}, compared to the reference solution.}
\label{tab:turek_HF_results}
\end{table}

The numerical results for all three meshes demonstrate good agreement with the reference solution. However, the mean drag force for \texttt{mesh0} exhibits a relatively larger deviation compared to the reference value, while the results for \texttt{mesh1} and \texttt{mesh2} align closely with the benchmark.

\subsubsection{Performance of the ROM-ROM solver}

\begin{figure}[H]
\centering

\subfloat[]{\includegraphics[width=0.33\textwidth]{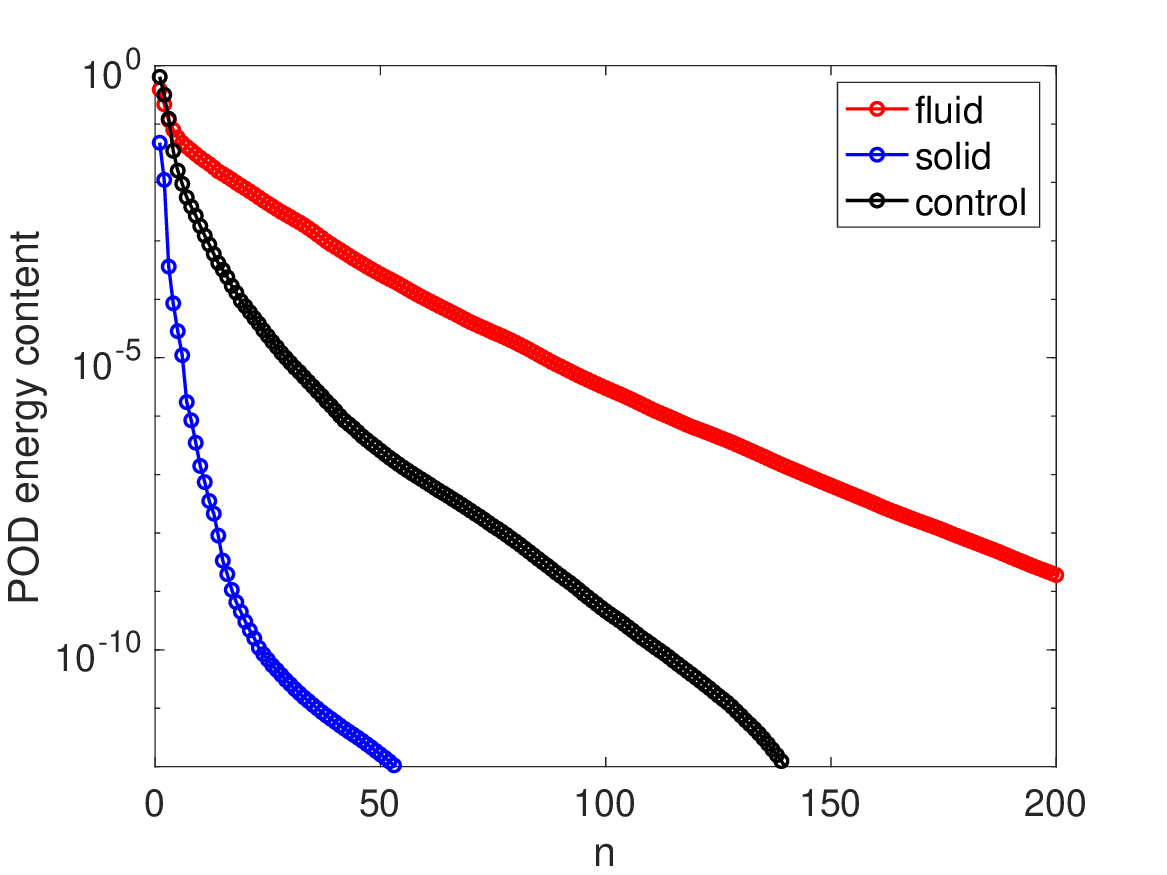}}
~~
\subfloat[]{\includegraphics[width=0.33\textwidth]{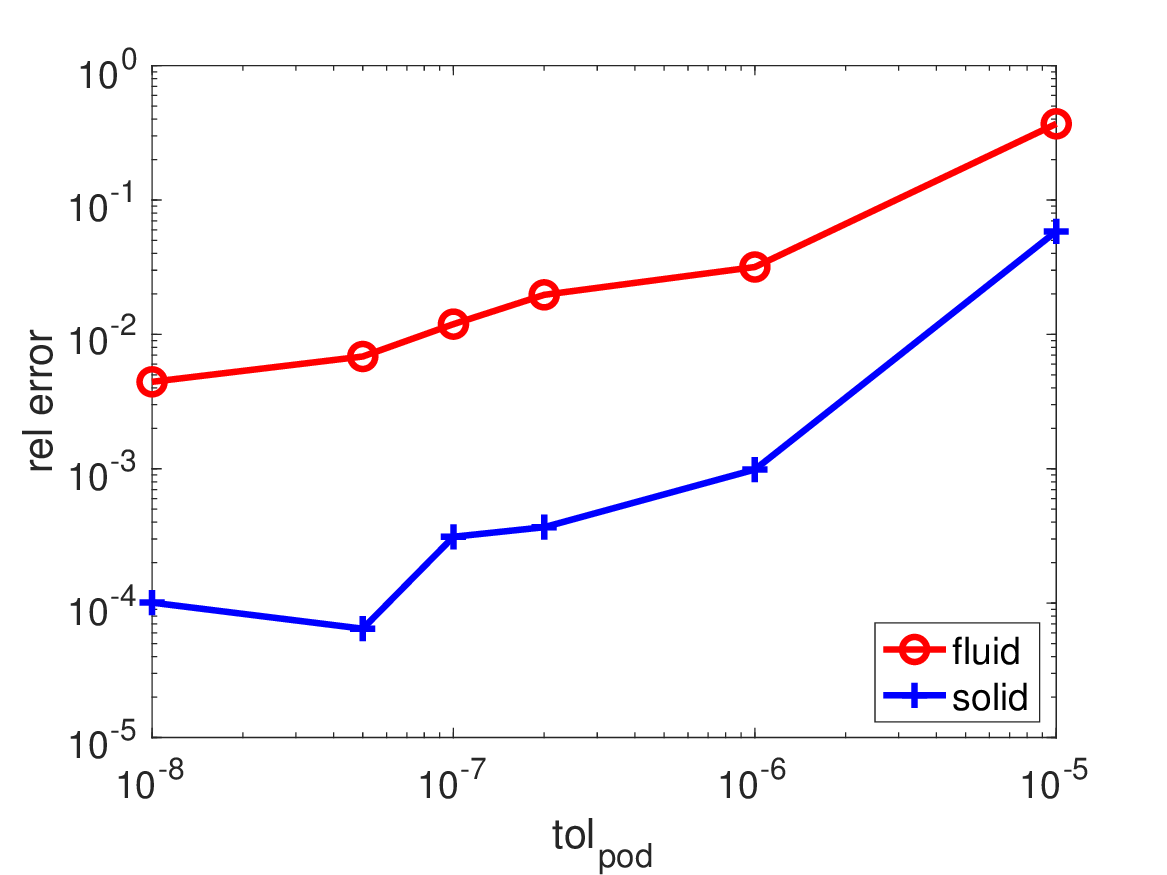}}
~~
\subfloat[]{\includegraphics[width=0.33\textwidth]{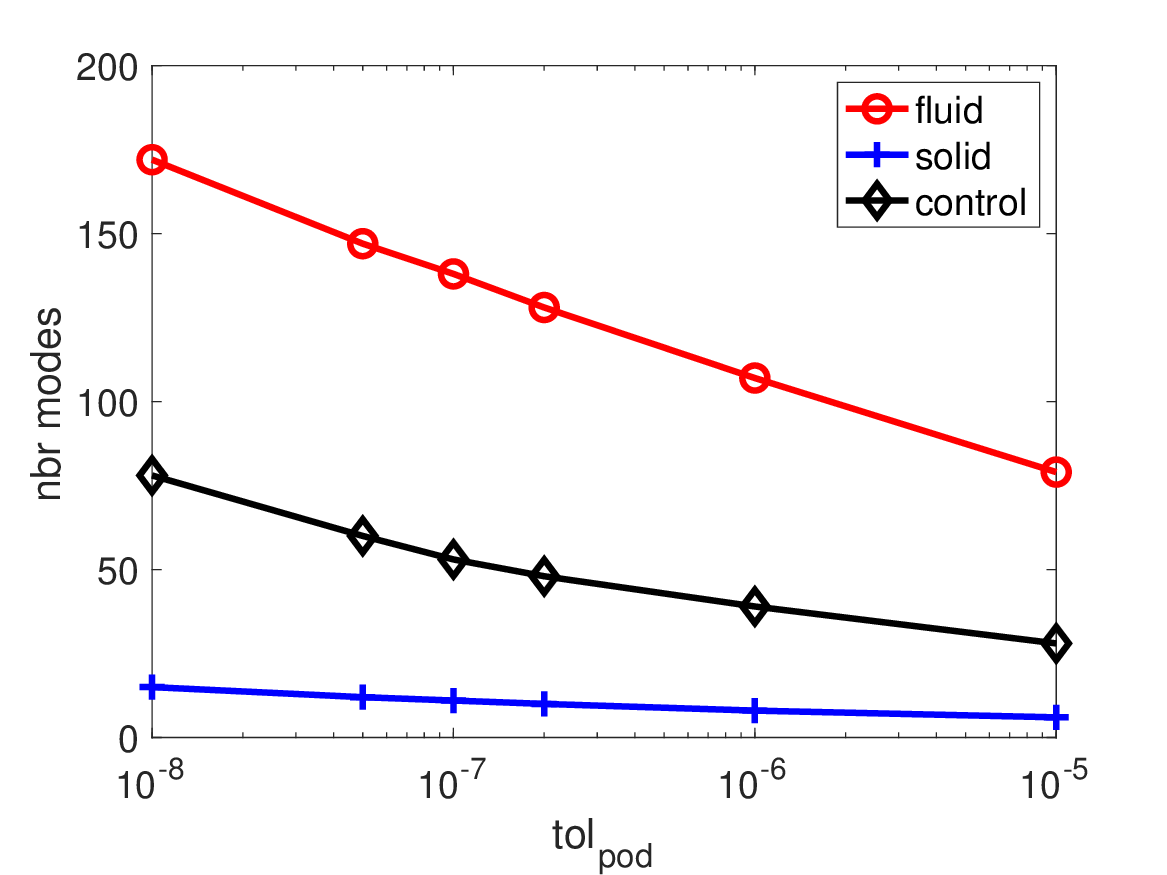}}
\caption{Turek problem; performance of the global ROM. 
(a) energy content of the discarded POD modes.
(b)-(c) relative error and number of modes for different choices of the tolerance ${\rm tol}_{\rm pod}$ with ${\rm tol}_{\rm en}=0.1$.}
\label{fig:turek_MOR1}
\end{figure}

We reproduce the HF solution associated with \texttt{mesh0} for a time duration of $3\,\mathrm{s}$. 
 Figure \ref{fig:turek_MOR1}(a) displays the behavior of the POD eigenvalues for the fluid state, the solid state, and the control, respectively. The eigenvalues associated with  the solid state decrease rapidly, while those of  the fluid state   diminish  much more slowly. 
We conclude that the solid subproblem and  (to a lesser  degree) the control are much more amenable for (linear)  model reduction: this finding motivates the development of the hybrid model of section \ref{sec:res_hybrid_solver}.
Figures \ref{fig:turek_MOR1}(b) and (c) show the behavior of the average relative error (cf. equation \eqref{eq:rel_error} ) for the fluid and the solid states
and the corresponding number of modes, for  different choices of ${\rm tol}_{\rm pod}$. 

\begin{figure}[H]
  \centering
  
\subfloat[$t=30.2$]{\includegraphics[width=0.33\textwidth]{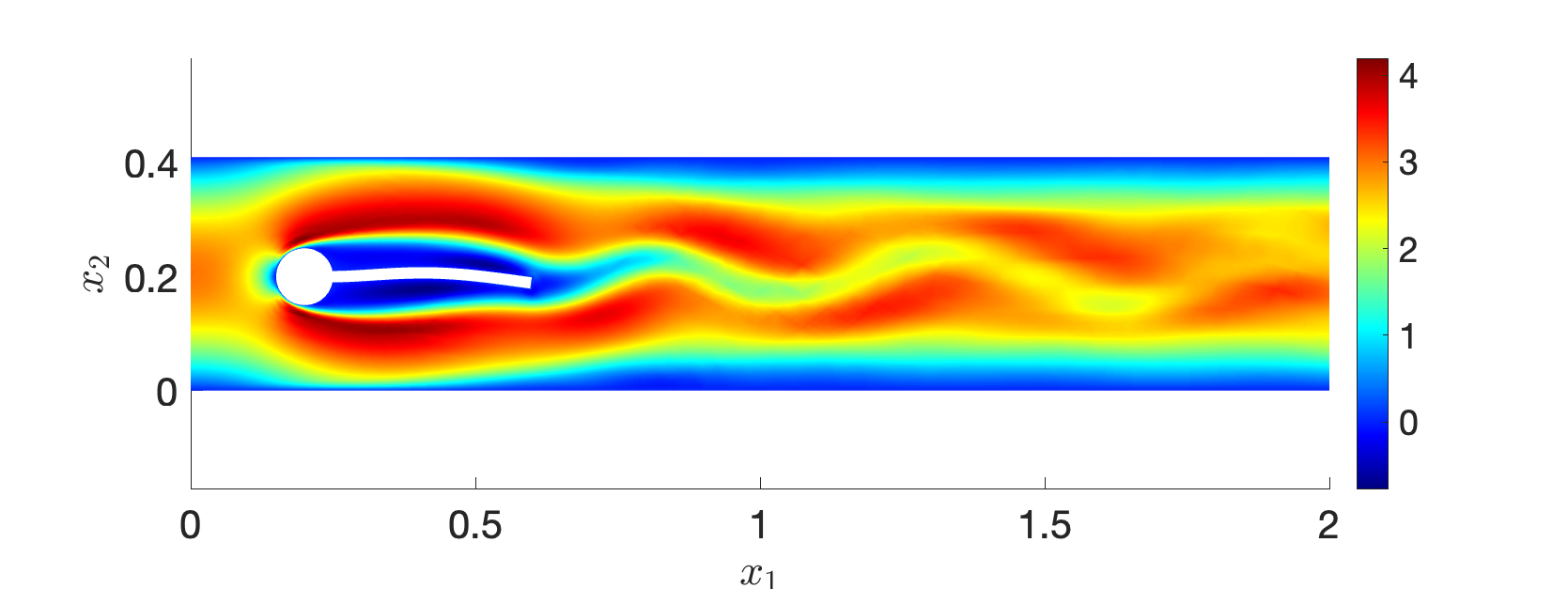}}
~~
\subfloat[$t=30.8$]{\includegraphics[width=0.33\textwidth]{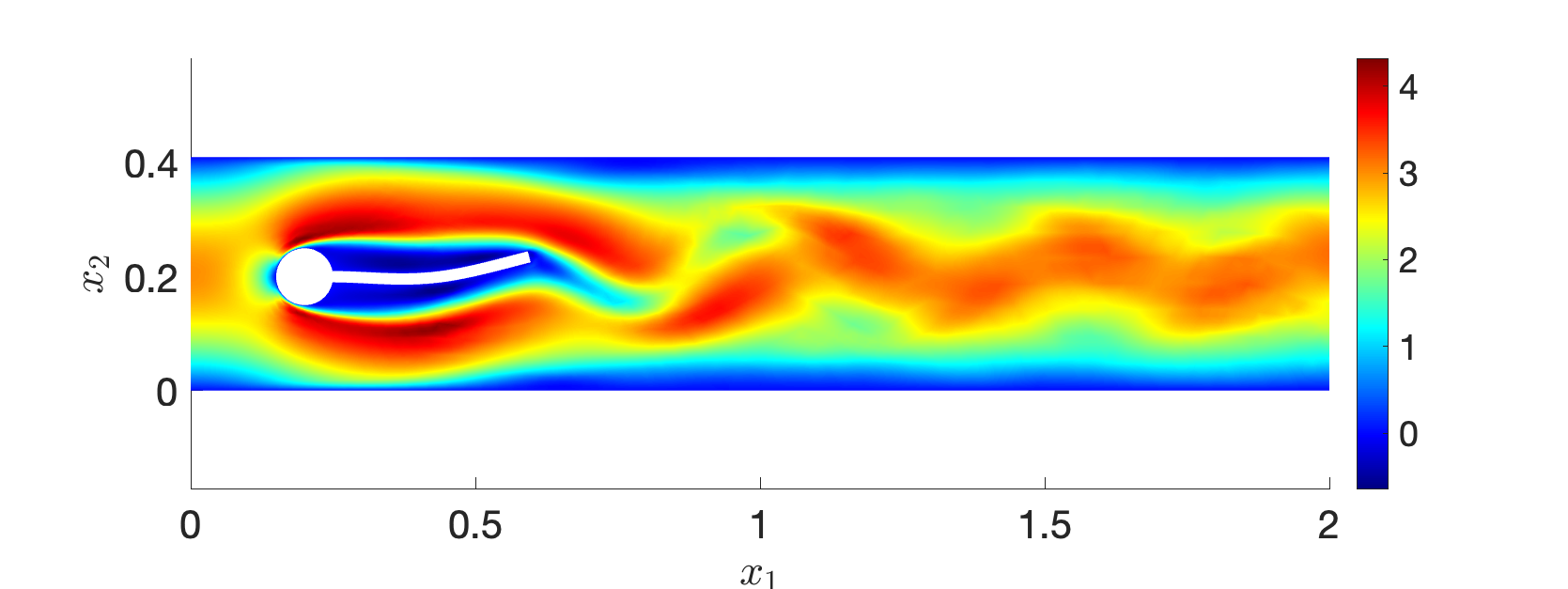}}
~~
\subfloat[$t=33$]{\includegraphics[width=0.33\textwidth]{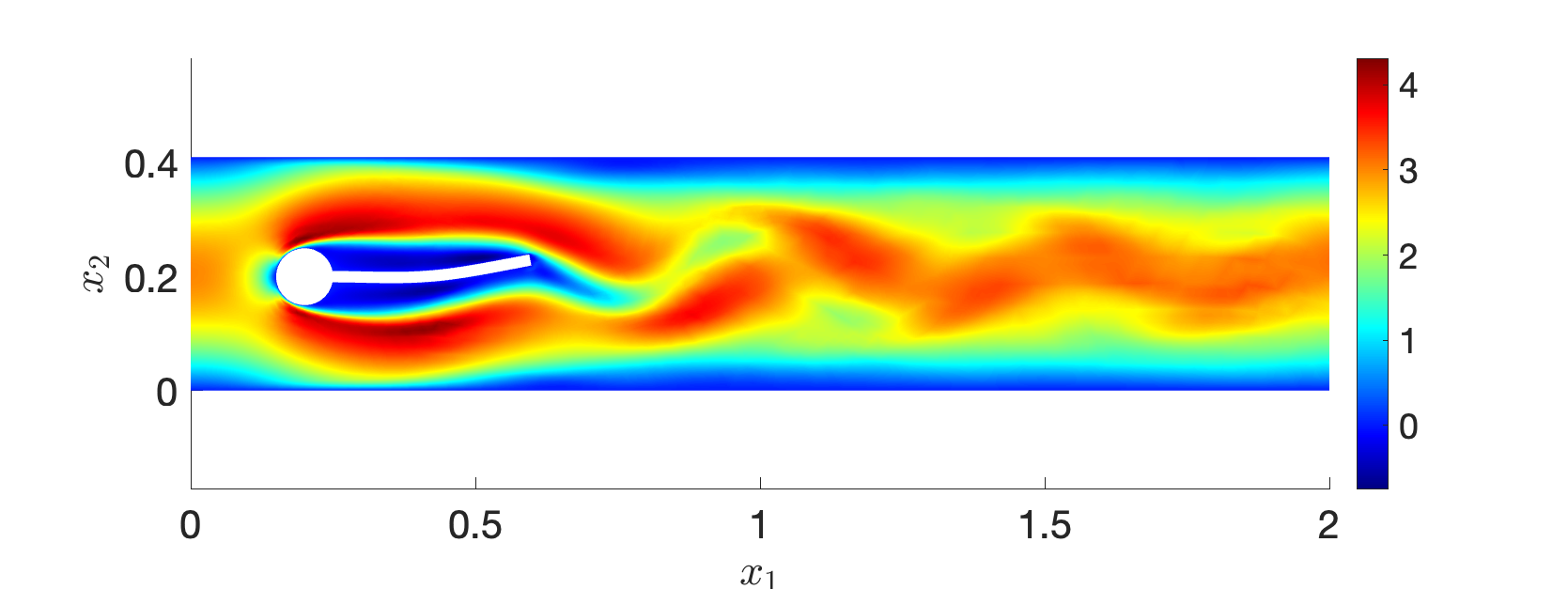}}

\subfloat[$t=30.2$]{\includegraphics[width=0.33\textwidth]{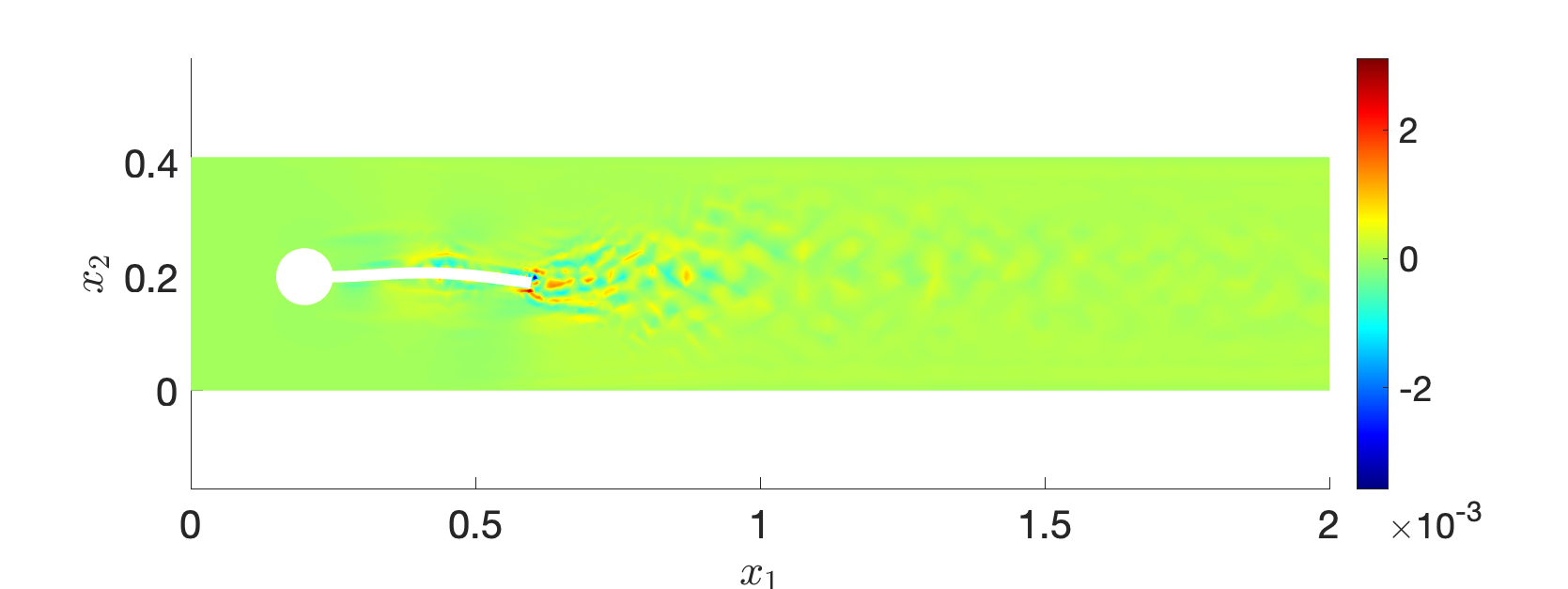}}
~~
\subfloat[$t=30.8$]{\includegraphics[width=0.33\textwidth]{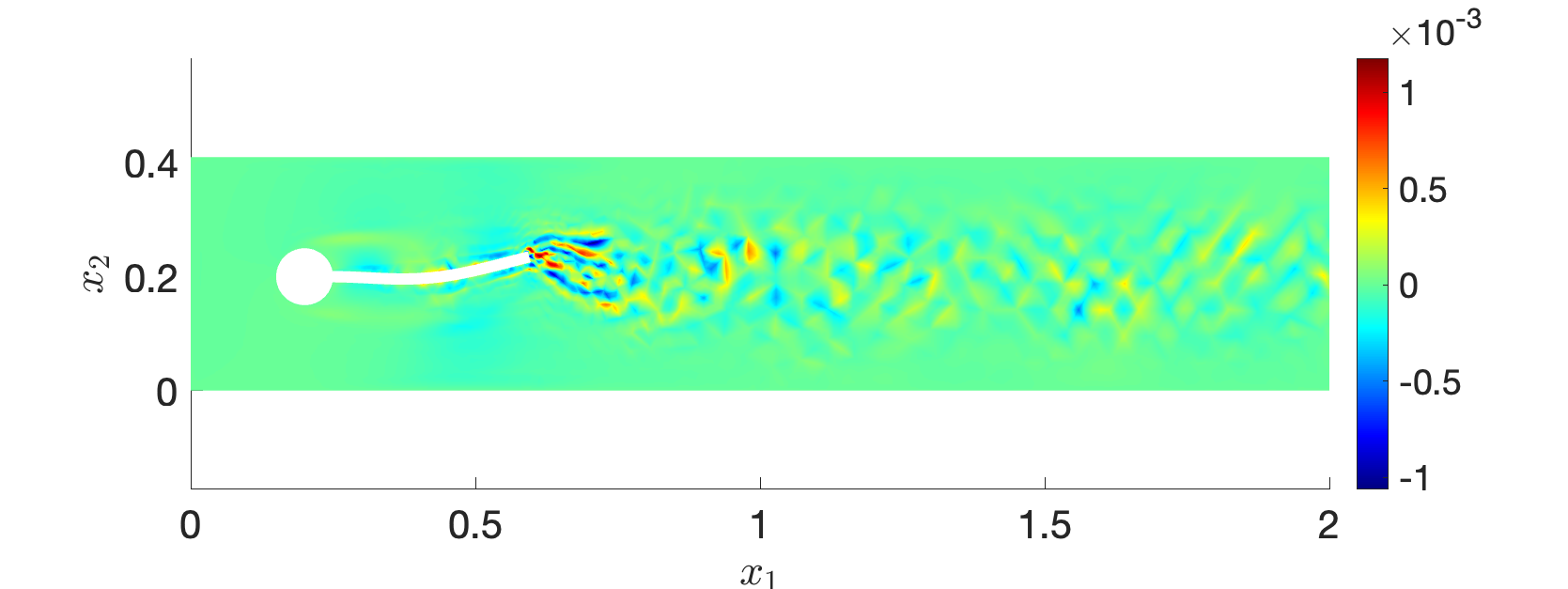}}
~~
\subfloat[$t=33$]{\includegraphics[width=0.33\textwidth]{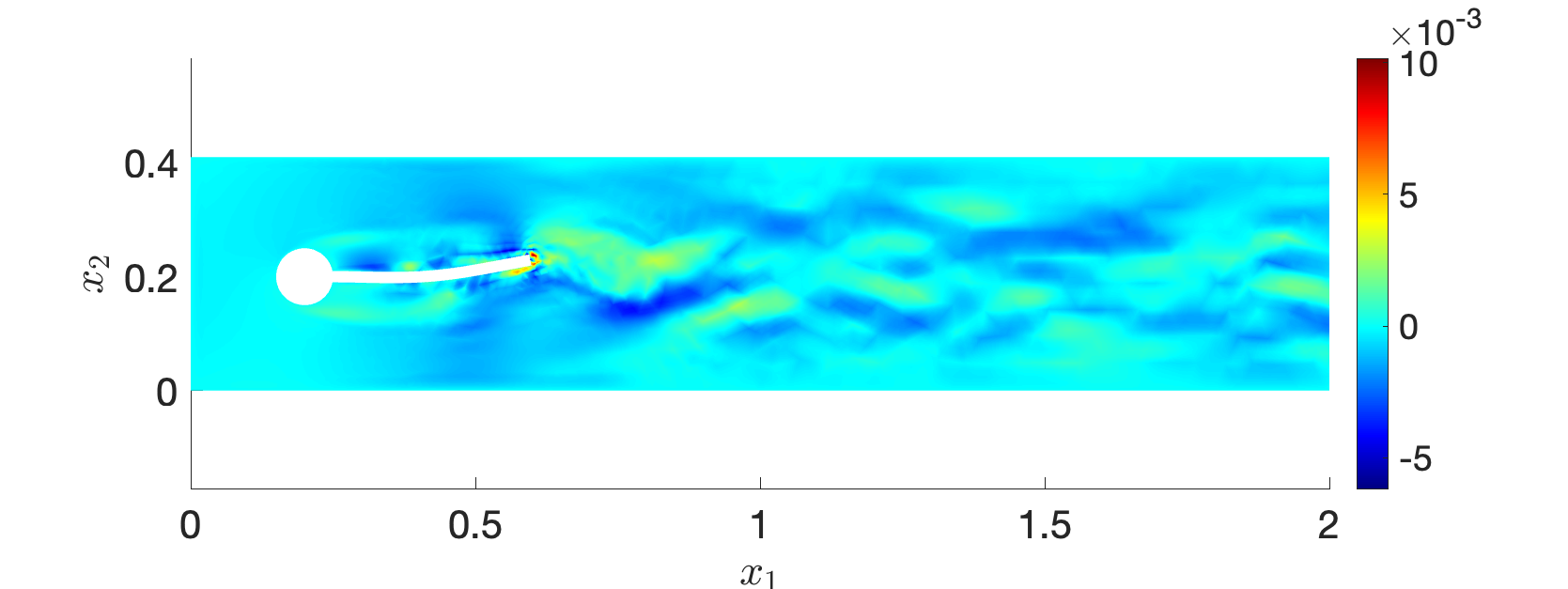}}
  \caption{Turek: ROM results for streamwise velocity at $3$ different time instants, (a)-(b)-(c) ROM results, (d)-(e)-(f) pointwise error of the ROM results relative to the HF results.}
  \label{fig:turek_rom_u}
\end{figure}

In Figure \ref{fig:turek_rom_u}, we 
show the estimates of the streamwise velocity at the same time instants shown in Figure \ref{fig:turek_hf_u} for the POD tolerance 
$10^{-6}$. We also report the behavior of the error field. Results of the ROM are in good agreement with the ones of the HF model.

\begin{figure}[H]
\centering
\includegraphics[width=8cm]{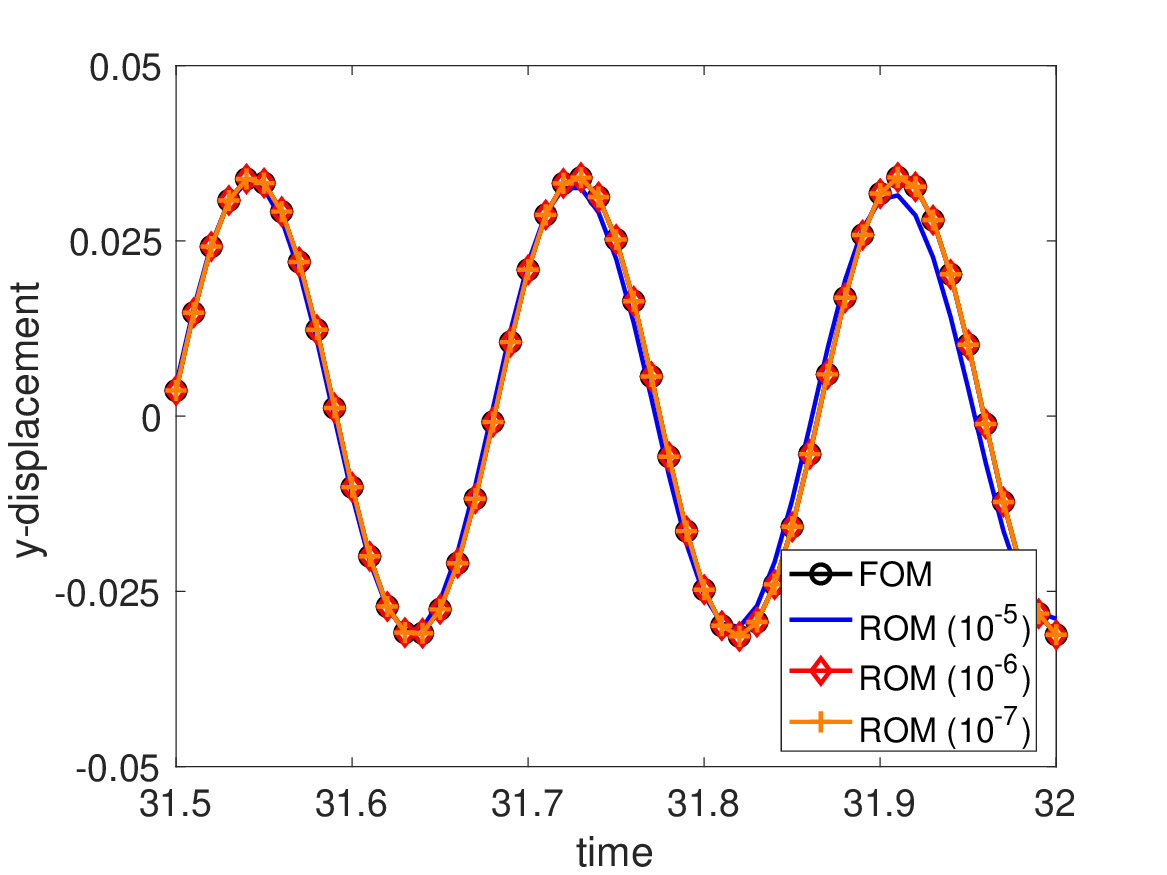}
\caption{Turek problem; vertical displacement of the control point $A$ for the partitioned ROM with different POD tolerances.}
\label{fig:turek_MOR2}
\end{figure}

In Figure \ref{fig:turek_MOR2}, we display the vertical displacement of the control point for different POD tolerances. 
The results of the ROM are in good agreement with the  HF results for all tolerances considered. We notice, however, that 
the use of larger tolerances (e.g., $\text{tol}_{\text{pod}}=10^{-4}$) leads to an unstable ROM.

\begin{figure}[H]
\centering
\subfloat[]{\includegraphics[width=0.45\textwidth]{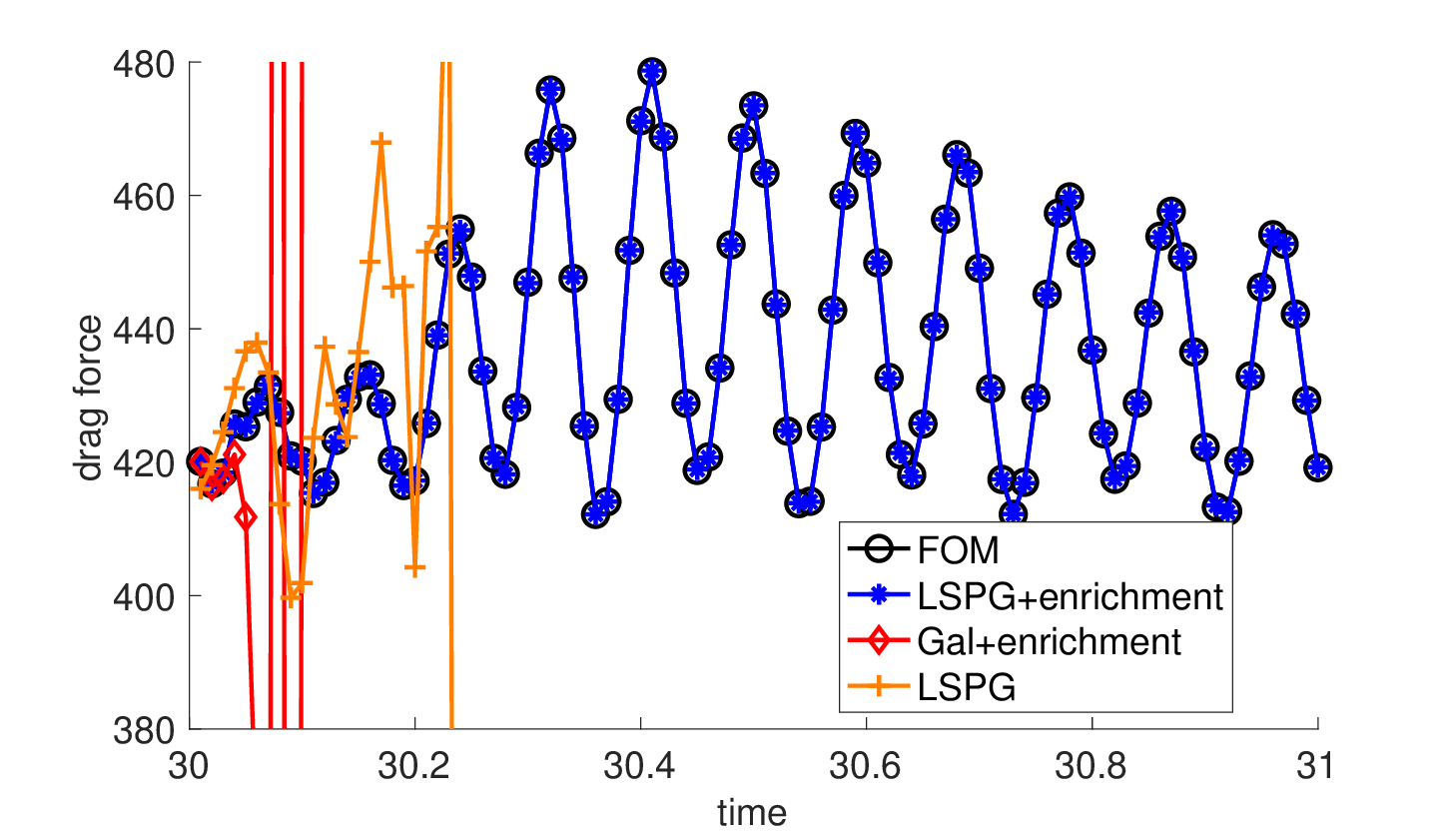}}
~~
\subfloat[]{\includegraphics[width=0.45\textwidth]{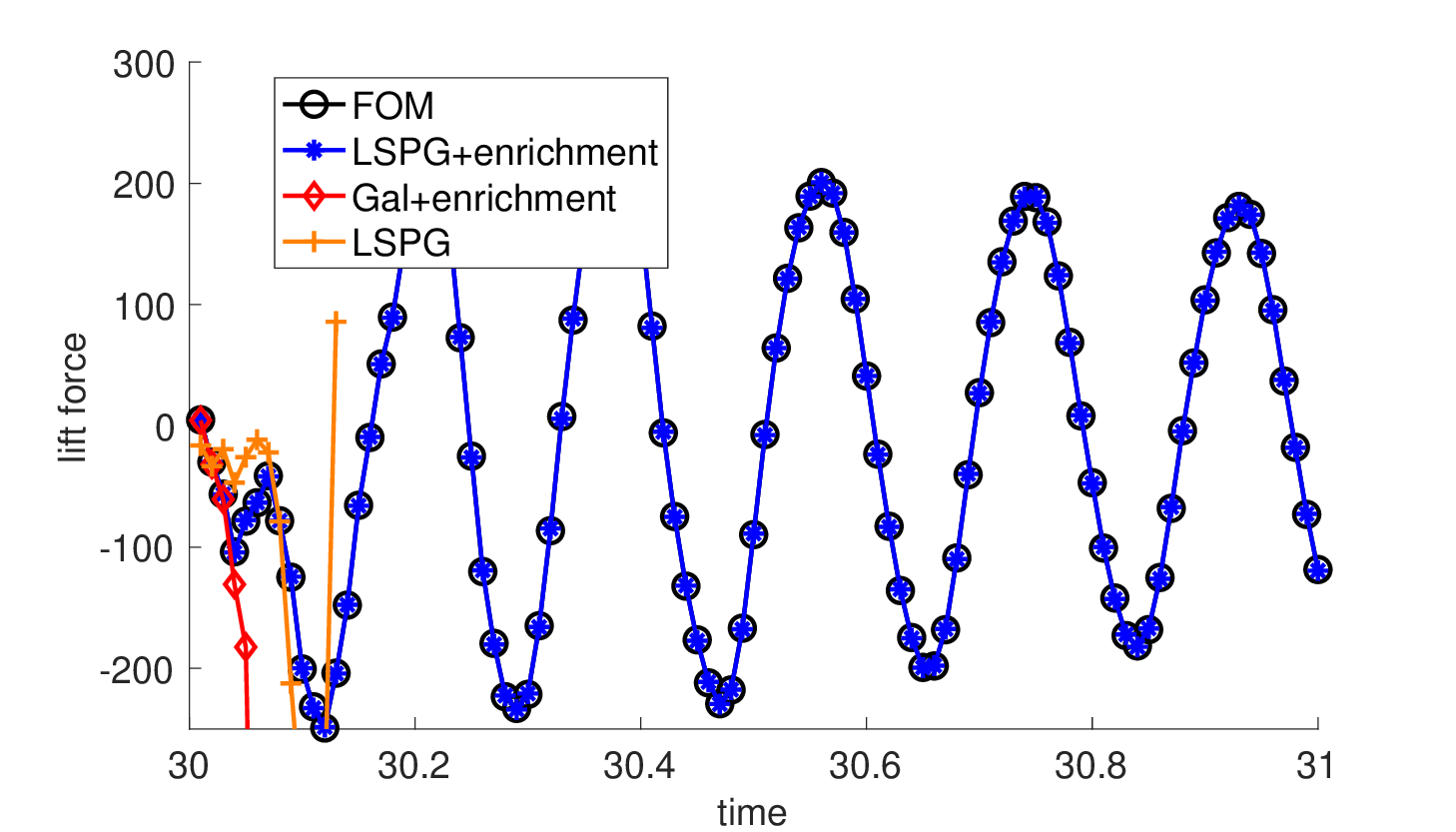}}
\caption{Turek problem; importance of enrichment and Petrov-Galerkin projection for the fluid subproblem. 
(a)-(b) estimates of the  drag and the lift forces (${\rm tol}_{\rm pod}=10^{-6}$).}
\label{fig:turek_MOR3}
\end{figure}

Figure \ref{fig:turek_MOR3} shows
the importance of the enrichment strategy and Petrov-Galerkin projection.
In more detail, we report the behavior of the drag and the lift forces predicted by the HF model  (FOM), the LSPG ROM with enrichment, the Galerkin ROM with enrichment, and the LSPG ROM without enrichment --- in all numerical tests, we consider a Galerkin ROM for the structure. We consider the POD tolerance $\text{tol}_{\text{pod}}=10^{-6}$.
We notice that the Galerkin ROM for the fluid is unstable for the very first time steps; similarly, the LSPG ROM shows instabilities at the very early stage of the simulation.
  
\subsubsection{Performance of the hybrid (ROM-FOM) solver}
\label{sec:res_hybrid_solver}

\begin{figure}[H]
\centering
\subfloat[\texttt{mesh0}]{\includegraphics[width=0.45\textwidth]{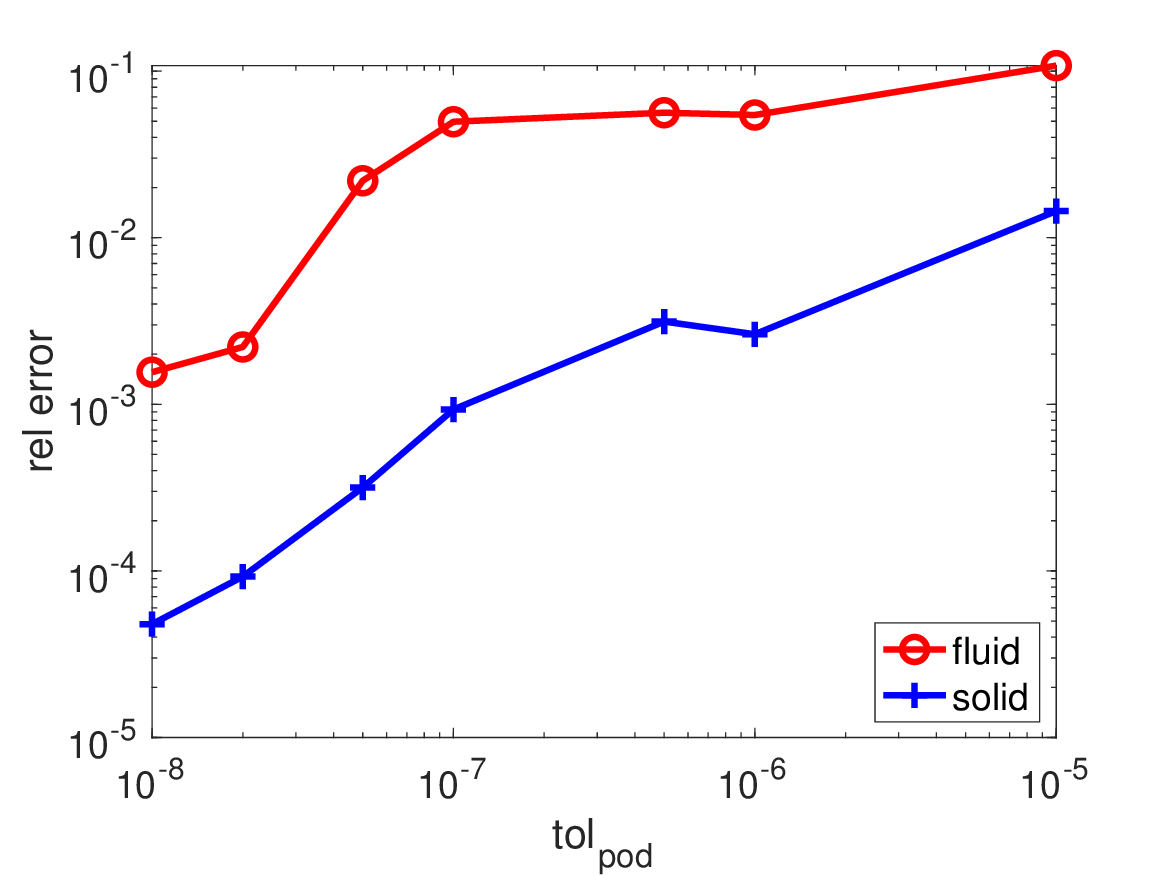}}
~~
\subfloat[\texttt{mesh1}]{\includegraphics[width=0.45\textwidth]{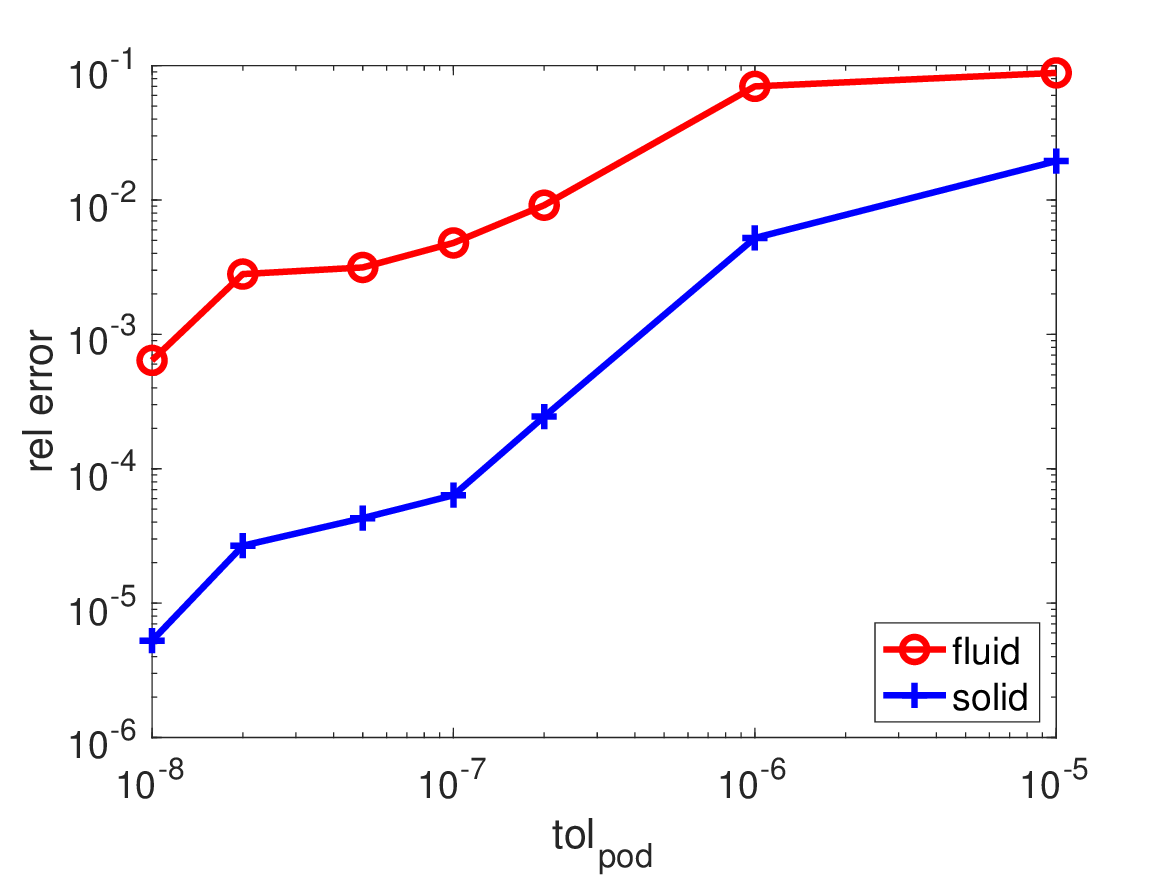}}

\caption{Turek problem; performance of the hybrid solver for several choices of the POD tolerance and for two meshes. }
\label{fig:turek_rel_error_hybrid}
\end{figure}

We report performance of the hybrid model where the FE basis is used for the fluid, and the POD bases are used for the solid state and the control. Simulations are performed over a  time interval of $10$ seconds --- that is, $(30 \mathrm{s}, 40\mathrm{s})$.
Figure \ref{fig:turek_rel_error_hybrid} shows the relative error (cf. equation \eqref{eq:rel_error}) associated with different POD tolerances for \texttt{mesh0} and 
\texttt{mesh1}.  We observe that the hybrid model produces accurate reconstructions over a longer time interval.

\begin{figure}[H]
\centering
\subfloat[]{\includegraphics[width=0.45\textwidth]{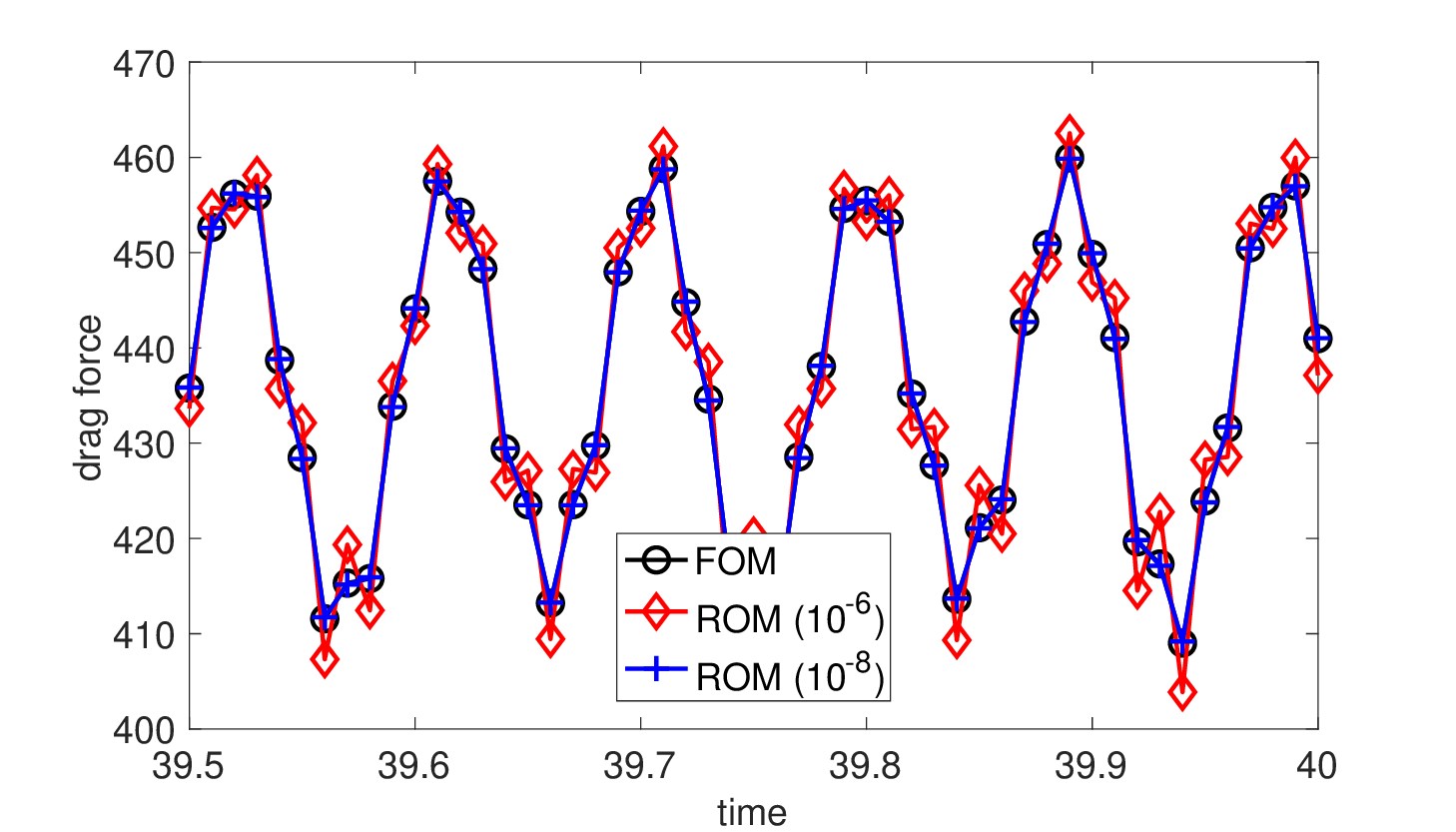}}
~~
\subfloat[]{\includegraphics[width=0.45\textwidth]{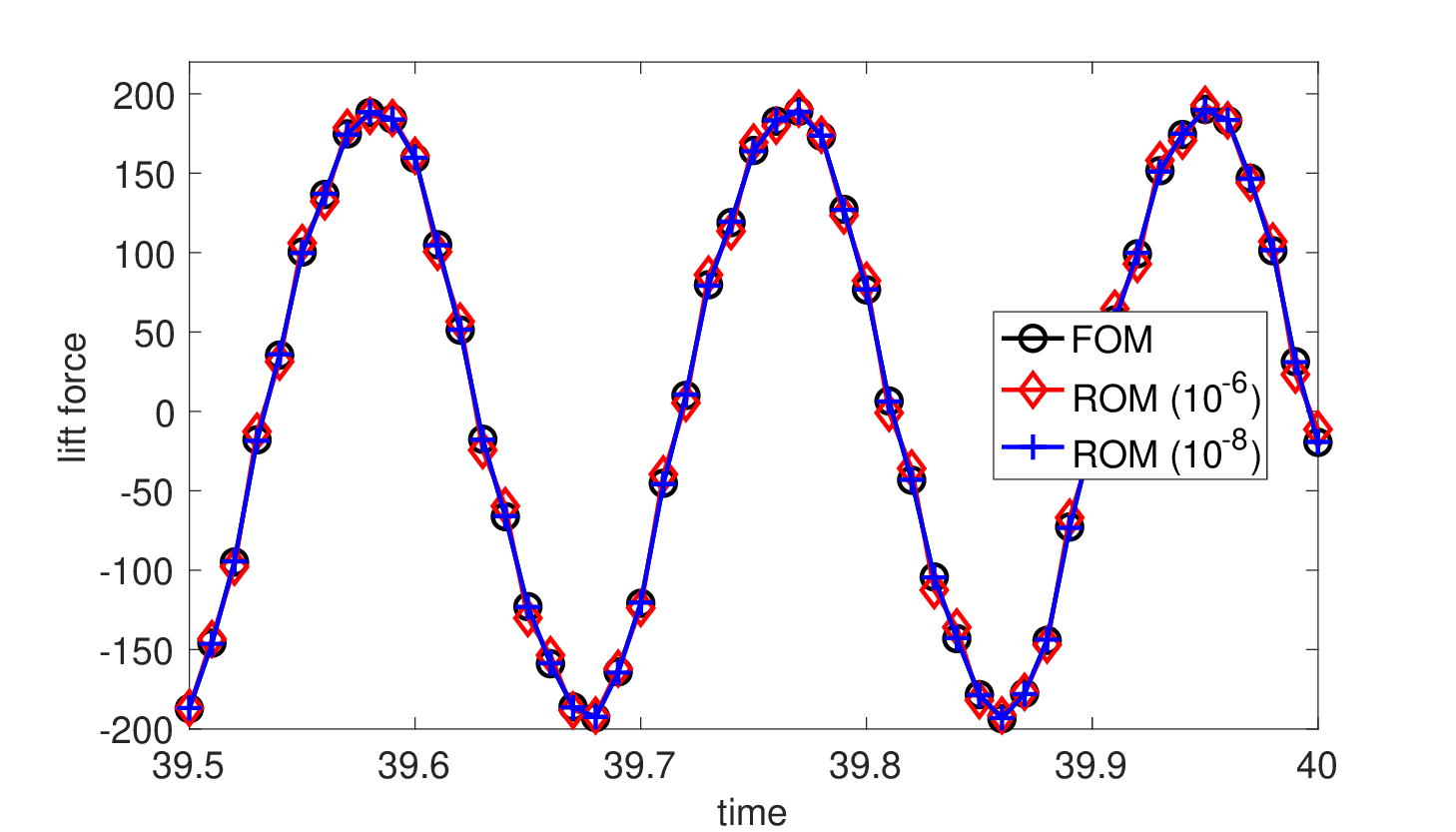}}

\subfloat[]{\includegraphics[width=0.45\textwidth]{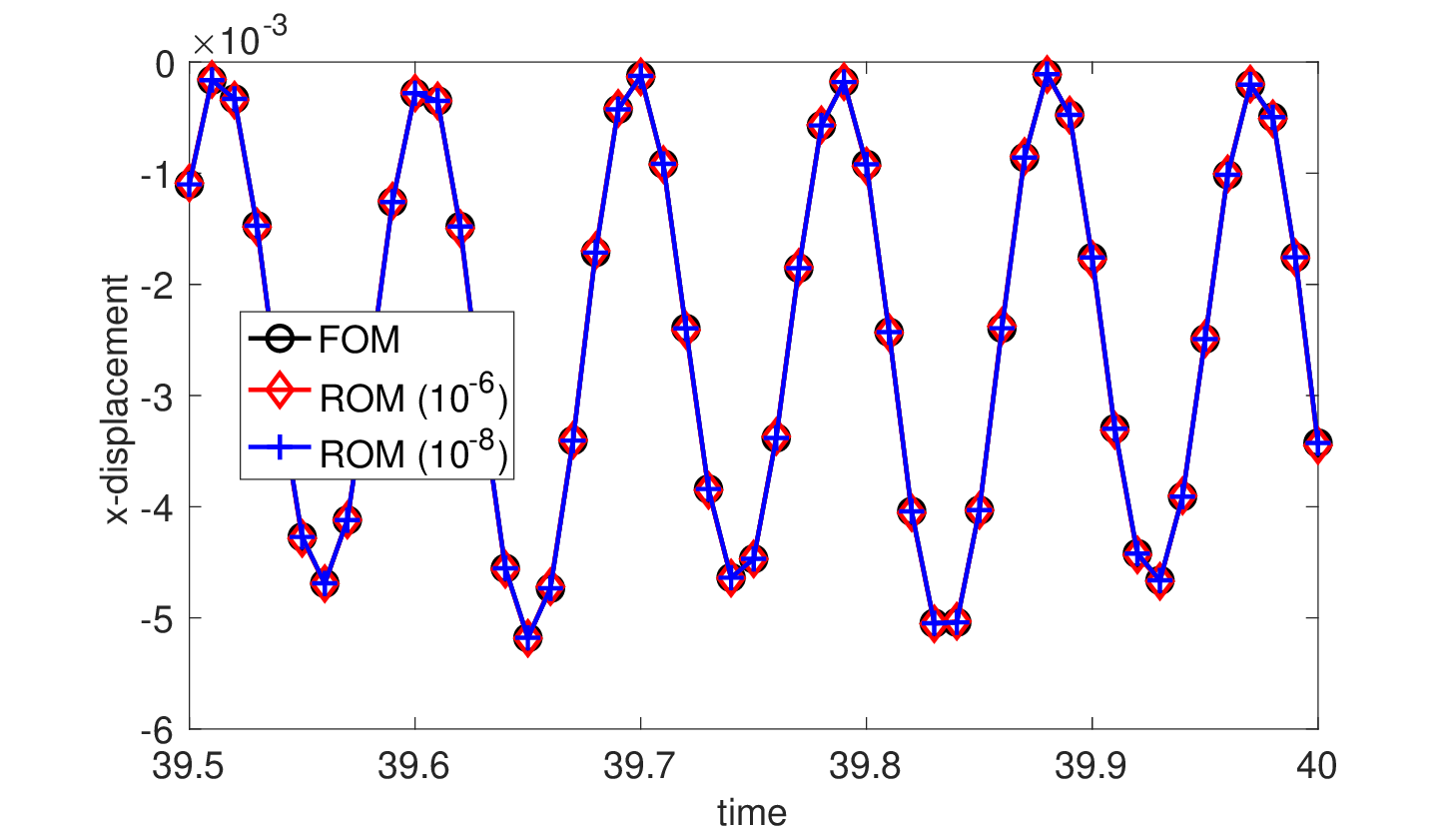}}
~~
\subfloat[]{\includegraphics[width=0.45\textwidth]{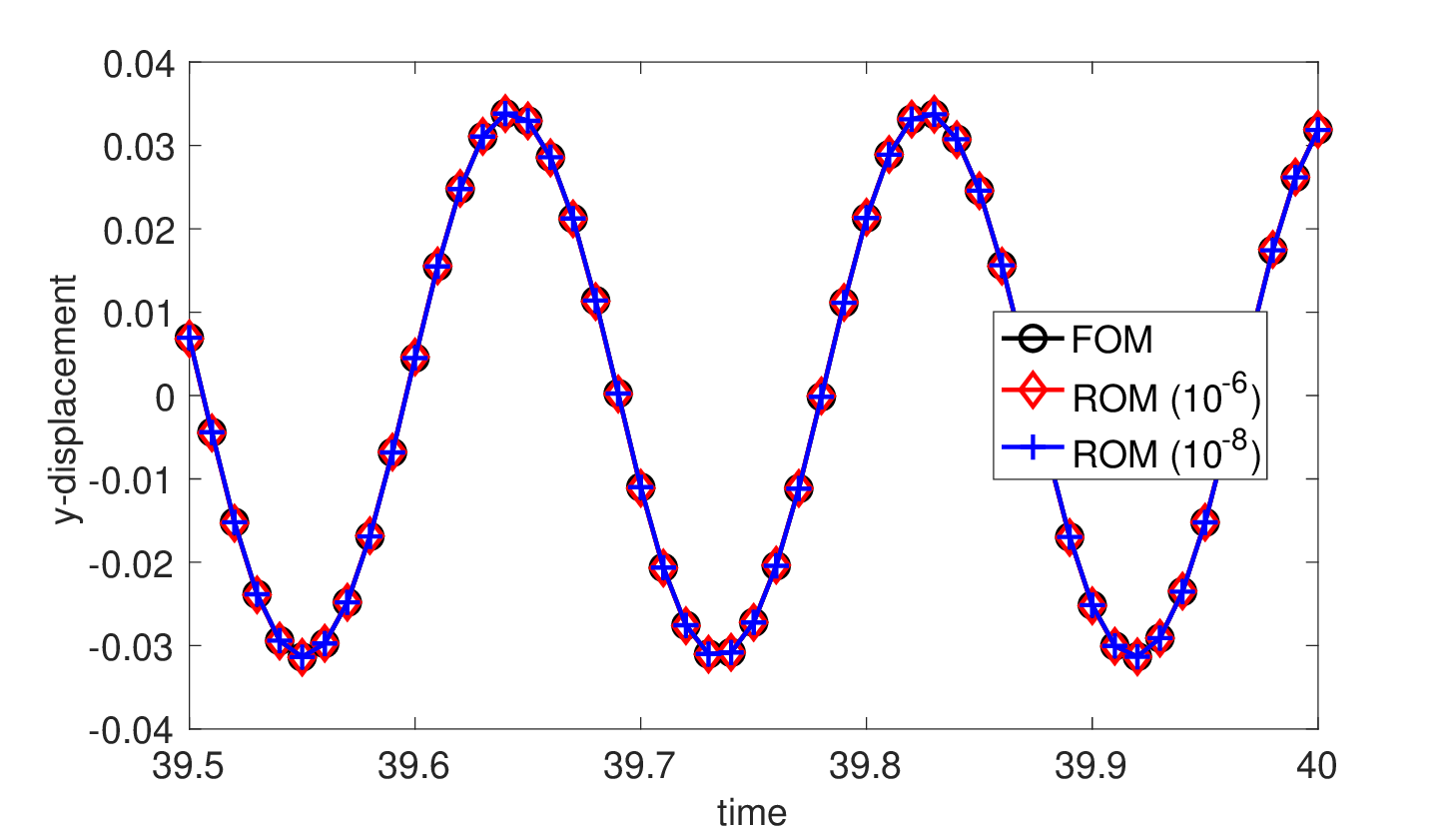}}

\caption{Turek problem; performance of the hybrid model.
(a)-(b) drag and lift forces.
(c)-(d) horizontal and vertical displacement of the control point $A$ (\texttt{mesh0}).}
\label{fig:turek_MOR_hybrid_mesh0}
\end{figure}

\begin{figure}[H]
\centering
\subfloat[]{\includegraphics[width=0.45\textwidth]{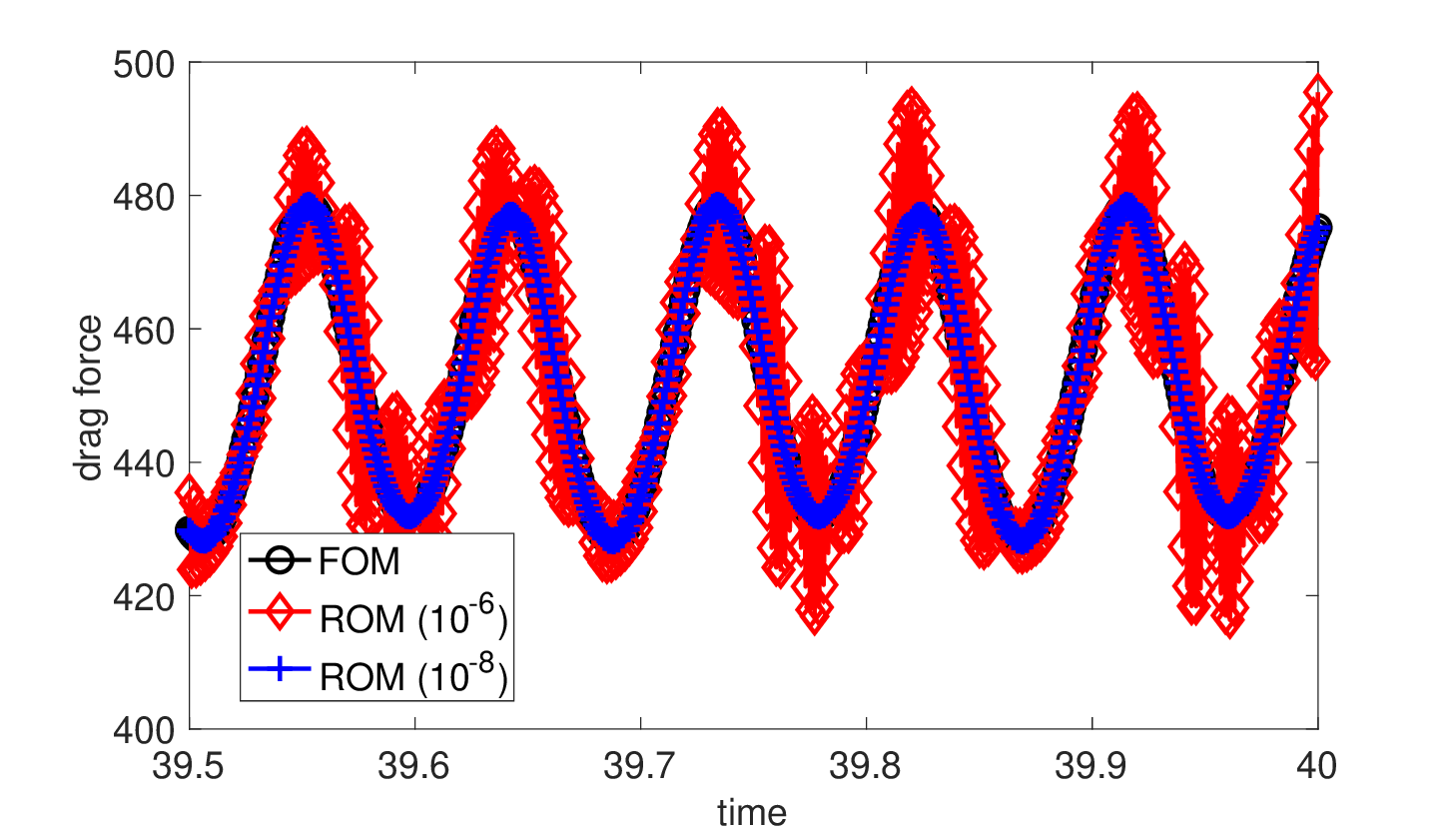}}
~~
\subfloat[]{\includegraphics[width=0.45\textwidth]{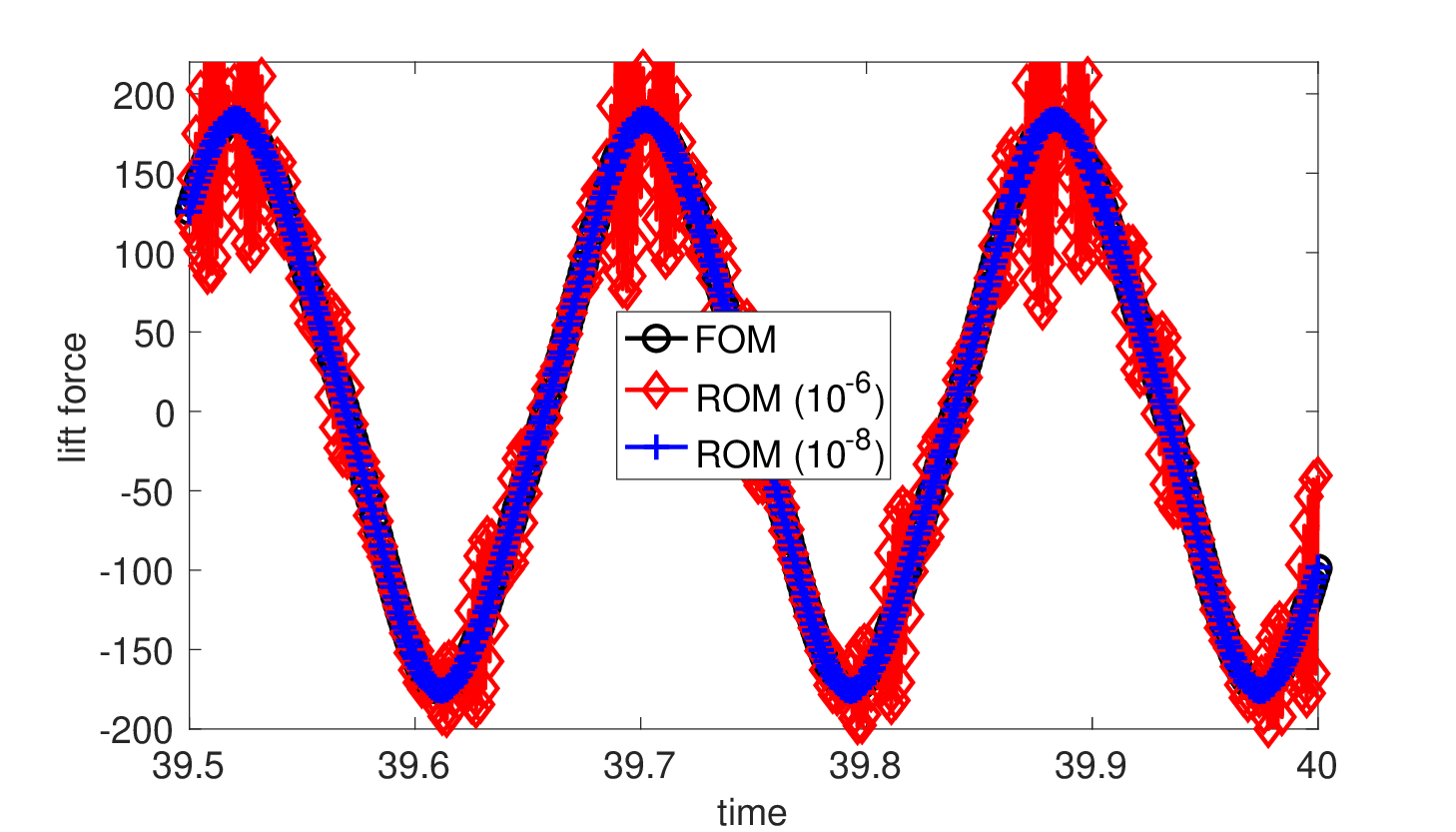}}

\subfloat[]{\includegraphics[width=0.45\textwidth]{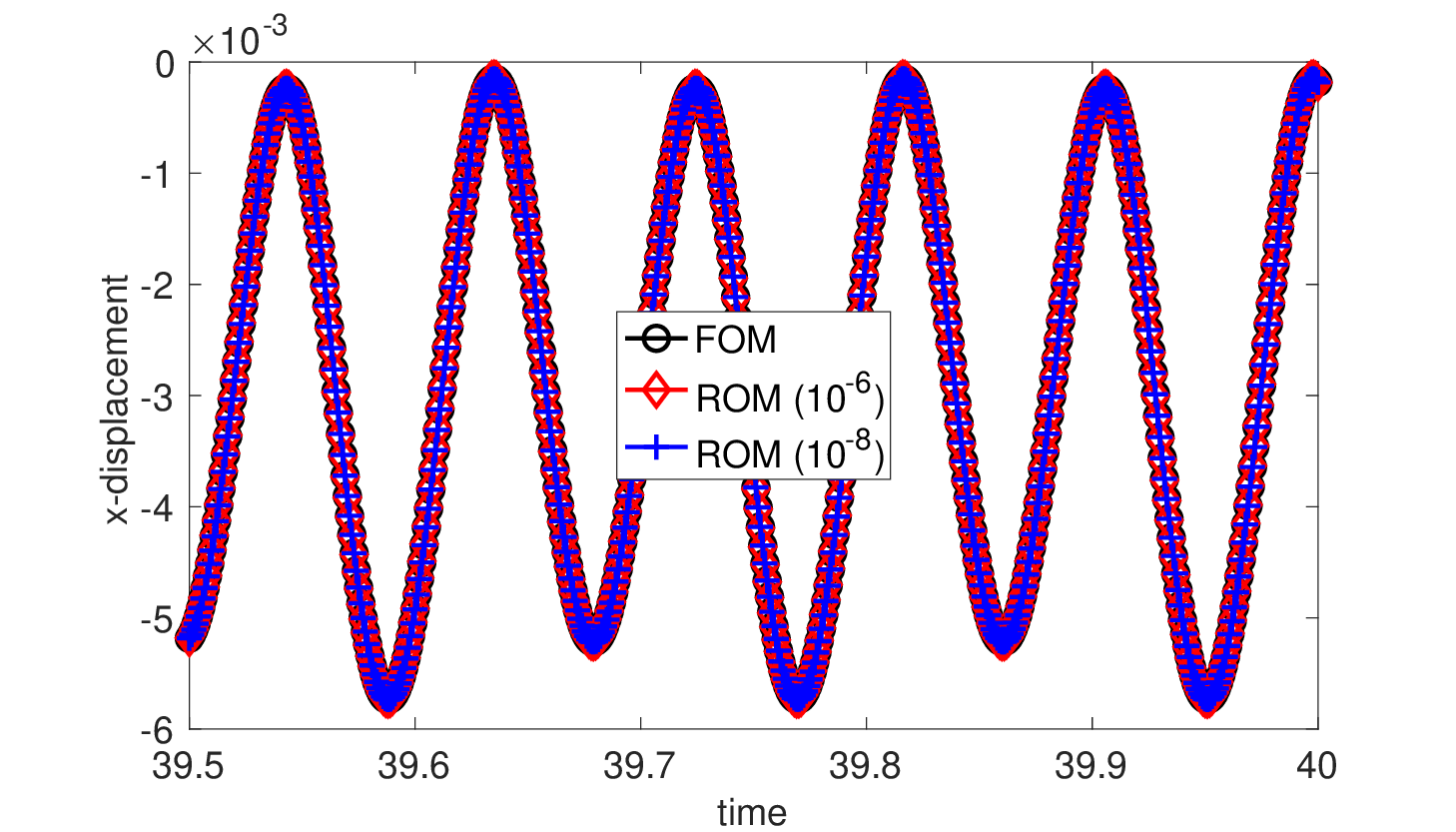}}
~~
\subfloat[]{\includegraphics[width=0.45\textwidth]{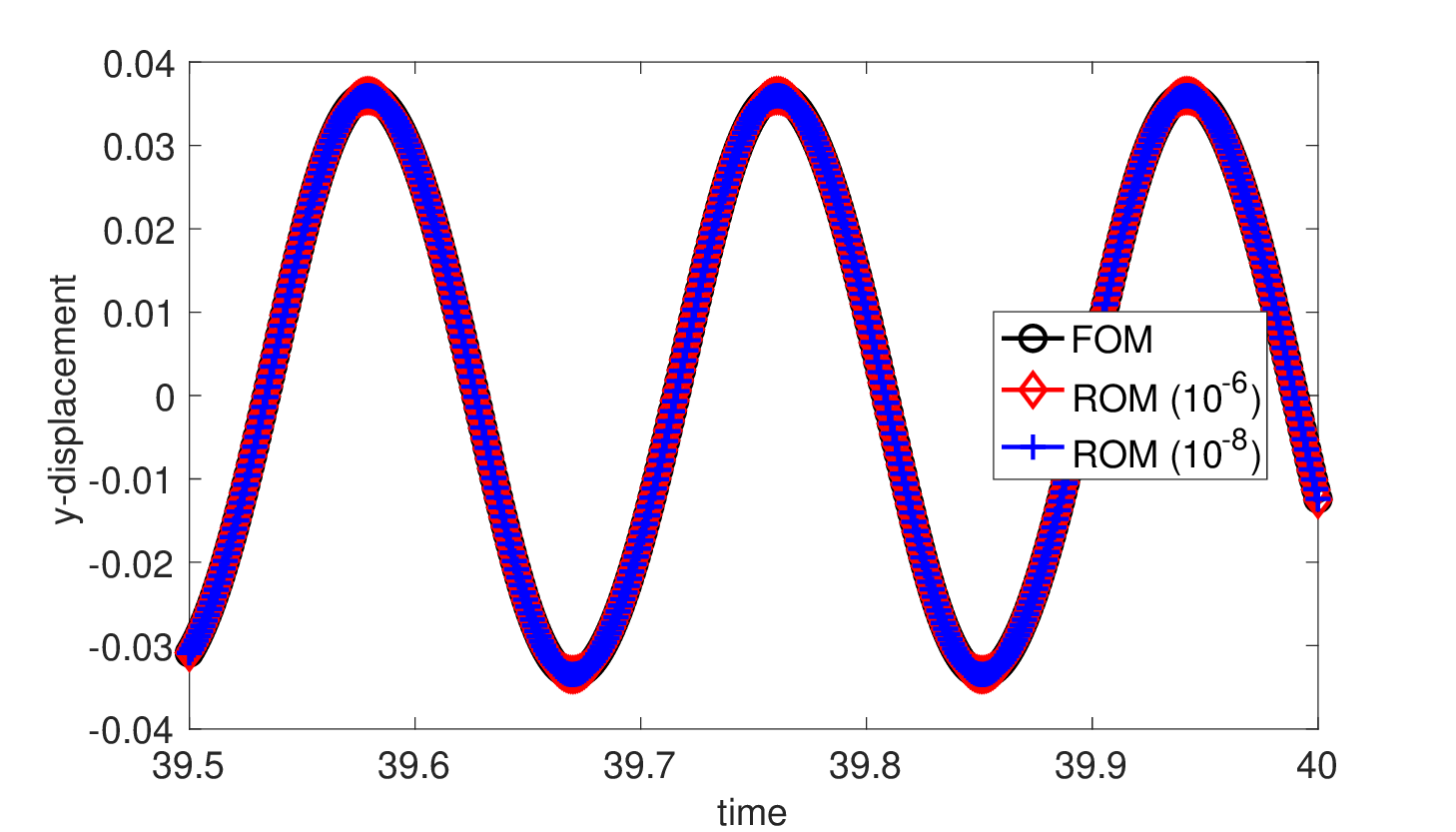}}
\caption{Turek problem; performance of the hybrid model.
(a)-(b) drag and lift forces.
(c)-(d) horizontal and vertical displacement of the control point $A$ (\texttt{mesh1}).}
\label{fig:turek_MOR_hybrid_mesh1}
\end{figure}

In Figure \ref{fig:turek_MOR_hybrid_mesh0}, we display the behavior of the  horizontal and the vertical displacements of the control point, along with the drag and lift forces, for the time interval $39.5$ to $40\,\mathrm{s}$ for the coarser discretization.
Figure \ref{fig:turek_MOR_hybrid_mesh1} replicates the results for the finer discretization \texttt{mesh1}.
The results indicate good accuracy for the lift force and the $x$- and $y$-displacements with both choices of the POD tolerance and both meshes;
on the other hand, the prediction of    the drag and the lift forces obtained with the coarser ROM on the finer mesh  exhibits 
spurious oscillations (cf. Figures \ref{fig:turek_MOR_hybrid_mesh1} (a)-(b)). We emphasize that these instabilities arise only when the integration is performed over a sufficiently long duration (approximately $7$ seconds).  

\section{Conclusions}
\label{sec:conclusions}
We developed an optimization-based method for  fluid structure interaction problems; 
our formulation, which relies on the ALE formalism, 
reads as a constrained optimization problem with equality constraints.
We introduced a control variable, which 
  corresponds to the normal flux at the fluid-structure interface,  to enable the effective solution to the optimization problem via static condensation of the local degrees of freedom.
We discussed an enrichment strategy of the local reduced spaces that ensures algebraic stability of the coupled problem.

We presented extensive numerical investigations to illustrate the performance of the method.
In more detail, 
the results of section \ref{sec:astorino_grandmont}   show the convergence rate of the FE solver for a problem with exact solution
(cf. Table \ref{tab:astorino_table});
the results of section \ref{sec:elastic_beam} show the 
benefit of SQP over the  standard Dirichlet-Neumann iterative method
for problems with control space of moderate size (cf. Table \ref{tab:vbeam_table});
finally, 
the results of section \ref{sec:turek} 
show the ability of the formulation to cope with hybrid (ROM-FOM) discretizations.

Our   numerical results  motivate further investigations to overcome the limitations of our method shown in  section 
\ref{sec:turek}.
First, 
we notice that  
the ROM might be prone to instabilities  for long-time integration and coarse
approximations (i.e., large POD tolerances).
To address this issue, we plan to devise 
new energy-conserving 
time integration schemes
  that satisfy  a 
  discrete counterpart of the energy balance  \eqref{eq:stability_properties_continuous};
  furthermore, following 
\cite{deparis2012stabilized,lovgren2006reduced}, we plan to resort to the Piola transformation to ensure the satisfaction of a discrete
incompressibility constraint. 
Second, the slow decay of the POD eigenvalues  for the  fluid subproblem (cf. Figure \ref{fig:turek_MOR1}) motivates 
the use of nonlinear approximations: we plan to apply registration-based methods \cite{taddei2020registration} that naturally fit in the ALE framework.
Third, we plan to develop space-time MOR formulations to ensure larger speedups for long-time integration
\cite{yano2014space} and enable more effective application of registration techniques \cite{taddei2021space}.
In this respect,  we remark that  space-time formulations of FSI problems are new in the MOR framework but have been previously considered  at the full-order level by Tezduyar  and collaborators
\cite{tezduyar1992new,tezduyar1992new2}.

\appendix

\section{Proof of Lemma \ref{th:stability_properties_continuous}}
\label{sec:proofs}
\begin{proof} We observe that:
$$
 \int_{\widetilde{\Omega}_{\rm s}} 
\rho_{\rm s} 
\frac{\partial^2 d_{\rm s}}{\partial t^2} \cdot 
\frac{\partial d_{\rm s}}{\partial t}
 \, dx  
= \frac{d}{dt} \int_{\widetilde{\Omega}_{\rm s}} 
\frac{\rho_{\rm s} }{2}
 \left|
\frac{\partial d_{\rm s}}{\partial t} 
  \right|^2 \, dx,
$$
and, recalling that 
$\sigma_{\rm s}=\frac{\partial W}{\partial F_{\rm s}}(F_{\rm s})$ and
$\frac{\partial F_{\rm s}}{\partial t} =  \nabla \big( \frac{\partial d_{\rm s}}{\partial t} \big)$, we find
$$
\int_{\widetilde{\Omega}_{\rm s}}
{\sigma}_{\rm s} \, : \,
\nabla \left( \frac{\partial d_{\rm s}}{\partial t} \right) \, dx
=
\int_{\widetilde{\Omega}_{\rm s}}
 \frac{\partial W}{\partial F_{\rm s}}     \, : \,
 \frac{\partial F_{\rm s}}{\partial t}  \, dx
 =
\int_{\widetilde{\Omega}_{\rm s}}
 \frac{\partial W}{\partial t}     \, dx
 =
 \frac{d}{dt} \left(
\int_{\widetilde{\Omega}_{\rm s}} \, W \, dx
 \right).
$$
Therefore, if we multiply 
\eqref{eq:weak_formulation_semidiscrete}$_3$ by $w=\frac{\partial d_{\rm s}}{\partial t}$, we find
\begin{equation}
\label{eq:proof_semidiscrete_1}
 \frac{d}{dt} \left(
\int_{\widetilde{\Omega}_{\rm s}}
\left( 
\frac{\rho_{\rm s} }{2}
 \Big|
\frac{\partial d_{\rm s}}{\partial t} 
 \Big|^2 
 \, + \, W \right) \, dx
 \right) \, + \, \int_{\widetilde{\Gamma}} g\cdot \frac{\partial d_{\rm s}}{\partial t}  \, dx = 0.
\end{equation}

For the fluid subproblem, exploiting the ALE transport formula \cite[Prop. 3.7]{formaggia2010cardiovascular} for the quantity $\rho_{\rm f}  \vert{u}_{\rm f} \vert^2$,
$$
\frac{d}{dt}\int_{\Omega_{\rm f}} \rho_{\rm f} \vert{u}_{\rm f}  \vert^2\,dx = 
\int_{\Omega_{\rm f}} \left(\rho_{\rm f} \frac{\partial \vert{u}_{\rm f} \vert^2}{\partial t}
\Big|_{\Phi_{\rm f}}
 + \rho_{\rm f}\vert{u}_{\rm f} \vert^2 \nabla\cdot \omega_{\rm f} \right) \,dx,
$$
we find that
\begin{equation}
\label{eq:proof_semidiscrete_2}
 \int_{\Omega_{\rm f}} 
\rho_{\rm f} \frac{\partial {u}_{\rm f}}{\partial t} \Big|_{\Phi_{\rm f}} \cdot {u}_{\rm f} \,dx
 \, = \,
  \int_{\Omega_{\rm f}} 
\frac{\rho_{\rm f}}{2}  \frac{\partial \big| {u}_{\rm f}  \big|^2 }{\partial t} \Big|_{\Phi_{\rm f}}   \,dx
 \, = \,
\frac{d}{dt}\int_{\Omega_{\rm f}} \frac{\rho_{\rm f}}{2} \vert{u}_{\rm f} \vert^2\,dx - \int_{\Omega_{\rm f}} \frac{\rho_{\rm f}}{2} \vert{u}_{\rm f} \vert^2 \nabla\cdot \omega_{\rm f}\,dx.  
\end{equation}
On the other hand, we have
$$
\begin{array}{rl}
\displaystyle{\int_{\Omega_{\rm f}} 
\rho_{\rm f} (u_{\rm f} - \omega_{\rm f}) \cdot \nabla u_{\rm f} \cdot u_{\rm f} \, dx
=} &
\displaystyle{\int_{\Omega_{\rm f}} 
\frac{\rho_{\rm f}}{2}
 (u_{\rm f} - \omega_{\rm f}) \cdot \nabla \big| u_{\rm f}\big|^2 \, dx}\\[3mm]
 = &
\displaystyle{
\int_{\partial \Omega_{\rm f}} 
\frac{\rho_{\rm f}}{2}
 (u_{\rm f} - \omega_{\rm f}) \cdot n_{\rm f} \big| u_{\rm f}\big|^2 \, dx 
\, - \,
\int_{  \Omega_{\rm f}} 
\frac{\rho_{\rm f}}{2}
 \big| u_{\rm f}\big|^2  \nabla \cdot u_{\rm f} \, dx
 +
 \int_{  \Omega_{\rm f}} 
\frac{\rho_{\rm f}}{2}
 \big| u_{\rm f}\big|^2  \nabla \cdot \omega_{\rm f} \, dx.}\\
\end{array}
$$
Recalling that $u_{\rm f}=0$ on $\partial \Omega_{\rm f} \setminus \Gamma$ and  $u_{\rm f}=\omega_{\rm f}$ on $\Gamma$, the boundary term vanishes and we obtain
\begin{equation}
\label{eq:proof_semidiscrete_3}
\int_{\Omega_{\rm f}} 
\rho_{\rm f} (u_{\rm f} - \omega_{\rm f}) \cdot \nabla u_{\rm f} \cdot u_{\rm f} \, dx
=
- \,
\int_{  \Omega_{\rm f}} 
\frac{\rho_{\rm f}}{2}
 \big| u_{\rm f}\big|^2  \nabla \cdot u_{\rm f} \, dx
 +
 \int_{  \Omega_{\rm f}} 
\frac{\rho_{\rm f}}{2}
 \big| u_{\rm f}\big|^2  \nabla \cdot \omega_{\rm f} \, dx.
\end{equation}

If we substitute \eqref{eq:proof_semidiscrete_2} and \eqref{eq:proof_semidiscrete_3} in \eqref{eq:weak_formulation_semidiscrete}$_3$ with $v = u_{\rm f}$, 
we find
\begin{equation}
    \label{eq:proof_semidiscrete_4}
\frac{d}{dt}\int_{\Omega_{\rm f}(t)} \frac{\rho_{\rm f}}{2} \vert{u}_{\rm f} \vert^2\,dx 
+\int_{\Omega_{\rm f}} 2\mu_{\rm f} \varepsilon_{\rm f}(u_{\rm f}) : \varepsilon_{\rm f}(u_{\rm f}) \, dx
- \,
\int_{\widetilde{\Gamma}} g\cdot u_{\rm f}\circ \Phi_{\rm f}
 \, dx
= 0.
\end{equation}
Since $u_{\rm f}\circ \Phi_{\rm f}
= \frac{\partial d_{\rm s}}{\partial t}$
on $\widetilde{\Gamma}$, if we sum \eqref{eq:proof_semidiscrete_1} and 
\eqref{eq:proof_semidiscrete_4}, we obtain the desired result.
\end{proof}

\section{Algebraic formulation of the FSI problem}
\label{sec:algebraic_formula}
We provide explicit expressions the terms in
\eqref{eq:weak_formulation_discrete_algebraic}.
In view of the discussion, we introduce several sets of indices.
We denote by   $\texttt{I}_{\rm f,\Gamma} =   \{ i_j^{\rm f} \}_{j=1}^{N_\Gamma} \subset \{1,\ldots,N_{\rm u}\}$  (resp., 
$\texttt{I}_{\rm s,\Gamma} =   \{ i_j^{\rm s} \}_{j=1}^{N_\Gamma} \subset \{1,\ldots,N_{\rm s}\}$) 
the  indices of the nodes  of the fluid (resp., solid) mesh on 
 the fluid-structure interface; 
 we denote by   $
\texttt{I}_{\rm f,dir} =  
 \{ i_j^{\rm f,dir} \}_{j=1}^{N_{\rm f,dir}} \subset \{1,\ldots,N_{\rm u}\}$  and
$\texttt{I}_{\rm s,dir} =  \{ i_j^{\rm s,dir} \}_{j=1}^{N_{\rm s,dir}} \subset \{1,\ldots,N_{\rm s}\}$  
 the Dirichlet nodes for the fluid and the solid domain, respectively;
 we denote by $
\texttt{I}_{\rm f,in} =  
 \{ i_j^{\rm f,in} \}_{j=1}^{N_{\rm f,in}} \subset \{1,\ldots,N_{\rm u}\}$  and
$\texttt{I}_{\rm s,in} =  
\{ i_j^{\rm s,in} \}_{j=1}^{N_{\rm s,in}} \subset \{1,\ldots,N_{\rm s}\}$  
 the complements of $\texttt{I}_{\rm f,dir}$ and  $\texttt{I}_{\rm s,dir}$, respectively.
 In the remainder, we use notation 
 $$
 \left\{
\begin{array}{l}
\displaystyle{
\mathbf{P}_{\rm f} \in \{0,1 \}^{N_{\Gamma}\times N_{\rm u}},
 \quad
 \mathbf{P}_{\rm s} \in \{0,1 \}^{N_{\Gamma}\times N_{\rm s}},
 \quad
  \mathbf{P}_{\rm f,dir} \in \{0,1 \}^{N_{\rm f,dir}\times N_{\rm u}},
 \quad
  \mathbf{P}_{\rm s,dir} \in \{0,1 \}^{N_{\rm s,dir}\times N_{\rm s}},
} \\
\displaystyle{
  \mathbf{P}_{\rm f,in} \in \{0,1 \}^{N_{\rm f,in}\times N_{\rm u}},
 \quad
  \mathbf{P}_{\rm s,in} \in \{0,1 \}^{N_{\rm s,in}\times N_{\rm s}},
}
\end{array}
 \right.
$$ 
  to refer to the mask matrices associated with $\texttt{I}_{\rm f,\Gamma}, \texttt{I}_{\rm s,\Gamma},
\texttt{I}_{\rm f,dir}, \texttt{I}_{\rm s,dir},  
  \texttt{I}_{\rm f,in}, \texttt{I}_{\rm s,in}$, respectively.

We define the FE fields $\{ w_j  \}_{j=1}^{N_\Gamma} \subset \widetilde{\mathcal{V}}_{\rm f}$  such that
$w_j$ is the FE solution to \eqref{eq:fsi_mesh_pseudo_elasticity} for $d_{\rm s} = \varphi_{i_j^{\rm f}}^{\rm f}$ for $j=1,\ldots, N_\Gamma$. Exploiting the linearity of \eqref{eq:fsi_mesh_pseudo_elasticity} and the first equation in 
  \eqref{eq:weak_formulation_discrete}$_4$, we find that
  $$
   \boldsymbol{\Phi}_{\rm f}
   (\mathbf{d} )  = \boldsymbol{\texttt{id}} + \mathbf{W} \mathbf{P}_{\rm s} \mathbf{d},
   \quad
   {\rm where} \;\; \mathbf{W} = [\mathbf{w}_1,\ldots,\mathbf{w}_{N_\Gamma}].
  $$
 Given the displacement $d\in \mathcal{V}_{\rm s}$, we define $\Omega_{\rm f}(d) = \Phi_{\rm f}(\widetilde{\Omega}_{\rm f}, d)$.

The time integration formulas \eqref{eq:newmark} and \eqref{eq:BDF} can be expressed as 
 $\texttt{D}_{{\rm s},\Delta t} \mathbf{d}_{\rm s}^{(k+1)} = \alpha_{\rm s,v}^{(k+1)} \mathbf{d}_{\rm s}^{(k+1)} + \mathbf{b}_{\rm s,v}^{(k+1)}$,
  $\texttt{D}_{{\rm s},\Delta t}^2 \mathbf{d}_{\rm s}^{(k+1)} = \alpha_{\rm s,a}^{(k+1)} \mathbf{d}_{\rm s}^{(k+1)} + \mathbf{b}_{\rm s,a}^{(k+1)}$
 for proper choices of $ \alpha_{\rm s,v}^{(k+1)},
\alpha_{\rm s,a}^{(k+1)} 
 \in \mathbb{R}$, $\mathbf{b}_{\rm s,v}^{(k+1)},
\mathbf{b}_{\rm s,a}^{(k+1)} 
 \in \mathbb{R}^{N_{\rm s}}$, and
 $\texttt{D}_{{\rm f},\Delta t} \mathbf{u}_{\rm f}^{(k+1)} = \alpha_{\rm f}^{(k+1)} \mathbf{u}_{\rm f}^{(k+1)} + \mathbf{b}_{\rm f}^{(k+1)}$. Exploiting these definitions, we can define the ALE velocity for the candidate displacement $d$ as 
$$
\boldsymbol{\omega}_{\rm f}^{(k+1)}  = 
  \mathbf{W} \mathbf{P}_{\rm s} 
\left(  
\alpha_{\rm s,v}^{(k+1)} \mathbf{d} + \mathbf{b}_{\rm s,v}^{(k+1)}
   \right).
$$

We have now the elements to introduce the discrete operators for the fluid and solid equations. First, we introduce the mass matrices
$$
\left( 
\mathbf{M}_{\rm f}(\mathbf{d})
\right)_{i,j}
=
\int_{\Omega_{\rm f}(d)} \phi_j^{\rm f}  \cdot  \phi_i^{\rm f} \, dx,
\quad
\left( 
\mathbf{M}_{\rm s} 
\right)_{i',j'}
=
\int_{\widetilde{\Omega}_{\rm s}} \varphi_{j'}^{\rm s}  \cdot  \varphi_{i'}^{\rm s} \, dx,
$$
for $i,j=1,\ldots,N_{\rm u}$ and 
 $i',j'=1,\ldots,N_{\rm s}$ --- note that in the first expression we omitted the dependence of the basis functions on $d$ to shorten notation.
Second, we introduce the residual in the fluid momentum equation
$$
\begin{array}{rl}
\displaystyle{
\left(
\widetilde{\mathbf{R}}_{\rm f}^{(k+1)}
\left(\mathbf{u}, \mathbf{p}, \mathbf{d}\right)
\right)_j
=} & 
\displaystyle{
\int_{\Omega_{\rm f}(d)} 
\left(
\sigma_{\rm f}(u,p) \, : \, \nabla \phi_j   \, + \,
\rho_{\rm f} (u - \omega_{\rm f}(d)) \cdot \nabla u \cdot \phi_j 
+
\frac{\rho_{\rm f}}{2}
(\nabla \cdot u) u\cdot \phi_j 
\right) \, dx
}
\\[3mm]
 & 
 \displaystyle{
-
\int_{\Gamma_{\rm f}^{\rm neu}(d)} g_{\rm f}^{\rm neu}(t^{(k+1)}) \cdot \phi_j  \, dx,
};
 \\
\end{array}
$$
and the matrix $\mathbf{B}_{\rm f} \in \mathbb{R}^{N_{\rm p}\times N_{\rm p}}$ associated with the fluid continuity equation
$$
\left(
\mathbf{B}_{\rm f}(\mathbf{d}
) \right)_{i,j} = - 
\int_{\Omega_{\rm f}(d)}
\, 
\nabla \cdot \phi_j^{\rm f} \cdot \theta_i^{\rm f} \, dx
\quad
i=1,\ldots,N_{\rm p}, \;\;
j=1,\ldots,N_{\rm u}.
$$
Third, we define the residual associated with  the solid equilibrium equation:
$$
\left(
\widetilde{\mathbf{R}}_{\rm s}^{(k+1)}
\left(\mathbf{d}\right)
\right)_j
=
\int_{\widetilde{\Omega}_{\rm s}} 
\left(
\sigma_{\rm s}(d): \nabla  \varphi_j^{\rm s}  \, - \,
f_{\rm s}^{(k+1)}   \cdot  \varphi_j^{\rm s}
\right) \, dx
-
\int_{\widetilde{\Gamma}_{\rm s}^{\rm neu}} g_{\rm s}^{{\rm neu},(k+1)} \cdot \varphi_j^{\rm s} \, dx,
\quad
j=1,\ldots,N_{\rm u}.
$$
Finally,  we introduce the matrices 
$\widetilde{\mathbf{E}}_{\rm f}\in \mathbb{R}^{N_{\rm u} \times N_{\Gamma}}$ and
$\widetilde{\mathbf{E}}_{\rm s}\in \mathbb{R}^{N_{\rm u} \times N_{\Gamma}}$ associated with the control:
$$
\left(
\widetilde{\mathbf{E}}_{\rm f}
\right)_{i,\ell}
=
\int_{\widetilde{\Gamma}} \varphi_{i_\ell^{\rm f}}^{\rm f} \cdot \varphi_i^{\rm f} \, dx,
\quad
\left(
\widetilde{\mathbf{E}}_{\rm s}
\right)_{j,\ell}
=
\int_{\widetilde{\Gamma}} \varphi_{i_\ell^{\rm f}}^{\rm f} \cdot \varphi_i^{\rm s} \, dx;
$$
note that $\widetilde{\mathbf{E}}_{\rm f}$ and $\widetilde{\mathbf{E}}_{\rm s}$ are time-independent and can hence be precomputed at the beginning of the iterative procedure.
In conclusion, 
we can state the fluid equations as in \eqref{eq:weak_formulation_discrete_algebraic}$_2$ with 
$$
\mathbf{R}_{\rm f}^{(k+1)}(\mathbf{u}, \mathbf{p}, \mathbf{d})   :=
\left[
\begin{array}{l}
\displaystyle{\mathbf{P}_{\rm f,in}^\top \left(
\mathbf{M}_{\rm f}(\mathbf{d}) \left(
\alpha_{\rm f}^{(k+1)} \mathbf{u} + \mathbf{b}_{\rm f}^{(k+1)}
\right)
\,+ \, 
\widetilde{\mathbf{R}}_{\rm f}^{(k+1)}
\left(\mathbf{u}, \mathbf{p}, \mathbf{d}\right)
\right)}
\\
\displaystyle{
\mathbf{B}_{\rm f}(\mathbf{d}) \mathbf{u}_{\rm f} 
} \\
\end{array}
\right],
\quad
\mathbf{E}_{\rm f}:=
\left[
\begin{array}{l}
\displaystyle{\mathbf{P}_{\rm f,in}^\top 
\widetilde{\mathbf{E}}_{\rm f} 
}
\\
\displaystyle{
0
} \\
\end{array}
\right],
$$
and the solid equations as  in \eqref{eq:weak_formulation_discrete_algebraic}$_3$ with 
$$
\mathbf{R}_{\rm s}^{(k+1)}(\mathbf{d}) 
=
\mathbf{P}_{\rm s,in}^\top \left(
\mathbf{M}_{\rm s} \left(
\alpha_{\rm s,a}^{(k+1)} \mathbf{d} + \mathbf{b}_{\rm s,a}^{(k+1)}
\right)
\,+ \, 
\widetilde{\mathbf{R}}_{\rm s}^{(k+1)}
\left(\mathbf{d}\right)\right),
\quad
\mathbf{E}_{\rm f}   :=
\mathbf{P}_{\rm s,in}^\top
\widetilde{\mathbf{E}}_{\rm s}. 
$$

\bibliographystyle{abbrv}
\bibliography{all_refs}

\end{document}